\documentclass[11pt,reqno]{amsart}
\usepackage{amsfonts}
\usepackage{eurosym}
\usepackage{amssymb}
\usepackage{amsmath}
 \usepackage{amsaddr} 
\usepackage{bm}
\usepackage{cite}
\usepackage{mathrsfs}
\usepackage{xcolor}
\usepackage[OT1]{fontenc}
\usepackage[left=2cm, right=2cm, top=3cm]{geometry}
\usepackage{hyperref}
\hypersetup{colorlinks=true, linkcolor=blue, citecolor=red}

\usepackage{times}

\usepackage[shortlabels]{enumitem}
\usepackage{enumitem}
\newcommand{\subscript}[2]{$#1 _ #2$}
\numberwithin{equation}{section}
\usepackage[square,numbers,sectionbib]{natbib}
\usepackage{varioref}

\usepackage{listings}
\usepackage{color}
\definecolor{dkgreen}{rgb}{0,0.6,0}
\definecolor{gray}{rgb}{0.5,0.5,0.5}
\definecolor{mauve}{rgb}{0.58,0,0.82}
\makeatletter
\@namedef{subjclassname@2020}{\textup{2020} Mathematics Subject Classification}
\makeatother

\allowdisplaybreaks
\newtheorem{theorem}{Theorem}[section]

\newtheorem{lemma}[theorem]{Lemma}

\theoremstyle{remark}
\newtheorem{remark}[theorem]{Remark}

\def \r {\rangle}

\newcommand{\numberset}{\mathbb}
\newcommand{\N}{\numberset{N}}
\newcommand{\R}{\numberset{R}}

\everymath{\displaystyle}

\def \no#1#2#3 {{\bf #1} (#3), #2.}
\def \eds#1#2#3 {#1, #2, #3.}

\begin{document}
	
	\title[	A phase-field system arising from
	multiscale modeling of thrombus
	biomechanics in blood vessels]
	{\large{A phase-field system arising from
			multiscale modeling of \\ thrombus
			biomechanics in blood vessels: \\ local well-posedness in dimension two}
 \\[2ex]
	{\normalfont \normalsize  \lowercase{\textit{In memoriam \uppercase{G}unduz \uppercase{C}aginalp, phase-field modeling pioneer and friend}}
	}
}

	\author[Grasselli \& Poiatti]{Maurizio Grasselli \& Andrea Poiatti}

	\address{Dipartimento di Matematica\\
		Politecnico di Milano\\
		Milano 20133, Italy\\
		\href{mailto:maurizio.grasselli@polimi.it}{maurizio.grasselli@polimi.it}\quad
		\href{mailto:andrea.poiatti@polimi.it}{andrea.poiatti@polimi.it}}
	
\maketitle	
	\begin{abstract} We consider a phase-field model which describes the interactions between the blood flow and the thrombus. The latter is supposed to be a viscoelastic material. The potential describing the cohesive energy of the mixture is assumed to be of Flory-Huggins type (i.e. logarithmic). This ensures the boundedness from below of the dissipation energy. In the two dimensional case, we prove the local (in time) existence and uniqueness of a strong solution, provided that the two viscosities of the pure fluid phases are close enough. We also show that the order parameter remains strictly separated from the pure phases if it is so at the initial time.
		\end{abstract}

	\section{Introduction}
	A thrombus is the final product of the blood coagulation step in hemostasis. In physiological conditions, blood clots form in the area of a vascular injury to prevent blood
	lremark and when the blood flow is high enough, the microthrombi do not adhere strongly
	to the walls of the vessels and can be destroyed by the fibrinolytic system \cite{Chapin, Colman, Formaggia}. In pathological cases the behavior is very different. Quoting \cite{Formaggia}, ``when the blood flow is low, the competition between the self-activation and
	inhibition of coagulation can favour thrombi formation and the thrombus
	can develop to block the vessel lumen inducing ischaemia of the tissues that
	it irrigates. The thrombi can also detach to form an embolus, which can
	lodge in a smaller vessel somewhere downstream from the site of thrombosis.
	These blockages are called infarcts and the results can be serious when they
	occur in the brain (stroke), the heart (heart attack) or the lungs (pulmonary
	embolism)''.
    Thrombus formation can be interpreted as a multiscale process (cf. \cite{Fedosov,Xu}), as well described in \cite{Karniadakis}: upon a vascular injury, freely flowing platelets in the blood become activated
	and form aggregates that cover the interested region, called thrombogenic area, within tens of seconds. Then, a
	series of biochemical reactions (see, e.g., \cite{Mack}) contribute to the formation of fibrin from the fibrinogen naturally present in blood. Then fibrin aggregates form a network, which facilitates the maturation of thrombus within minutes up to hours, say a much longer time than the first phase. If not
 degraded, thrombus can undergo a progressive remodeling with replacement of fibrin
	with collagen fibers (see, for instance, \cite{Fineschi} and the references therein), causing a rapid increase in stiffness, strength, and stability, as observed in
	\cite{em,Lee}.  Therefore, as pointed out in \cite{Karniadakis}, ``the structural constituents of a thrombus and the corresponding mechanical properties are vital to predict thrombus shape and deformation under various hemodynamic conditions and to evaluate the risk of thromboembolism and other
	pathological consequences''.
	The model proposed in \cite{Karniadakis} has the capacity of capturing the whole process of thrombus formation and its deformation with
	generation of emboli under dynamic shear flow conditions. After a first simulation of platelet aggregation, when the aggregates have become stable, the coarse-grained platelet distribution is converted into a continuum field for the estimation of the clot volume fraction and then the simulation of thrombus
	deformation within a blood flow is carried on using a phase-field model (see, for instance, the pioneering \cite{Cag}) , which is the objective of our analysis.
A lot of different models are used to describe interactions between fluids
and viscoelastic solids, such as the arbitrary Lagrangian Eulerian formulation (ALE) and
the immersed boundary method as well as level set, volume-of-fluid, and phase-field methods (see, e.g., \cite{Dong,Zeng} and their references). The main advantage of the phase-field modeling is the description of the system by means of the free energy. Different works have used this approach: in \cite{xu, Karniadakis} the authors model the interactions between flowing blood and thrombi, whereas in \cite{tierra} the interaction between blood and biofilms is studied. The phase-field method is characterized by the notion of diffuse interface (see, for instance, \cite{AMW,E} and references therein). This means
that the transition layer between the
phases has a finite size. Thus there is no tracking mechanism for the interface, but the phase state is
incorporated into the governing equations. The (diffuse) interface is associated with a smooth, but highly localized variation of the so-called phase-field variable $\phi$. In \cite{xu, Karniadakis} the domain $\Omega$ is a vascular lumen containing a thrombus, $\phi = 1$ represents the blood, $0 < \phi < 1$ represents a mixture of blood and thrombus, and $\phi=0$ represents thrombus only. Note that in our study we consider $\phi\in[-1,1]$, with $\phi=-1$ representing the "pure" thrombus. A harmless and usual redefinition of the order parameter as a (relative) concentration difference.

As typical in these cases, the equations are derived by minimizing the total energy $\mathcal{E}$ of the
blood-thrombus system, which is the sum of the kinetic energy, cohesive energy of the mixture, and the elastic energy of the fibrin network and platelets. The cohesive energy $\mathcal{E}_{coh}$ of the mixture is given by (see, e.g., \cite{CH}):
$$
\mathcal{E}_{coh}:= \int_\Omega\left(\frac{\sigma}{2}\vert\nabla\phi\vert^2+\Psi(\phi)\right)dx,
$$
with $\sigma>0$ is a capillary coefficient and $\Psi$ is the so-called Flory-Huggins potential given by (see, for instance, \cite{BTP} and references therein)
\begin{equation}
\label{FH}
\Psi(s)=\frac{\theta}{2}((1+s)\text{ln}(1+s)+(1-s)\text{ln}(1-s)) - \frac{\theta_0}{2} s^2, \quad \forall s\in (-1,1),
\end{equation}
where $\theta>0$ and $\theta_0>0$ are the absolute temperature and the critical temperature, respectively. We recall that $\theta<\theta_0$ in order to have a double-well structure for $\Psi$ which gives rise to phase separation.

The elastic energy $\mathcal{E}_{ela}$ arises from the induced extra stress tensor. To account for the viscoelastic rheology of the fibrin network in thrombus a Kelvin-Voigt type model is adopted, where a neo-Hookean relation
is used to describe the elastic behavior (Oldroyd-B model) and a Newtonian damper is used to describe the viscous behavior. In particular, the model follows the formulation presented in \cite{Lin}. The elastic energy thus reads
\begin{equation}
\label{eneela}
\mathcal{E}_{ela}=\int_\Omega \frac{\lambda(\phi)}{2}\text{tr}(\textbf{F}^T\textbf{F}-\textbf{I})dx,
\end{equation}
where $\lambda$ is the elastic relaxation time and $\textbf{F}$ is the deformation-gradient tensor. Notice that the dependence of $\lambda$ on $\phi$ is a critical issue. In \cite{Karniadakis} the function is given by $\lambda(s)=\lambda_e s,\ s\in \R$, with $\lambda_e>0$, but this choice seems to be mathematically hard to treat. Indeed, considering $\phi\in[0,1]$, the case of pure blood flow, i.e., $\phi=0$, would correspond to $\mathcal{E}_{ela}=0$, and this degeneracy in the energy does not allow to find appropriate estimates on a suitable approximating solution. Therefore, we assume $\lambda(s)\geq \lambda_*>0$ (see \eqref{lambda}). Note that this is physically reasonable, since the blood flow in the thrombogenic areas is often low, thus also the ``pure'' blood can be seen as presenting a viscoelastic behavior (see, e.g., \cite{Anand}). In an Eulerian reference frame, the deformation gradient tensor $\textbf{F}(x)=\partial \hat{x}/\partial X$, $\hat{x}$ and $X$ being the current and initial locations of a material point, satisfies the evolution equation (see, e.g., \cite{Lin}):
\begin{equation}
\label{defgrad}
\partial_t\textbf{F}+\textbf{u}\cdot \nabla\textbf{F}=\nabla\textbf{u}\cdot \textbf{F},
\end{equation}
where $\textbf{u}$ is the (volume averaged) fluid velocity and the dot stands for the matrix product. Note that (see again \cite{Lin}) if we assume $\text{div} \textbf{F}(0)=0$ then $\text{div}\textbf{F}(t)=0$ for all $t>0$. Thus, in the two dimensional case, this entails that $\textbf{F}$ can be expressed as $\textbf{F}=\nabla\times \psi$ and the corresponding transport equation becomes
\begin{equation}
\label{potential}
\partial_t \psi +\textbf{u}\cdot \nabla\psi=\textbf{0}.
\end{equation}
In conclusion, the last component of the total energy is the kinetic energy, which accounts for the mixture transport. This is given by
$$
\mathcal{E}_{kin}=\int_\Omega\frac{1}{2}\vert\textbf{u}\vert^2dx.
$$
Hence, the total energy of the system is
$$
\mathcal{E}=\mathcal{E}_{kin}+\mathcal{E}_{coh}+\mathcal{E}_{ela}.
$$
On the other hand, the combination of the macroscopic dissipation and microscopic dissipation of the system is defined by the following functional
\begin{equation}
\label{dissip}
\mathcal{D}=\int_\Omega\left[\nu(\phi) \vert D\textbf{u}\vert^2+\frac{\nu(\phi)}{2k(\phi)}(1-\phi)\vert\textbf{u}\vert^2+\vert\nabla\mu\vert^2\right]dx,
\end{equation}
where $D\textbf{u}$ is the symmetrized strain tensor, $\nu$ is the kinematic viscosity and $k$ is the permeability. The latter accounts for friction induced by the porous media flow inside the clot.
Concerning this point, we recall that ``a thrombus can be considered as
porous medium exhibiting viscoelastic behavior with a fibrous network contributing to both the viscous and the elastic properties''(\cite{Karniadakis}). Notice that both $\nu$ and $k$ depend on $\phi$, as introduced in \cite{Karniadakis}, since they are considered as spatially variable (clearly $k$ and $\nu$ must be strictly positive). Moreover, in two dimensions, the chemical potential $\mu$ is defined by
\begin{equation}
\label{chempot}
\mu:=\frac{\delta\left(\mathcal{E}_{coh}+\mathcal{E}_{ela}\right)}{\delta\phi}=-\sigma\Delta\phi+\Psi^\prime(\phi)+\frac{\lambda^\prime(\phi)}{2}\vert\nabla \psi\vert^2.
\end{equation}
Applying then Onsager's principle and force balancing, one obtains the following system of partial differential equations
	\begin{equation}
	\begin{cases}
	\partial_t \textbf{u}+(\textbf{u}\cdot\nabla)\textbf{u}+\nabla \pi-\text{div}(\nu(\phi) D\textbf{u})+\text{div}(\lambda(\phi)\nabla\psi^T\nabla\psi)\\+\text{div}(\nabla\phi\otimes\nabla\phi)+\nu(\phi)\frac{(1-\phi)\textbf{u}}{2k(\phi)}=\textbf{0},\\
	\text{div}\,\textbf{u}=0,\\
	\partial_t \psi +\textbf{u}\cdot \nabla\psi=\textbf{0},\\
	\partial_t \phi+\textbf{u}\cdot \nabla\phi-\Delta\mu=0\\
	\mu=-\Delta\phi+\Psi^\prime(\phi)+\frac{\lambda^\prime(\phi)}{2}\vert\nabla \psi\vert^2,\\
\end{cases}
	\label{syst}
	\end{equation}
in $\Omega\times(0,T)$, for a given $T>0$, endowed with the following initial and boundary conditions
\begin{equation}
	\begin{cases}
\textbf{u}(0)=\textbf{u}_0, \quad {\phi}(0)={\phi}_0, \quad {\psi}(0)={\psi}_0, \quad \text{ in }\Omega,\\
\textbf{u} = \textbf{0}, \quad \partial_\textbf{n} \phi= \partial_\textbf{n}\mu=0, \quad\text{ on } \partial\Omega\times (0,T).\\
	\end{cases}
	\label{IBC}
	\end{equation}
	Here $\pi$ denotes the pressure of the fluid mixture and some positive constants have been set equal to the unity.

Observe that, in absence of the terms related to $\psi$ and the permeability term, the above system is the well known model H (see \cite{GPV}, see also \cite{AMW} and references therein) which has been
studied by several authors also in the case of non-constant viscosity (see, for instance \cite{Abels, gal2,Giorginitemam} and references therein; see also \cite{ADG,AF,Giorgini} for more general models with unmatched densities).
Concerning the Oldroyd-B model and viscoelastic fluids, corresponding to the coupling without the order parameter $\phi$, we refer to \cite{Lin}, where the authors established local and global existence of classical solutions without using an artificial damping mechanism. Note that the global solution is obtained under a suitable smallness condition on the initial data (see \cite[Thm.2.3,][]{Lin}). Moreover, in \cite{Lin2} the authors show the local well-posedness of the initial-boundary value problem of the viscoelastic fluid system related to the Oldroyd model and they prove the global existence of regular solutions provided that the initial data are close enough to the equilibrium state.

The only work related to the above model is the recent contribution \cite{KTT}. The authors analyze the three-dimensional model. Therefore, in their case, equation \eqref{potential} is replaced by \eqref{defgrad} and,
as a consequence, the chemical potential \eqref{chempot} depends on $\textbf{F}$. They establish the local (in time) existence and uniqueness of a strong solution. Concerning $\Psi$, following \cite{xu,Karniadakis}, the authors take a double-well fourth-order polynomial approximation of the Flory-Huggins potential \eqref{FH}. Nevertheless, in that case, it is not longer guaranteed that the phase-field $\phi$ takes values in the physical range $[-1,1]$ (or in $[0,1]$),  due to the lack of comparison principles for the Cahn-Hilliard equation, namely, the equation which rules the evolution of $\phi$.
Consequently, the dissipation functional \eqref{dissip} is no longer bounded from below. Also, a smooth approximation of the potential allows the author to handle
a possibly degenerating elastic energy, taking $\lambda$ as in \cite{Karniadakis} (see \eqref{eneela}). However, in this case one cannot ensure the non-negativity of
the elastic energy. These considerations motivate our decision to investigate the case of logarithmic potential \eqref{FH}. However, this choice is a source of several
technical difficulties which seem hard to overcome in the three dimensional case so far. For instance, the strong solution found in \cite{KTT} is stronger than ours and is obtained through classical energy estimates performed in a full Galerkin scheme.
In particular, the authors differentiate the Navier-Stokes system and multiply by $\partial_t\textbf{u}$. In the present case we cannot proceed in that way and we need to test the Navier-Stokes system with suitable negative powers of the Stokes operator applied to $\partial_t\textbf{u}$. Moreover, here we need to use $L^p$ estimates and the construction of an approximating solution is far more complicated. All these technical difficulties arise from the choice of a singular potential.

Summing up, here we prove the existence of a local (in time) existence and uniqueness of a strong solution to problem \eqref{syst}-\eqref{IBC}, provided that the viscosities of the two pure components of the fluid are sufficiently close to one another. Moreover, we show that the order parameter $\phi$ stays uniformly away from the pure phases $\pm1$, having assumed that the initial datum $\phi_0$ is itself strictly separated (see Remark \ref{separated}).

The plan of the paper goes as follows.
In Section \ref{preliminary} we introduce the functional setup and the main assumptions. Section \ref{main2} contains the statement of our main result. Section \ref{tools} we present a number of auxiliary results which are necessary to perform the proof of the main theorem (see Section \ref{proof}). Finally, in Section \ref{proofs}, we prove the auxiliary results.

\section{Notation and basic assumptions}
\label{preliminary}
	\subsection{Notation and function spaces}
	\label{math}
	Let $\Omega$ be a bounded domain of $\R^2$ with a smooth boundary.
	For the velocity field we set, for $p\in(2,\infty)$,
	\begin{align*}
	&\textbf{H}_\sigma=\overline{\{\textbf{u}\in\ [C^\infty_0(\Omega)]^2:\ \text{div}\ \textbf{u}=0\}}^{[L^2(\Omega)]^2}\ \ \ \ \ 	\textbf{V}_\sigma=\overline{\{\textbf{u}\in\ [C^\infty_0(\Omega)]^2:\ \text{div}\ \textbf{u}=0\}}^{[H^1(\Omega)]^2}\\&
	\textbf{L}^p_\sigma=\overline{\{\textbf{u}\in\ [C^\infty_0(\Omega)]^2:\ \text{div}\ \textbf{u}=0\}}^{[L^p(\Omega)]^2}\quad\quad \textbf{W}_\sigma=[{H}^2(\Omega)]^2\cap\textbf{V}_\sigma.
	\end{align*}
    We also introduce 
	$\textbf{V}_\sigma^s=D(\textbf{A}^{s/2})$, for any $s\in \R$. Here $\textbf{A}$ is the usual Stokes operator. We observe that the classical $H^s$-norm in $\textbf{V}_\sigma^s$ is equivalent to the norm $\Vert\cdot\Vert_{\textbf{V}_\sigma^s}$ for $s>-1/2$ (see, e.g., \cite{salgado}).
    Then we set
	\begin{equation*}
	V=H^1(\Omega),\ \ \ \ V_2=\{v\in H^2(\Omega):\ \partial_\textbf{n}v=0\ \text{a.e. on }\partial\Omega\}.
	\end{equation*}
	We denote by $(\cdot, \cdot)$ the standard inner product in $\textbf{H}_\sigma$ (or in $H$) and by $\Vert   \cdot\Vert$ the induced norm. In $\textbf{V}_\sigma$ we define the inner product $(\textbf{u},\textbf{v})_{\textbf{V}_\sigma}=(\nabla\textbf{u},\nabla\textbf{v})$ whose induced norm is $\Vert   \textbf{v}\Vert   _{\textbf{V}_\sigma}=\Vert   \nabla\textbf{v}\Vert$. We then indicate by $(\cdot,\cdot)_V$ the canonical inner product in $V$, while $\Vert   \cdot\Vert_V$ stands for its induced norm. Moreover, we define $\Vert  \textbf{v}\Vert  ^2_{\textbf{W}_\sigma}:=(\textbf{Av},\textbf{Av})$.
    Recalling \cite[Appendix A,][]{Giorginitemam}, we introduce $V_0:=\{v\in V: \overline{v}=0\}$ and we denote its dual by $V_0'$. The restriction ${A}$ to $V_0$
    of $A_0: V \to V^\prime$, $<A_0u,v>=(\nabla u,\nabla v)$ for all $v\in V$, is an isomorphism from $V_0$ onto $V_0'$. Thus ${A}^{-1}$ stands for its inverse map
    and we set $\Vert  f\Vert  _{\sharp}:=\Vert  \nabla {A}^{-1}f\Vert$ which is a norm in $V_0'$ that is equivalent to the canonical one.
	\subsection{Main assumptions}
\label{Assumptions}
	On account of \eqref{FH}, we suppose more in general
	\begin{enumerate}[label=(\subscript{H}{{\arabic*}})]
		\item 	$\Psi(s)=F(s)-\frac{\alpha_0}{2}s^2$,
		where $F\in C([-1, 1]) \cap C^3 (-1,1)$ fulfills
		\begin{equation*}
		\lim_{s\to-1} F'(s)=-\infty, \quad \lim_{s\to1} F'(s)=+\infty,\quad F''(s)\geq\overline{\alpha}, \quad\forall\ s\in(-1,1),
		\end{equation*}
		assuming $\tilde{\alpha}=\alpha_0-\overline{\alpha}>0$ we have that $\Psi$ is double well and
		\begin{equation}
		\Psi''(s)\geq -\tilde{\alpha},\quad\forall s\in(-1,1)
        \label{quasiconv}
		\end{equation}
		We also extend $F(s)=+\infty$  for any $s\notin[-1, 1]$.
		Notice that the above assumptions imply that there
		exists $s_0\in(-1, 1)$ such that $F'(s_0 )=0$. Without lremark of generality, we can take $s_0 = 0$ and
		$F(s_0 ) = 0$. In particular, this entails that $F (s) \geq 0$ for any $s \in[-1, 1]$. Moreover we require that $F''$ is convex and
		\begin{equation}
		F''(s)\leq Ce^{C \vert F'(s) \vert }, \quad \forall s\in(-1,1)
		\label{growth}
		\end{equation}
		for some positive constant C. Also, we assume that there exists $\gamma\in(0,1)$ such that $F''$ is nondecreasing in $[1-\gamma,1)$ and nonincreasing in $(-1,-1+\gamma]$.
		\label{h3}
	\end{enumerate}
	Concerning $\lambda$, we assume that $\lambda\in C^{2}(\R)$ and convex. Moreover, we ask
    \begin{equation}
	0<\lambda_*\leq \lambda(s)\leq \lambda^*,\qquad\lambda_{1,*}\leq \lambda^\prime(s)\leq \lambda_1^*, \quad\forall s\in \R,
	\label{lambda}
	\end{equation}
	for some positive values $\lambda_*,\lambda^*,\lambda_{1,*},\lambda_1^*$.
    We assume
    \begin{align}
    \nu(s)=\nu_1\frac{1+s}{2}+\nu_2\frac{1-s}{2},\quad s\in[-1,1],
    \label{nu1}
    \end{align}
    which can be easily extended to the whole $\R$  in such a way to comply
    \begin{equation}
    0<\nu_*\leq \nu(s)\leq \nu^*,\qquad \forall s\in \R,
    \label{nuA}
    \end{equation}
    for some positive values $\nu_\ast, \nu^\ast$. Also, we suppose $k\in W^{1,\infty}(\R)$ and
        \begin{equation}
    0<k_*\leq k(s)\leq k^*,\qquad \forall s\in \R.
    \label{k}
    \end{equation}
	The assumptions on the initial data are:
	\begin{enumerate}[label=(\subscript{H}{{\arabic*}})]
		\item $\textbf{u}_0\in \textbf{V}_\sigma^{3/2+2\widetilde{\delta}}$ for some $\widetilde{\delta}\in (0,1/4]$; \label{Hbegin}
		\item\label{Hmedium}  $\psi_0\in [W^{2,\frac{4}{1+2\varepsilon}}(\Omega)]^2$ for some $\varepsilon\in (0,1/6)$;
		\item  $\phi_0\in H^2(\Omega)$ such that $\partial_\textbf{n}\phi_0=0$ on $\partial\Omega$ and $\mu_0:=-\Delta\phi_0+\Psi^\prime(\phi_0)\in V$, with $-1\leq \phi_0\leq 1$ for any $x\in\overline{\Omega}$
      and $\overline{\phi}_0\in(-1,1)$. \label{Hend} 
	\end{enumerate}

\begin{remark}
As already noticed, e.g., in \cite[Remark 3.9]{GP}, being $\mu_0\in V$, recalling \cite[Thm.A.2,][]{Giorginitemam}) and setting $f=\mu_0+\alpha_0\phi_0\in V$, we have $\phi_0\in W^{2,p}(\Omega)$ for every $p\geq2$ and $F'(\phi_0)\in L^p(\Omega)$ for every $p\geq2$. Moreover, from property (\ref{growth}) of $F$ and again by \cite[Thm.A.2,][]{Giorginitemam}, being $f\in V$, we deduce that $F''(\phi_0)\in L^p(\Omega)$ for every $p\geq2$. Since $\phi_0$ belongs to $V$ and $F'\in C^2((-1,1))$, we can apply the chain rule to obtain
	\begin{equation*}
	\nabla F'(\phi_0)=F''(\phi_0)\nabla\phi_0,
	\end{equation*}
	but then $\Vert \nabla F'(\phi_0)\Vert_{L^p(\Omega)} \leq \Vert F''(\phi_0) \Vert_{L^{2p}(\Omega)}\Vert \nabla\phi_0 \Vert_{L^{2p}(\Omega)}<\infty$, since $  \nabla\phi_0\in W^{1,q}(\Omega)$ for every $q\geq2$. Therefore we get $F'(\phi_0)\in W^{1,p}(\Omega)$ for every $p\geq2$, implying that $F'(\phi_0)\in L^\infty(\Omega)$ and thus we obtain that the initial field $\phi_0$ is necessarily strictly separated, namely, there exists $\tilde{\delta}>0$ such that
	\begin{equation*}
	\Vert \phi_0\Vert_{C(\overline{\Omega})}\leq 1-\tilde{\delta}.
	\end{equation*}	
	This result will be coherent with the fact that the strong solution $\phi$ will be strictly separated on the entire interval of existence (see \eqref{separ}).
	Also, thanks to standard elliptic regularity, we have that $\phi_0\in H^3(\Omega)$, being $\mu_0+\Psi^\prime(\phi_0)\in V$.
	\label{separated}
\end{remark}
\section{Main result}
\label{main2}
We shall prove the following
\begin{theorem}
	\label{main}
	(Local well-posedness in 2D). Let $\Omega$ be a bounded domain in $\R^2$ with a smooth boundary and suppose that the assumptions stated in Subsection \ref{Assumptions} hold. Moreover, assume
	\begin{align}
	\frac{\vert \nu_1-\nu_2\vert }{\nu_1+\nu_2}\leq\frac{1}{C},
	\label{small}
	\end{align}
	for a positive constant $C$ depending only on the domain $\Omega$. Then, for any initial
	datum $(\textbf{u}_0,{\psi}_0,\phi_0)$ satisfying \ref{Hbegin}-\ref{Hend}, there exists a positive time $\overline{T}_M=\overline{T}_M(\varepsilon)$, which depends only on suitable norms of the initial data and on $\varepsilon$, and a solution $(\textbf{u},\phi,\psi)$ to problem \eqref{syst}-\eqref{IBC} on $[0, T_M ]$, for any $T_M<\overline{T}_M$, such that
	\begin{align}
	\begin{cases}
	\textbf{u}\in L^\infty(0,T_M;\textbf{V}_\sigma)\cap H^1(0,T_M;\textbf{V}_\sigma^{\frac{1}{2}-\varepsilon})\cap W^{1,\infty}(0,T_M;\textbf{V}_\sigma^{-\frac{1}{2}-\varepsilon})\cap L^{\frac{p}{p-1}}(0,T_M;[W^{2,p}(\Omega)]^2),\\
	\phi\in L^\infty(0,T_M;W^{2,p}(\Omega))\cap  H^1(0,T_M;V),\\
	{\psi}\in L^\infty(0,T_M; [W^{2,p}(\Omega)]^2)\cap W^{1,\infty}(0,T_M;[H^1(\Omega)]^2),\\
	{\mu}\in L^\infty(0,T_M;V)\cap L^2(0,T_M;H^3(\Omega)),\\
	\vert\phi(x,t)\vert\leq 1\quad\text{in }\Omega\times[0,T_M],
	\end{cases}
	\label{regularitiess}
	\end{align}
	for any $p\in\left[2,\frac{4}{1+2\varepsilon}\right]$.
	Moreover, $\phi$ is strictly separated from the pure phases, i.e., there exists $\hat{\delta}>0$ such that
	\begin{equation}
	\max_{t\in[0,T_M]}\Vert \phi\Vert_{C(\overline{\Omega})}\leq 1-\hat{\delta}.
	\label{separ}
	\end{equation}	
	If, in addition, $\lambda^{\prime\prime}\in L^\infty(\R)$, then the solution is unique and depends continuously on the initial data norms in $\textbf{H}_\sigma\times H\times [H^1(\Omega)]^2$.
	\end{theorem}

\begin{remark}
\label{pressure2}
On account of the regularity properties of the local solution, arguing as in \cite{temam}, we can obtain the existence of the pressure $\pi$ up to a constant. Indeed, we have
	\begin{equation*}
	\textbf{f}=-\partial_t\textbf{u}-(\textbf{u}\cdot\nabla)\textbf{u}+\text{div}(\nu(\phi) D\textbf{u})-\text{div}(\lambda(\phi)\nabla\psi^T\nabla\psi)-\text{div}(\nabla\phi\otimes\nabla\phi)-\nu(\phi)\frac{(1-\phi)\textbf{u}}{2k(\phi)},
	\end{equation*}
	with $\textbf{f}\in L^2(0,T;[L^2(\Omega)]^2)$. Then we deduce that there exists $\pi\in L^2(0,T;V)$, with $\overline{\pi}\equiv 0$, and it satisfies $\nabla\pi=\textbf{f}$ almost everywhere in $\Omega\times(0,T)$. Thus the local solution satisfies the following problem
\begin{equation*}
\begin{cases}
\partial_t \textbf{u}+(\textbf{u}\cdot\nabla)\textbf{u}+\nabla \pi-\text{div}(\nu(\phi) D\textbf{u})+\text{div}(\lambda(\phi)\nabla\psi^T\nabla\psi)\\
+\text{div}(\nabla\phi\otimes\nabla\phi)+\nu(\phi)\frac{(1-\phi)\textbf{u}}{2k(\phi)}=\textbf{0},\\
\text{div}\,\textbf{u}=0,\\
\partial_t \psi +\textbf{u}\cdot \nabla\psi=\textbf{0},\\
\partial_t \phi+\textbf{u}\cdot \nabla\phi-\Delta\mu=0,\\
\mu=-\Delta\phi+\Psi^\prime(\phi)+\frac{\lambda^\prime(\phi)}{2}\vert\nabla \psi\vert^2,
\end{cases}
\end{equation*}
almost everywhere in $\Omega \times (0,T_M)$, with initial and boundary conditions
\begin{equation*}
\begin{cases}
\textbf{u}(0)=\textbf{u}_0,\quad {\phi}(0)={\phi}_0 ,\quad {\psi}(0)={\psi}_0, \quad\text{ a.e. in }\Omega,\\
\textbf{u} = \textbf{0}, \quad \partial_\textbf{n} \phi= \partial_\textbf{n}\mu=0, \quad\text{ a.e. on } \partial\Omega\times (0,T_M).
\end{cases}
\end{equation*}

\end{remark}

\begin{remark}
		Notice that the condition $\varepsilon>0$ cannot be relaxed to $\varepsilon=0$, due to the critical fact that $\textbf{V}_\sigma^{1/2}\subsetneqq [H^{1/2}(\Omega)]^2\cap \textbf{H}_\sigma$ (see, e.g., \cite{salgado} and the fourth estimate in Subsection \ref{ests} below). Therefore, even in the case $\psi\in [W^{2,4}(\Omega)]^2$ we could not reach the regularity $\textbf{u}\in H^1(0,T_M;\textbf{V}_\sigma^{\frac{1}{2}})\cap W^{1,\infty}(0,T_M;\textbf{V}_\sigma^{-\frac{1}{2}})$ and thus neither the integrability exponent $p=4$ can be achieved.
	\end{remark}
\begin{remark}
	\label{rmm}
Notice that  $\phi \in C([0,T_M];C(\overline{\Omega}))$. Indeed, $\phi\in L^2(0,T_M; H^2(\Omega))\cap H^1(0,T_M;V)$, so that, by \cite[Thm.II.5.14,][]{Boyerlibro}, we get $\phi \in C([0,T_M];H^{3/2}(\Omega))$ but $H^{3/2}(\Omega)\hookrightarrow C(\overline{\Omega})$.
\end{remark}
\begin{remark}
If a strong solution exists, then for proving uniqueness we just need $\lambda\in C^{1,1}(\R)$.
\end{remark}

\section{Technical tools}
\label{tools}
\paragraph{General agreement.} The symbol $C>0$ will denote a generic constant, depending
only on the structural parameters of the problem and the maximal time ($T_M$ or $T_0$, accordingly), but independent of $t$ and of the approximating indices (e.g., $\alpha,\xi$). The dependence of these indices will be explicitly indicated by a corresponding symbol, e.g., $C_\alpha, C_\xi$.
\subsection{Main $\alpha$ approximations}
We introduce a suitable mollification of $\lambda$, say, $\lambda_\alpha$, $\alpha>0$, such that $\lambda_\alpha,\lambda^\prime_\alpha, \lambda^{\prime\prime}_\alpha$ are bounded and $\lambda_\alpha,\lambda^\prime_\alpha$ are globally Lipschitz in $\R$. In particular, by a well-known property, for $i=0,1$, we have
$$
\Vert\lambda_\alpha^{(i)}\Vert_{L^\infty(\R)}\leq \Vert\lambda^{(i)}\Vert_{L^\infty(\R)}.
$$
Thus
\begin{align}
\Vert\lambda_\alpha\Vert_{L^\infty(\R)}\leq \lambda_0:=\max\{\lambda_*,\lambda^*\},\qquad  \Vert\lambda_\alpha^\prime\Vert_{L^\infty(\R)}\leq \lambda_1:=\max\{\lambda_{1,*},\lambda^{1,*}\},
\label{lambdaest}
\end{align}
Moreover, being the mollifier nonnegative, we have
\begin{align}
0<\tilde{\lambda}_*\leq \lambda_\alpha(s),\qquad 0\leq \lambda_\alpha^{\prime\prime}(s)\qquad\forall s\in \R,
\label{lambdamin}
\end{align}
for some $\tilde{\lambda}_*>0$ independent of $\alpha$.
Notice also that, being $\lambda\in C^2(\R)$, we have that $\lambda_\alpha,\lambda^\prime_\alpha,\lambda^{\prime\prime}_\alpha$ converge to $\lambda,\lambda^\prime,\lambda^{\prime\prime}$ uniformly on compact sets of $\R$.

We also set, for $\alpha>0$, $\textbf{v}=(I+\alpha \textbf{A})^{-2}\textbf{u}$ and $\textbf{w}=(I+\alpha \textbf{A})^{-1}\textbf{u}$, i.e.
\begin{align}
\textbf{v}+\alpha\textbf{A}\textbf{v}=\textbf{w},\qquad \textbf{w}+\alpha\textbf{A}\textbf{w}=\textbf{u}.
\label{w}
\end{align}
We know that in general, for any $p\in [2,\infty)$ (see estimates for the resolvent of the Stokes operator, with $\lambda=\frac{1}{\alpha}$, e.g., in \cite[Ch.3, Thm.8]{Saal})
\begin{align}
\Vert\textbf{w}\Vert_{[L^p(\Omega)]^2}\leq C\Vert\textbf{u}\Vert_{[L^p(\Omega)]^2},
\label{lp1}
\end{align}
with $C>0$ \textit{independent} of $\alpha\in (0,1)$.

If $\textbf{u}\in D(\textbf{A})$, on account of the boundary condition, we have $\textbf{w}\in D(\textbf{A})$. Thus we immediately deduce from \eqref{w} that $\textbf{A}\textbf{w}\in D(\textbf{A})$. Therefore, we get
$$
\textbf{A}\textbf{w}+\alpha\textbf{A}^2\textbf{w}=\textbf{A}\textbf{u}.
$$
and setting $\textbf{z}=\textbf{Aw}$ we come back to
$$
\textbf{z}+\alpha\textbf{A}\textbf{z}=\textbf{A}\textbf{u}.
$$
Hence, using again the estimates for the resolvent of the Stokes operator, we find
$$
\Vert\textbf{z}\Vert_{[L^p(\Omega)]^2}\leq C\Vert\textbf{Au}\Vert_{[L^p(\Omega)]^2},
$$
so that
\begin{align}
\Vert\textbf{w}\Vert_{[W^{2,p}(\Omega)]^2}\leq C\Vert\textbf{u}\Vert_{[{W}^{2,p}(\Omega)]^2}.
\label{w2p}
\end{align}
Similarly, by considering the equation for $\textbf{z}=\textbf{A}^{3/4-\varepsilon/2}\textbf{w}$, for any $0\leq \varepsilon\leq 3/2$, we get
$$
\textbf{z}+\alpha\textbf{A}\textbf{z}=\textbf{A}^{3/4-\varepsilon/2}\textbf{u}.
$$
This gives
\begin{align}
\Vert{\textbf{A}}^{3/4-\varepsilon/2}\textbf{w}\Vert_{[L^p(\Omega)]^2}\leq C\Vert\textbf{A}^{3/4-\varepsilon/2}\textbf{u}\Vert_{[L^p(\Omega)]^2}.
\label{H32}
\end{align}
Moreover, starting from \eqref{w}, multiplying it by $\textbf{A}\textbf{w}$ and integrating over $\Omega$, we get, by Young's inequality,
$$
\frac{1}{2}\Vert\textbf{w}\Vert_{\textbf{V}_\sigma}^2+{\alpha} \Vert\textbf{A}\textbf{w}\Vert^2\leq \frac{1}{2}\Vert\textbf{u}\Vert^2_{\textbf{V}_\sigma},
$$
i.e.,
\begin{align}
\Vert\textbf{w}\Vert_{\textbf{V}_\sigma}\leq C\Vert\textbf{u}\Vert_{\textbf{V}_\sigma}.
\label{h1}
\end{align}
Notice that the constants $C$ are independent of $\alpha$.

We know that, in general, for some constant depending on $\alpha$, for any $k\geq0$ and for any $p\in [2,\infty)$, the following inequalities hold
$$
\Vert\textbf{w}\Vert_{[H^{k+2}(\Omega)]^2}\leq C_\alpha\Vert\textbf{u}\Vert_{[H^k(\Omega)]^2},
$$
$$
\Vert\textbf{w}\Vert_{[W^{k+2,p}(\Omega)]^2}\leq C_\alpha\Vert\textbf{u}\Vert_{[W^{k,p}(\Omega)]^2},
$$
Therefore, applying the same results to $\textbf{v}=(I+\alpha \textbf{A})^{-1}\textbf{w}$, we deduce
\begin{align}
&	\Vert\textbf{v}\Vert_{[L^{p}(\Omega)]^2}\leq C\Vert\textbf{w}\Vert_{[L^p(\Omega)]^2}\leq C\Vert\textbf{u}\Vert_{[L^p(\Omega)]^2},
\label{lp}\\&
\Vert\textbf{v}\Vert_{\textbf{V}_\sigma}\leq C\Vert\textbf{w}\Vert_{\textbf{V}_\sigma}\leq C\Vert\textbf{u}\Vert_{\textbf{V}_\sigma},
\label{h1v}\\&
\Vert\textbf{v}\Vert_{[W^{2,p}(\Omega)]^2}\leq C\Vert\textbf{w}\Vert_{[W^{2,p}(\Omega)]^2}\leq C\Vert\textbf{u}\Vert_{[{W}^{2,p}(\Omega)]^2},
\label{w2pv}\\&
\Vert{\textbf{A}}^{3/4-\varepsilon/2}\textbf{v}\Vert_{[L^p(\Omega)]^2}\leq C\Vert\textbf{A}^{3/4-\varepsilon/2}\textbf{w}\Vert_{[L^p(\Omega)]^2}\leq C\Vert\textbf{A}^{3/4-\varepsilon/2}\textbf{u}\Vert_{[L^p(\Omega)]^2}.\label{H3_2}
\end{align}
In addition, for any $k\geq0$ and for any $p\in [2,\infty)$, we get
\begin{align}
&\Vert\textbf{v}\Vert_{[H^{k+4}(\Omega)]^2}\leq C_\alpha\Vert\textbf{w}\Vert_{[H^{k+2}(\Omega)]^2}\leq C_\alpha\Vert\textbf{u}\Vert_{[H^k(\Omega)]^2},
\label{calpha1}\\&
\Vert\textbf{v}\Vert_{[W^{k+4,p}(\Omega)]^2}\leq C_\alpha\Vert\textbf{w}\Vert_{[W^{k+2,p}(\Omega)]^2}\leq C_\alpha\Vert\textbf{u}\Vert_{[W^{k,p}(\Omega)]^2},
\label{calpha2}
\end{align}
and, observing that $\Vert\textbf{w}\Vert_{\textbf{V}^\prime_\sigma}\leq C\Vert\textbf{u}\Vert_{\textbf{V}^\prime_\sigma}$, we have
\begin{align}
\Vert\textbf{v}\Vert_{\textbf{V}_\sigma}\leq C_\alpha \Vert\textbf{u}\Vert_{\textbf{V}^\prime_\sigma}.
\label{dual}
\end{align}
Note that the same results also hold also if we consider $\partial_t\textbf{v}$ and $\partial_t\textbf{u}$, passing through the difference quotients
(see, e.g., \eqref{dtu} below).

We now consider the linear unbounded operator $A_1  =  - \Delta  +  I$
on $H$ with domain $D(A_1) = \{ u \in  H^2 (\Omega ) : \partial_{\textbf{n}} u = 0 \text{ on } \partial \Omega \}$. It is well-known that
$A_1$ is positive and self-adjoint in $H$ with compact inverse (see, e.g.,
\cite[Ch. II, sec.2.2]{temam2}).
Similarly to the Stokes operator, for any given $\alpha>0$ and any given $\phi\in D(A_1)$, we set $\varphi_\alpha=(I+\alpha A_1)^{-1}\phi$, i.e.,
\begin{align}
\alpha A_1\varphi_\alpha+\varphi_\alpha=\phi,
\label{phialpha}
\end{align}
which can be rewritten as
$$
-\alpha\Delta\varphi_\alpha+(1+\alpha)\varphi_\alpha = \phi.
$$
It is immediate to deduce, for any $p\in[2,\infty)$, by multiplying by $\vert\varphi_\alpha\vert^{p-2}\varphi_\alpha$ and integrating over $\Omega$, that
\begin{align}
\Vert\varphi_\alpha\Vert_{L^p(\Omega)}\leq \Vert\phi\Vert_{L^p(\Omega)},
\label{lpphi}
\end{align}
with $C>0$ independent of $p$ and $\alpha$. This implies
\begin{align}
\Vert\varphi_\alpha\Vert_{L^\infty(\Omega)}\leq \Vert\phi\Vert_{L^\infty(\Omega)}.
\label{limit}
\end{align}
From \eqref{phialpha} we infer that $A_1\varphi_\alpha\in D(A_1)$. Therefore, owing to \eqref{lpphi}, we infer
$$
\Vert A_1\varphi_\alpha\Vert_{L^p(\Omega)}\leq C\Vert A_1\phi\Vert_{L^p(\Omega)},
$$
i.e.,
\begin{align}
\Vert\varphi_\alpha\Vert_{W^{2,p}(\Omega)}\leq C\Vert\phi\Vert_{W^{2,p}(\Omega)},
\label{lpphi2}
\end{align}
with $C>0$ independent of $\alpha$ and $p$. In conclusion, by standard elliptic regularity, we get, for any $k\in \N\cup \{0\}$,
\begin{align}
\Vert\varphi_\alpha\Vert_{H^{k+2}(\Omega)}\leq C_\alpha\Vert\phi\Vert_{H^k(\Omega)},
\label{elliptic}
\end{align}
with $C_\alpha>0$ depending on $\alpha$.

Here below we report some technical results which will be useful in the sequel.
\subsection{The stationary Stokes problem}
\textbf{The homogeneous Stokes problem with nonconstant viscosity.}
We consider the homogeneous Stokes problem with nonconstant viscosity
depending on a given measurable function $\varphi$. The system reads as follows
\begin{align}
\begin{cases}
-\text{div}(\nu(\varphi)D\textbf{u})+\nabla \pi=\textbf{f},\quad \text{in }\Omega \\
\text{div}\,\textbf{u}=0,\quad \text{in }\Omega,
\end{cases}
\label{NavierStokes}
\end{align}
with no-slip boundary conditions, where the coefficient $\nu$  fulfils the assumptions stated in \eqref{nuA}.
We have the following Theorem, whose straightforward proof can be found, e.g., in \cite[Thm A.1]{GGW}.
\begin{theorem}
	Let  $\textbf{f}\in [{L}^p(\Omega)]^2$ for some $p\in(1,\infty)$ and $\varphi\in W^{1,r}(\Omega)$, with $r>2$. Consider the (unique) weak solution $\textbf{u}\in \textbf{V}_\sigma$ such that $(\nu(\varphi)D\textbf{u},D\textbf{v})=(\textbf{f},\textbf{v})$ for all $\textbf{v}\in \textbf{V}_\sigma$. Suppose also that $\textbf{u}\in \textbf{V}_\sigma\cap [W^{1,s}(\Omega)]^2$ with $s\geq2$ such that $\frac{1}{p}=\frac{1}{s}+\frac{1}{r}$, $r\geq \frac{2s}{s-1}$. Then there exists $C=C(p,s,\Omega)>0$ such that
	\begin{align}
	\Vert\textbf{u}\Vert_{[W^{2,p}(\Omega)]^2}\leq C\Vert\textbf{f}\Vert_{[L^p(\Omega)]^2}+C\Vert\nabla\varphi\Vert_{[{L}^r(\Omega)]^2}\Vert D\textbf{u}\Vert_{[L^{s}(\Omega)]^4}.
	\label{H2}
	\end{align}
	\label{stokes}
\end{theorem}
\textbf{On the Stokes operator.}
First we state a useful Lemma, which is an adaptation of \cite[Lemma B.2]{Giorginitemam}. We know that, for any $\textbf{g}\in \textbf{H}_\sigma$, there exist $\textbf{A}^{-1}\textbf{g}\in D(\textbf{A})$ and $p\in V$ that solve
    \begin{align}
    -\Delta\textbf{A}^{-1}\textbf{g}+\nabla p=\textbf{g}.
    \label{st}
    \end{align}
    In turn this entails that $\textbf{A}\textbf{A}^{-1}\textbf{g}=\textbf{g}$. Then we have
    \begin{lemma}
    	Let $d=2,3$ and $\textbf{f}\in \textbf{V}_\sigma^{-1/2}$. For any $\varepsilon>0$, there exists a positive constant $C$, depending on $\varepsilon$ but
    independent of $\textbf{f}$, such that
    	\begin{align}
    	\Vert p\Vert\leq C\left(\Vert\textbf{f}\Vert_{\textbf{V}_\sigma^{-1}}+\Vert\textbf{f}\Vert_{\textbf{V}_\sigma^{-1/2+\varepsilon}}\right),
    	\end{align}
    	where $p$ is defined by \eqref{st} with $\textbf{g}=\textbf{f}$.
    	\label{press}
    \end{lemma}

    \begin{proof}
    	First, in order to have a function in $\textbf{H}_\sigma$, we consider the approximation
    	\begin{align}
    	\textbf{g}=(\textbf{I}+\alpha\textbf{A})^{-1}\textbf{f}\in
    	\textbf{H}_\sigma,\qquad \alpha>0,
    	\label{app}
    	\end{align}
    	which satisfies \eqref{st}. We need to estimate the ${H}^{-1}-$norm of $\textbf{g}$. Take $\textbf{v}\in[{H}^1_0(\Omega)]^2$ with $\Vert\textbf{v}\Vert_{[{H}^1_0(\Omega)]^2}\leq 1$. By integration by parts we have
    	\begin{align*}
    	&(\textbf{g},\textbf{v})=(\textbf{P}(-\Delta)\textbf{A}^{-1}\textbf{g},\textbf{v})=((-\Delta)\textbf{A}^{-1}\textbf{g},\textbf{P}\textbf{v})\\&=(\nabla\textbf{A}^{-1}\textbf{g},\nabla\textbf{P}\textbf{v})-\int_{\partial\Omega}\nabla\textbf{A}^{-1}\textbf{g}\ \textbf{n}\cdot \textbf{Pv}d\sigma.
    	\end{align*}
     Here $\textbf{P}: [L^2(\Omega)]^2\to \textbf{H}_\sigma$ is the Leray projector.    	
     We have that, by the properties of $\textbf{P}$ and \cite[Thm. 9.4]{Lions}, for any $\varepsilon>0$,
    	\begin{align*}
    	&\int_{\partial\Omega}\nabla\textbf{A}^{-1}\textbf{g}\ \textbf{n}\cdot \textbf{Pv}d\sigma\leq \Vert\nabla\textbf{A}^{-1}\textbf{g}\Vert_{[{L}^2(\partial\Omega)]^4}\Vert\textbf{Pv}\Vert_{[{L}^2(\partial\Omega)]^2}\\&\leq C\Vert\nabla\textbf{A}^{-1}\textbf{g}\Vert_{[{H}^{1/2+\varepsilon}(\Omega)]^4}\Vert\textbf{v}\Vert_{[{H}^1_0(\Omega)]^2}\leq C\Vert\textbf{A}^{-1}\textbf{g}\Vert_{[{H}^{3/2+\varepsilon}(\Omega)]^2}\Vert\textbf{v}\Vert_{[{H}^1_0(\Omega)]^2}\\&\leq C\Vert\textbf{A}^{3/4+\varepsilon/2}\textbf{A}^{-1}\textbf{g}\Vert\Vert\textbf{v}\Vert_{[{H}^1_0(\Omega)]^2}=C\Vert\textbf{g}\Vert_{\textbf{V}_\sigma^{-1/2+\varepsilon}}\Vert\textbf{v}\Vert_{[{H}^1_0(\Omega)]^2}.
    	\end{align*}
    	Here we have used the equivalence of the norms in $\textbf{H}^{3/2+\varepsilon}(\Omega)$, being $\textbf{A}^{-1}\textbf{g}\in D(\textbf{A}^{3/4+\varepsilon/2})$. The presence of $\varepsilon>0$ is essential to guarantee the validity of the estimate (see, e.g., \cite[Thm. 9.5]{Lions})
    	
    	Therefore, being $(\nabla\textbf{A}^{-1}\textbf{g},\nabla\textbf{P}\textbf{v})\leq C\Vert\nabla\textbf{A}^{-1}\textbf{g}\Vert\Vert\textbf{v}\Vert_{[{H}^1_0(\Omega)]^2}$,
    	we infer
    	\begin{align}
    		\Vert\textbf{g}\Vert_{[{H}^{-1}(\Omega)]^2}\leq C\left(\Vert\nabla\textbf{A}^{-1}\textbf{g}\Vert+\Vert\textbf{g}\Vert_{\textbf{V}_\sigma^{-1/2+\varepsilon}}\right).
    		\label{g}
    	\end{align}
    	In order to pass to the limit as $\alpha\rightarrow0$, we rename $\textbf{g}$ as $\textbf{g}_\alpha$ and observe that, by \eqref{app},
    		\begin{align}
    	\textbf{f}=\textbf{g}_\alpha+\alpha\textbf{A}\textbf{g}_\alpha.
    	\label{app2}
    	\end{align}
    	Thus multiplying \eqref{app2}  by $\textbf{A}^{-5/2}(\textbf{f}-\textbf{g}_\alpha)$ and integrating over $\Omega$, we get
    	$$
    	\Vert\textbf{A}^{-5/4}(\textbf{f}-\textbf{g}_\alpha)\Vert\leq \alpha \Vert\textbf{A}^{-1/4}\textbf{g}_\alpha\Vert^2.
    	$$
    	Moreover, if we multiply \eqref{app2} by $\textbf{A}^{-1/2}\textbf{g}_\alpha$ then we get
    	\begin{align}	
        \frac{1}{2}\Vert\textbf{A}^{-1/4}\textbf{g}_\alpha\Vert^2+\alpha\Vert\textbf{A}^{1/4}\textbf{g}_\alpha\Vert^2\leq \frac{1}{2}\Vert\textbf{A}^{-1/4}\textbf{f}\Vert^2.
    	\label{bound}
    		\end{align}
    	Therefore we have
        $$
    	\Vert\textbf{A}^{-5/4}(\textbf{f}-\textbf{g}_\alpha)\Vert^2\leq \alpha \Vert\textbf{A}^{-1/4}\textbf{f}\Vert^2,
    	$$
    	and, being $\textbf{f}\in {\textbf{V}_\sigma^{-1/2}}$, we deduce that
    	$$
    	\textbf{g}_\alpha\rightarrow\textbf{f}\quad\text{ in }\textbf{V}_\sigma^{-5/2}, \quad \alpha \to 0.
    	$$
    	This, together with \eqref{bound} and the embedding $\textbf{V}^{-1/2+\varepsilon}_\sigma\hookrightarrow \textbf{V}^{-5/2}_\sigma$, gives, up to subsequences, as $\alpha\to 0$,
    		\begin{align}
    	\textbf{A}^{-1/4+\varepsilon/2}\textbf{g}_\alpha\rightarrow\textbf{A}^{-1/4+\varepsilon/2}\textbf{f}\quad\text{ in }[{L}^2(\Omega)]^2,
    	\text{ i.e., }\;
    	\textbf{g}_\alpha\rightarrow\textbf{f}\quad\text{ in }\textbf{V}_\sigma^{-1/2+\varepsilon},
    	\label{conv1}
    	\end{align}
    	and, as a consequence,
    		\begin{align}
    	\nabla\textbf{A}^{-1}\textbf{g}_\alpha\rightarrow	\nabla\textbf{A}^{-1}\textbf{f}\quad\text{ in }[{L}^2(\Omega)]^{2\times 2},
    	\text{ i.e., }\;
    	\textbf{g}_\alpha\rightarrow\textbf{f}\quad\text{ in }\textbf{V}_\sigma^{-1}.
    	\label{conv2}
    	\end{align}
    	 We can then pass to the limit by a standard procedure in \eqref{g}. Indeed, we see by \eqref{g} and \eqref{conv1} that $\{\textbf{g}_\alpha\}_\alpha$ is a Cauchy sequence in $[{H}^{-1}(\Omega)]^2$, which implies $\textbf{g}_\alpha\rightarrow\tilde{\textbf{f}}$ in $[{H}^{-1}(\Omega)]^2$ for some $\tilde{\textbf{f}}\in [{H}^{-1}(\Omega)]^2$. Due to the estimate $\Vert\textbf{h}\Vert_{\textbf{V}_\sigma^{-1}}\leq C\Vert\textbf{h}\Vert_{[{H}^{-1}(\Omega)]^2}$, we then deduce, thanks to \eqref{conv1} and by uniqueness of the limit in $\textbf{V}_\sigma^{-1}$, that $\tilde{\textbf{f}}=\textbf{f}$. Passing then to the limit as $\alpha\to0$ in \eqref{g}, we infer that
    	\begin{align}
    	\Vert\textbf{f}\Vert_{[{H}^{-1}(\Omega)]^2}\leq C\left(\Vert\nabla\textbf{A}^{-1}\textbf{f}\Vert+\Vert\textbf{f}\Vert_{\textbf{V}_\sigma^{-1/2+\varepsilon}}\right)\leq C\left(\Vert\textbf{f}\Vert_{\textbf{V}_\sigma^{-1}}+\Vert\textbf{f}\Vert_{\textbf{V}_\sigma^{-1/2+\varepsilon}}\right).
    	\label{end}
    	\end{align}
    	To conclude the proof, we recall that $\Vert p\Vert\leq C\Vert \textbf{f}\Vert_{[{H}^{-1}(\Omega)]^2}$ and we apply \eqref{end}.
    \end{proof}

We also need to collect some useful estimates which will be needed below. We recall the interpolation (see, e.g., \cite{Lions}): given $a<b<c$, we have, for any $\textbf{f}\in \textbf{V}_\sigma^{c}$,
$$
\Vert \textbf{f}\Vert_{\textbf{V}_\sigma^b}\leq C\Vert\textbf{f}\Vert_{\textbf{V}_\sigma^c}^{\overline{\omega}}\Vert\textbf{f}\Vert_{\textbf{V}_\sigma^a}^{1-\overline{\omega}},
$$
where $\overline{\omega}=\frac{b-a}{c-a}$. From this result we get, for $\partial_t\textbf{u}$ sufficiently regular (e.g., belonging to $\textbf{V}_\sigma$ for almost any $t>0$), that
\begin{align}
\Vert \partial_t\textbf{u}\Vert\leq C\Vert\partial_t\textbf{u}\Vert^{1/2-\varepsilon}_{\textbf{V}_\sigma^{-1/2-\varepsilon}}\Vert\partial_t\textbf{u}\Vert^{1/2+\varepsilon}_{\textbf{V}_\sigma^{1/2-\varepsilon}},
\label{interp1}
\end{align}
for any $0\leq\varepsilon\leq 1/2$. Similarly, we obtain
\begin{align}
\Vert \partial_t\textbf{u}\Vert_{\textbf{V}_\sigma^{-2\varepsilon}}\leq C\Vert\partial_t\textbf{u}\Vert^{1/2+\varepsilon}_{\textbf{V}_\sigma^{-1/2-\varepsilon}}\Vert\partial_t\textbf{u}\Vert^{1/2-\varepsilon}_{\textbf{V}_\sigma^{1/2-\varepsilon}}.
\label{interp5}
\end{align}
Moreover, by Sobolev-Gagliardo-Nirenberg's inequality, the equivalence of the $H^s$-norms and \eqref{interp1}, for any $0\leq\varepsilon\leq 1/6$ we have
\begin{align}
\Vert\partial_t\textbf{u}\Vert_{[L^3(\Omega)]^2}\leq C\Vert\partial_t\textbf{u}\Vert_{\textbf{V}_\sigma^{1/2-\varepsilon}}^{\frac{2}{3-6\varepsilon}}\Vert\partial_t\textbf{u}\Vert^{\frac{1-6\varepsilon}{3-6\varepsilon}}\leq C\Vert\partial_t\textbf{u}\Vert_{\textbf{V}_\sigma^{1/2-\varepsilon}}^{\frac{5-12\varepsilon^2-4\varepsilon}{2(3-6\varepsilon)}}\Vert\partial_t\textbf{u}\Vert_{\textbf{V}_\sigma^{-1/2-\varepsilon}}^{\frac{1+12\varepsilon^2-8\varepsilon}{2(3-6\varepsilon)}}.
\label{interp2}
\end{align}
Moreover, recalling \eqref{H3_2}, we find
\begin{align}
&\nonumber\Vert\nabla(I+\alpha\textbf{A})^{-2}\textbf{A}^{-1/2}\partial_t\textbf{u}\Vert_{[L^3(\Omega)]^4}\leq \Vert(I+\alpha\textbf{A})^{-2}\textbf{A}^{-1/2}\partial_t\textbf{u}\Vert_{[W^{1,3}(\Omega)]^2}\\&\leq C\Vert(I+\alpha\textbf{A})^{-2}\textbf{A}^{-1/2}\partial_t\textbf{u}\Vert_{\textbf{V}_\sigma^{3/2-\varepsilon}}^{\frac{2}{3-6\varepsilon}}\Vert(I+\alpha\textbf{A})^{-2}\textbf{A}^{-1/2}\partial_t\textbf{u}\Vert_{\textbf{V}_\sigma}^{\frac{1-6\varepsilon}{3-6\varepsilon}}\nonumber\\&\leq\nonumber C\Vert\textbf{A}^{-1/2}\partial_t\textbf{u}\Vert_{\textbf{V}_\sigma^{3/2-\varepsilon}}^{\frac{2}{3-6\varepsilon}}\Vert\textbf{A}^{-1/2}\partial_t\textbf{u}\Vert_{\textbf{V}_\sigma}^{\frac{1-6\varepsilon}{3-6\varepsilon}}\nonumber\\&=C\Vert\partial_t\textbf{u}\Vert_{\textbf{V}_\sigma^{1/2-\varepsilon}}^{\frac{2}{3-6\varepsilon}}\Vert\partial_t\textbf{u}\Vert^{\frac{1-6\varepsilon}{3-6\varepsilon}}\leq C\Vert\partial_t\textbf{u}\Vert_{\textbf{V}_\sigma^{1/2-\varepsilon}}^{\frac{5-12\varepsilon^2-4\varepsilon}{2(3-6\varepsilon)}}\Vert\partial_t\textbf{u}\Vert_{\textbf{V}_\sigma^{-1/2-\varepsilon}}^{\frac{1+12\varepsilon^2-8\varepsilon}{2(3-6\varepsilon)}},
\label{l3}
\end{align}
independently of $\alpha$.
In conclusion, we have, by {\color{blue}2D} Sobolev embeddings and the equivalence of the norms in $H^{1/2-\varepsilon}$ and $\textbf{V}_\sigma^{1/2-\varepsilon}$,
\begin{align}
\Vert\partial_t\textbf{u}\Vert_{[L^{\frac{4}{1+2\varepsilon}}(\Omega)]^2}\leq C\Vert\partial_t\textbf{u}\Vert_{\textbf{V}_\sigma^{1/2-\varepsilon}},
\label{interp3}
\end{align}
for any $\varepsilon\geq0$. Thus, we infer
\begin{align}
\Vert \textbf{A}^{-1/2-\varepsilon}\partial_t\textbf{u}\Vert_{[L^{\frac{4}{1+2\varepsilon}}(\Omega)]^2}\leq C\Vert \textbf{A}^{-1/2-\varepsilon}\partial_t\textbf{u}\Vert_{\textbf{V}_\sigma^{1/2-\varepsilon}}\leq C\Vert\partial_t\textbf{u}\Vert_{\textbf{V}_\sigma^{-1/2-3\varepsilon}},
\label{interp4}
\end{align}
and
\begin{align}
&\nonumber\Vert \nabla\textbf{A}^{-1/2-\varepsilon}\partial_t\textbf{u}\Vert_{[L^{\frac{4}{1+2\varepsilon}}(\Omega)]^4}\leq C\Vert \textbf{A}^{-1/2-\varepsilon}\partial_t\textbf{u}\Vert_{[W^{1,\frac{4}{1+2\varepsilon}}(\Omega)]^2}\\&\leq C\Vert \textbf{A}^{-1/2-\varepsilon}\partial_t\textbf{u}\Vert_{\textbf{V}_\sigma^{3/2-\varepsilon}}= C\Vert\partial_t\textbf{u}\Vert_{\textbf{V}_\sigma^{1/2-3\varepsilon}}\leq C\Vert\partial_t\textbf{u}\Vert_{\textbf{V}_\sigma^{1/2-\varepsilon}} .
\label{interp4bis3}
\end{align}
\subsection{Main results on the decoupled regularized problems}
\label{three}
This subsection contains some results on the three decoupled problems which will be exploited in the fixed point argument (see Subsection \ref{fix} below). More precisely, they are related to the transport problem for $\psi$, the Navier-Stokes equations for $\textbf{u}$ and a generalized Cahn-Hilliard equation for $\phi$. For the ease of readability, their proofs are postponed in Section \ref{proofs}.
    \subsubsection{The advection equation}
    First we recall a result which will be used when dealing with pure transport equations (see, e.g., \cite[Lemma A.4,][]{brezisbourg}).
    \begin{lemma}
    	Let $\beta\geq2$ be an integer, and let $1\leq p\leq q\leq +\infty$ such that $\beta>\frac{2}{q}+1$. Let $P\in W^{\beta,p}(\Omega)$ and $G\in W^{\beta,q}(\Omega;\R^2)$ and $G$ is a $C^1$ diffeomorphism from $\overline{\Omega}$ onto $\overline{\Omega}$.Then $P\circ G\in W^{\beta,p}(\Omega)$ and
    	\begin{equation}
    	\Vert P\circ G\Vert_{W^{\beta,p}(\Omega)}\leq C\Vert P\Vert_{W^{\beta,p}(\Omega)}\dfrac{1}{\inf_{\Omega} \vert \textbf{J}_G\vert^{1/p}}(\Vert G\Vert^\beta_{W^{\beta,q}(\Omega)}+1),
    	\label{bourg}
    	\end{equation}
    	where $C$ depends on $\beta, p,q,\Omega$ and $\textbf{J}_G$ is the Jacobian of $G$.
    \end{lemma}
The following lemma about the transport equation can be proved by slightly adapting \cite[Lemma 1.2,][]{Lady} (see Section \ref{proof1}).
	\begin{lemma}
		Let $\textbf{v}\in C([0,T];C^2(\overline{\Omega}))$ for some $T>0$ and $\psi_0\in C^\infty(\overline{\Omega})$. Then, for any $0\leq t\leq T$ there exists a unique $X(\cdot,t,\cdot)\in C^1([0,T];C^2(\overline{\Omega}))$, which is a $C^1$-diffeomorphism from $\overline{\Omega}$ onto $\overline{\Omega}$, such that
		$$
		X(s,t,x)=x+\int_t^s\textbf{v}(X(\tau,t,x),\tau)d\tau,
		$$
		for any $s,t\in[0,T]$. Moreover, $X(0,\cdot,\cdot)\in C([0,T];C^2(\overline{\Omega}))\cap C^1([0,T];C^1(\overline{\Omega})) $. Therefore, defining
		\begin{align}
		\psi(x,t)=\psi_0(X(0,t,x)),
		\label{composition}
		\end{align}
		we have,
		\begin{align}
		\psi\in C([0,T];C^2(\overline{\Omega}))\cap C^1([0,T];C^1(\overline{\Omega})),
		\label{psireg}
		\end{align}
		\begin{align}
		\vert \psi(x,t)\vert \leq C\qquad\forall(x,t)\in\overline{\Omega}\times[0,T]
		\label{est}
		\end{align}
		and $\psi$ solves, for any $(x,t)\in \overline{\Omega}\times [0,T]$
		\begin{align*}
		\begin{cases}
			\partial_t \psi +\textbf{v}\cdot \nabla\psi=\textbf{0}\\
			\psi(0)=\psi.
		\end{cases}
		\end{align*}
		Furthermore, we have,
		     \begin{align}
		\left(\sum_{j,k=1}^2\vert\partial_jX_k(0,t,x)\vert^2\right)^{1/2}\leq Ce^{\int_0^t\Vert\textbf{v}(s)\Vert_{[W^{1,\infty}(\Omega)]^2}ds},
		\label{estW1p}
		\end{align}
		and
		          \begin{align}
		&\nonumber\left(\sum_{i,j,k=1}^2\vert\partial_i\partial_jX_k(0,t,x)\vert^2\right)^{1/2}\\&\leq Ce^{\int_0^t\Vert\textbf{v}(\tau)\Vert_{[W^{1,\infty}(\Omega)]^2}}\int_0^t\Vert\textbf{v}(s)\Vert_{[W^{2,\infty}(\Omega)]^2}e^{\int_0^s\Vert\textbf{v}(\tau)\Vert_{[W^{1,\infty}(\Omega)]^2}d\tau}ds.
		\label{estW2p}
		\end{align}
		In conclusion, if $\textbf{v}	\in C([0,T];C^3(\overline{\Omega}))$, we also have
				\begin{align}
		\psi\in  C^1([0,T];C^2(\overline{\Omega})).
		\label{psireg2}
		\end{align}
		\label{transport}
		In this case, for any $p\in(2,\infty)$, we get
		\begin{align}
			\frac{d}{dt}\frac{1}{2}\Vert{\psi}\Vert_{[W^{2,p}(\Omega)]^2}^2\leq C\Vert{\textbf{v}}\Vert_{[W^{2,p}(\Omega)]^2}\Vert{\psi}\Vert_{[W^{2,p}(\Omega)]^2}^2.
			\label{psiw2p}
		\end{align}
	
	\end{lemma}
\begin{remark}
	Note that estimate \eqref{psiw2p} is valid only with some additional regularity on the velocity field $\textbf{v}$, i.e., $C([0,T];C^3(\overline{\Omega}))$, whereas for the lower regularity case in the first part of the lemma we will resort to Lemma \ref{bourg} to get a similar estimate (see Subsection \ref{fix} below).
\end{remark}

\subsubsection{Regularized Navier-Stokes equations}
	Let $\varphi$, $\textbf{v}$, $g$ and \textbf{h} be given. Consider the following system for a given $\alpha>0$
	\begin{equation}
	\begin{cases}
	\partial_t\textbf{u}-\text{div}(\nu(\varphi)D\textbf{u})+(\textbf{v}\cdot\nabla) \textbf{u}+\nabla\pi+g\textbf{u}=(I+\alpha^2\textbf{A})^{-2}\textbf{h},\\
	\text{div }\textbf{u}=0,
	\end{cases}
	\label{NS1}
	\end{equation}
	in $\Omega\times(0,T)$, $T>0$, subject to the boundary and initial conditions
	\begin{equation}
	\begin{cases}
	\textbf{u}=0, \quad \text{ on }\partial\Omega\times(0,T),\\
	\textbf{u}(0)=\textbf{u}_0, \quad\text{ in } \Omega.
	\end{cases}
	\label{NS2}
	\end{equation}
	Then the following theorem holds.
\begin{theorem}
\label{ns}
	Let
    \begin{align*}
    &\varphi\in L^\infty(0,T;V)\cap L^4(0,T;{H}^2(\Omega)),\quad \textbf{v}\in L^2(0,T;[L^\infty(\Omega)]^2),\\
    &g \in L^\infty( 0, T ; L^\infty(\Omega)), \quad\textbf{h} \in L^2 (0, T ; \textbf{V}_\sigma^\prime ).
    \end{align*}
    Suppose also that $g\geq 0$ almost everywhere in $\Omega\times(0,T)$. Then, for every
	$\textbf{u}_0 \in \mathbb{H}^1$, problem \eqref{NS1}-\eqref{NS2} possesses a unique strong solution $\textbf{u}$ on the interval $[0, T ]$ such that	
	$$
	\textbf{u}\in L^\infty(0,T;\textbf{V}_\sigma)\cap L^2(0,T;\textbf{W}_\sigma)\cap H^1(0,T;\textbf{H}_\sigma),
	$$
	and there is $C > 0$, depending on $\Vert\textbf{u}_0\Vert_{\textbf{V}_\sigma}$, $\Vert\varphi\Vert_{L^4(0,T;{H}^2(\Omega))}$, $\Vert\varphi\Vert_{L^\infty(0,T;V)}$,
    $\Vert g\Vert_{L^\infty( 0, T ; L^\infty(\Omega))}$,  $\Vert\textbf{v}\Vert_{L^2(0,T;[L^\infty(\Omega)]^2)}$, $\Vert\textbf{h}\Vert_{ L^2 (0, T ; \textbf{V}_\sigma^\prime )}$, such that
	\begin{align}
	&\Vert\textbf{u}\Vert_{L^\infty(0,T;\textbf{V}_\sigma)}+\Vert\textbf{u}\Vert_{L^2(0,T;\textbf{W}_\sigma)}+\Vert\textbf{u}\Vert_{ H^1(0,T;\textbf{H}_\sigma)}\leq C(T).
	\end{align}	
\end{theorem}
\subsubsection{A generalized Cahn-Hilliard equation}
Here we consider the following boundary and initial value problem for a generalized Cahn-Hilliard equation:	
	\begin{equation}
	\begin{cases}
	\partial_t \phi+\textbf{v}\cdot \nabla\phi-\Delta\mu=0, \quad\text{in } \Omega\times(0,T),\\
	\mu=-\Delta\phi+\Psi^\prime(\phi)+\lambda_\alpha^\prime(\phi)f,\quad\text{in } \Omega\times(0,T),\\
	\partial_{\textbf{n}}\phi=\partial_{\textbf{n}}\mu=0,\quad\text{ on }\partial\Omega\times(0,T),\\
    {\phi}(0)={\phi}_0, \quad\text{ in }\Omega.\\
	\end{cases}
	\label{CH}
	\end{equation}
	The following result holds.
	\begin{theorem}
		\label{ch}
		Let $f\in L^2(0,T;V)\cap L^2(0,T;L^\infty(\Omega))$, $\textbf{v}\in L^2(0,T;[L^\infty(\Omega)]^2\cap \textbf{H}_\sigma)$ and $\phi_0\in V$ be such that $\Vert\phi_0\Vert_{L^\infty(\Omega)}\leq 1$ and $\vert\overline{\phi}_0\vert<1$. Then there exists a unique $\phi$ weak solution to \eqref{CH} satisfying the following properties
		\begin{align*}
		&\Vert\phi\Vert_{L^\infty(\Omega\times(0,T))}\leq 1,
		\\&
		\phi \in L^\infty(0,T;V)\cap L^4(0,T;H^2(\Omega))\cap H^1(0,T;V^\prime),\\&
		\mu\in L^2(0,T;V).
		\end{align*}
		If, in addition, $f\in  H^1(0,T;H)\cap L^\infty(0,T;V)$ and $f_0\in V$ are such that $f(0)=f_0$, together with $\textbf{v}\in H^1(0,T;\textbf{V}_\sigma)$, $\textbf{v}(0)=\textbf{v}_0\in \textbf{V}_\sigma$ and $\phi_0\in H^2(\Omega)$, with $\partial_\textbf{n}\phi_0=0$ almost everywhere on $\partial\Omega$, and $\mu_0:=-\Delta\phi_0+\Psi^\prime(\phi_0)\in V$, then, for any $q\in[2,\infty)$, we have
		\begin{align*}
		&\phi\in L^\infty(0,T;W^{2,q}(\Omega))\cap H^1(0,T;V),\\&
	    \mu\in L^\infty(0,T;V),\quad F^{\prime}(\phi)\in L^\infty(0,T;L^q(\Omega)),\quad F^{\prime\prime}(\phi)\in L^\infty(0,T;L^q(\Omega)).
		\end{align*}
	\end{theorem}
\begin{remark}
	Note that the assumptions $f(0)=f_0$ and $\textbf{v}(0)=\textbf{v}_0$ make sense, since $f\in C_w([0,T];V)$ (see, e.g., \cite[Lemma II.5.9]{Boyerlibro}) and $\textbf{v}\in C([0,T];\textbf{V}_\sigma)$.
\end{remark}
\section{Proof of Theorem \ref{main}}
The proof is organized in six subsections as follows. First we introduce an approximating problem which is solved (locally in time) through fixed point argument. The regularization is essentially based
on the Leray-$\alpha$ model. Then we perform some higher-order estimates on the approximating solution. These estimates allow us to find a series of bounds which are uniform with respect to the regularizing parameter $\alpha>0$. We can then pass to the limit letting $\alpha$ go to $0$. Finally, we establish a continuous
dependence estimate which entails uniqueness.

\label{proof}
\subsection{The approximating problem}
The approximating problem has the following form

\begin{equation}
\label{problem}
\begin{cases}
\partial_t \textbf{u}_\alpha+(\textbf{v}_\alpha \cdot\nabla)\textbf{u}_\alpha +\nabla \pi_\alpha-\text{div}(\nu(\varphi_\alpha) D\textbf{u}_\alpha )+(I+\alpha \textbf{A})^{-2} \textbf{P}\text{div}(\lambda_\alpha(\phi_\alpha )\nabla\psi_\alpha ^T\nabla\psi_\alpha )\\+(I+\alpha \textbf{A})^{-2}\textbf{P}\text{div}(\nabla\phi_\alpha \otimes\nabla\phi_\alpha )+\nu(\phi_\alpha )\frac{(1-\phi_\alpha )\textbf{u}_\alpha }{2k(\phi_\alpha )}=\textbf{0},\\
\text{div}\,\textbf{u}_\alpha =0,\\
\partial_t \psi_\alpha +\textbf{v}_\alpha \cdot \nabla\psi_\alpha =\textbf{0},\\
\partial_t \phi_\alpha +\textbf{v}_\alpha \cdot \nabla\phi_\alpha -\Delta\mu_\alpha=0,\\
\mu_\alpha=
-\Delta\phi_\alpha +\Psi^\prime(\phi_\alpha )+\frac{\lambda_\alpha^\prime(\phi_\alpha )}{2}\vert\nabla \psi_\alpha \vert^2,\\
\end{cases}
\end{equation}
in $\Omega\times(0,T)$, equipped with the following boundary and initial conditions
\begin{equation}
\label{problem_BIC}
\begin{cases}
\textbf{u}_\alpha =0, \quad
\partial_{\textbf{n}}\phi_\alpha =0, \quad
\partial_{\textbf{n}}\mu_\alpha=0, \quad\text{ on }\partial\Omega\times(0,T),\\
\textbf{u}_\alpha (0)=\textbf{u} _{0,\alpha},\quad
{\phi_\alpha }(0)={\phi_\alpha }_0,\quad
{\psi_\alpha }(0)={\psi}_{0,\alpha},\quad\text{ in }\Omega,\\
\end{cases}
\end{equation}
where
$$
\textbf{v}_\alpha =(I+\alpha\textbf{A})^{-2}\textbf{u}_\alpha,\qquad
\varphi_\alpha=(I+\alpha A_1)^{-1}\phi_\alpha.
$$
Moreover, $\psi_{0,\alpha}\in C^\infty(\overline{\Omega})$ is a suitable mollification of $\psi_0$, such that $\psi_{0,\alpha}\rightarrow\psi_0$ in $[W^{2,p}(\Omega)]^2$ as $\alpha\to 0$, and $\textbf{u}_{0,\alpha}$ is a sequence in $\textbf{W}_\sigma$ such that $\textbf{u}_{0,\alpha}\rightarrow\textbf{u}_0$ in $\textbf{V}_\sigma^{3/2+2\widetilde{\delta}}$ as $\alpha\to0$ (which exists, being $\textbf{W}_\sigma$ dense in $\textbf{V}_\sigma^{3/2+2\widetilde{\delta}}$, where $\widetilde{\delta}$ is defined in \ref{Hbegin}).

\begin{remark}
Note that $(I+\alpha\textbf{A})^{-2}\textbf{P}$ is self-adjoint. This fact allows us to obtain the energy identity. Indeed, for example,
	\begin{align*}
	((I+\alpha \textbf{A})^{-2}\textbf{P}\text{div}(\lambda(\phi_\alpha )\nabla\psi^T\nabla\psi),\textbf{u}_\alpha )&=(\text{div}(\lambda(\phi_\alpha )\nabla\psi^T\nabla\psi),(I+\alpha \textbf{A})^{-2}\textbf{u}_\alpha  )\\&=(\text{div}(\lambda(\phi_\alpha )\nabla\psi^T\nabla\psi),\textbf{v}_\alpha  ),
	\end{align*}
	which cancels out with the corresponding term obtained by handling properly the equation for $\psi$. Recall that $\textbf{P}(I+\alpha\textbf{A})^{-2}=(I+\alpha\textbf{A})^{-2}$.	
    Notice also that the logarithmic potential is not approximated as usual. This is due to the fact that we need to exploit the bound $\Vert\phi_\alpha \Vert_{L^\infty(\Omega\times(0,T))}\leq 1$.
\end{remark}

For the sake of simplicity, in the sequel we omit the index $\alpha$ unless it might create confusion.
\subsection{The fixed point argument}
\label{fix}
Let us now fix $\widetilde{\textbf{u}}\in L^\infty(0,T;\textbf{V}_\sigma)\cap H^1(0,T;\textbf{H}_\sigma) $ and $\Vert\textbf{u}_{0,\alpha}\Vert_{\textbf{V}_\sigma}=R$. Setting
$$
\widetilde{\textbf{v}}=(I+\alpha\textbf{A})^{-2}\widetilde{\textbf{u}}
$$
and recalling \eqref{calpha1}, we have
	\begin{equation}
	\widetilde{\textbf{v}}\in C([0,T]; [H^5(\Omega)]^2)\cap H^1(0,T;[H^4(\Omega)]^2)\hookrightarrow C([0,T]; C^2(\overline{\Omega}))
	\label{vtilde}
	\end{equation}
Then we consider the problem
\begin{equation}
\begin{cases}
\partial_t \widetilde{\psi}+\widetilde{\textbf{v}}\cdot \nabla\widetilde{\psi}=\textbf{0},\\
\partial_t \widetilde{\phi}+\widetilde{\textbf{v}}\cdot \nabla\widetilde{\phi}-\Delta\widetilde{\mu}=0,\\
\widetilde{\mu}=
-\Delta\widetilde{\phi}+\Psi^\prime(\widetilde{\phi})+\frac{\lambda_\alpha^\prime(\widetilde{\phi})}{2}\vert\nabla \widetilde{\psi}\vert^2,\\
\end{cases}
\label{sys0}
\end{equation}
in $\Omega \times (0,T)$, subject to
\begin{equation}
\begin{cases}
\partial_{\textbf{n}}\widetilde{\phi}= \partial_{\textbf{n}}\widetilde{\mu}=0,\quad\text{ on }\partial\Omega\times(0,T),\\
 {\widetilde{\phi}}(0)={\phi}_0,\quad {\widetilde{\psi}}(0)={\psi}_{0,\alpha}, \quad\text{ in }\Omega.\\
\end{cases}
\end{equation}
By \eqref{vtilde} and being $\psi_{0,\alpha}\in C^\infty(\overline{\Omega})$, in virtue of Lemma \ref{transport}, with $\widetilde{\textbf{v}}$ as velocity,
we know that there exists a unique $\widetilde{\psi}$ such that
		\begin{align}
		\widetilde{\psi}\in C([0,T];C^2(\overline{\Omega}))\cap C^1([0,T];C^1(\overline{\Omega})),
		\label{psiregtilde}
		\end{align}
which solves \eqref{sys0}$_1$ with its initial condition.

By \eqref{vtilde} and \eqref{psiregtilde} we have
	    $$f:=\frac{1}{2}\vert\nabla\widetilde{\psi}\vert^2\in L^\infty(0,T;V)\cap L^2(0,T;L^\infty(\Omega))\cap H^1(0,T;H),$$
	    and, being $\phi_0\in H^2(\Omega)$, with $\partial_\textbf{n}\phi_0=0$ almost everywhere on $\partial\Omega$, and $\mu_0:=-\Delta\phi_0+\Psi^\prime(\phi_0)\in V$, in virtue of Theorem \ref{ch} (see also \eqref{vtilde}) we deduce that there exists a unique $\widetilde{\phi}$ such that, for any $q\in[2,\infty)$,
	    \begin{align}
	    &\Vert\widetilde{\phi }\Vert_{L^\infty(\Omega\times(0,T))}\leq 1, \label{physbound}\\&
	    \widetilde{\phi}\in L^\infty(0,T;W^{2,q}(\Omega))\cap L^4(0,T;H^2(\Omega)) \cap H^1(0,T;V) ,\label{phitilde}\\&
	    \widetilde{\mu}\in L^\infty(0,T;V),\quad F^{\prime}(\widetilde{\phi})\in L^\infty(0,T;L^q(\Omega)),\quad F^{\prime\prime}(\widetilde{\phi})\in L^\infty(0,T;L^q(\Omega)),\label{potbound}
	    \end{align}
which solves \eqref{sys0}$_2$-\eqref{sys0}$_3$ with its boundary and initial conditions.

Then, on account of \eqref{vtilde}, \eqref{psiregtilde} and \eqref{phitilde}, we can apply Theorem \ref{ns} (with $\widetilde{\textbf{v}}$ as the advection velocity) and deduce that there exists a unique ${\textbf{u}}$
such that
	    \begin{align}
	\textbf{u}\in L^\infty(0,T;\textbf{V}_\sigma)\cap L^2(0,T;\textbf{W}_\sigma)\cap H^1(0,T;\textbf{H}_\sigma)
	    \label{nsfin}
	    \end{align}
which solves
\begin{equation}
\begin{cases}
\partial_t\textbf{u}+(\widetilde{\textbf{v}}\cdot\nabla)\textbf{u}+\nabla \pi-\text{div}(\nu(\widetilde{\varphi}) D\textbf{u})+(I+\alpha \textbf{A})^{-2}\textbf{P}\text{div}(\lambda_\alpha(\widetilde{\phi})\nabla\widetilde{\psi}^T\nabla\widetilde{\psi})\\+(I+\alpha \textbf{A})^{-2}\textbf{P}\text{div}(\nabla\widetilde{\phi}\otimes\nabla\widetilde{\phi})+\nu(\widetilde{\phi})\frac{(1-\widetilde{\phi})\textbf{u}}{2k(\widetilde{\phi})}=\textbf{0},
\quad\text{ in } \Omega\times(0,T)\\
\text{div}\,\textbf{u}=0,\quad\text{ in } \Omega\times(0,T)\\
\textbf{u}=0,\quad\text{ on }\partial\Omega\times(0,T),\\
\textbf{u}(0)=\textbf{u}_0,\quad\text{ in }\Omega,\\
\end{cases}
\label{sys}
\end{equation}
where
$$
	\widetilde{\varphi}=(I+\alpha A_1)^{-1}\widetilde{\phi}.
$$
Indeed the assumptions are verified. In particular, recalling \eqref{nuA} and \eqref{k}, we have
	   \begin{align*}
	   &g:=\nu(\widetilde{\phi })\frac{(1-\widetilde{\phi})}{2k(\widetilde{\phi })}\geq 0,\quad\text{a.e. in }\Omega\times(0,T),\quad \Vert g\Vert_{L^\infty(0,T;L^\infty(\Omega))}\leq C,\\&
	   \textbf{h}:=\textbf{P}\text{div}(\lambda(\widetilde{\phi})\nabla\widetilde{\psi}^T\nabla\widetilde{\psi})+\textbf{P}\text{div}(\nabla\widetilde{\phi}\otimes\nabla\widetilde{\phi}),
	   	\end{align*}
	   	and (see \eqref{phitilde})
	   	$$
	   	\Vert\textbf{P}\text{div}(\nabla\widetilde{\phi}\otimes\nabla\widetilde{\phi})\Vert_{\textbf{V}_\sigma^\prime}\leq \Vert\nabla\widetilde{\phi }\Vert_{[L^4(\Omega)]^2}^2\leq C\Vert\nabla\widetilde{\phi }\Vert\Vert\widetilde{\phi }\Vert_{H^2(\Omega)}\leq C\Vert\widetilde{\phi }\Vert_{H^2(\Omega)},
	   	$$
	   	entailing, together with the bound on $\lambda_\alpha$ and \eqref{psiregtilde},
	   	$$
	   	\textbf{h}\in L^2(0,T;\textbf{V}_\sigma^\prime).
	   	$$
We can now construct the fixed point map. Consider the closed convex set
$$
X_T:=\left\{\textbf{u}\in L^\infty(0,T;\textbf{V}_\sigma)\cap H^1(0,T;\textbf{H}_\sigma) : \Vert\textbf{u}\Vert_{L^\infty(0,T;\textbf{V}_\sigma)}^2+\Vert\partial_t\textbf{u}\Vert_{L^2(0,T;\textbf{H}_\sigma)}^2\leq M^2\right\},
$$
where $M$ will be chosen later on. Notice that Aubin-Lions Lemma yields
$$
L^\infty(0,T_0;\textbf{V}_\sigma)\cap H^1(0,T_0;\textbf{H}_\sigma)\hookrightarrow\hookrightarrow C([0,T_0];\textbf{H}_\sigma),
$$
which implies that $X_T$ is compact in $C([0,T];\textbf{H}_\sigma)$.
Therefore we can define a map ${S}: X_T \to C([0,T];\textbf{H}_\sigma)$ by setting $\textbf{u}={S}(\tilde{\textbf{u}})$ and it is clear that a fixed point of $S$ is a solution to \eqref{problem}-\eqref{problem_BIC}.

We first show that $S$ takes $X_T$ into itself for a sufficiently small $T=T_0$. On account of Lemma \ref{transport}, we can apply Lemma \ref{bourg} to \eqref{composition},
with $q=\infty$, $p\in[2,\infty)$ and $\beta=2$ to obtain
    \begin{align}
    &\nonumber\Vert\widetilde{\psi}\Vert_{[W^{2,p}(\Omega)]^2}\\
    &\nonumber\leq C\Vert\psi_{0,\alpha}\Vert_{[W^{2,p}(\Omega)]^2}(1+\Vert X(0,t,\cdot)\Vert_{W^{2,\infty}(\Omega)}^2)\\&\leq\nonumber C\Vert\psi_{0}\Vert_{[W^{2,p}(\Omega)]^2}\left(1+
 Ce^{2\int_0^t\Vert\widetilde{\textbf{v}}(\tau)\Vert_{[W^{1,\infty}(\Omega)]^2}}\left(\int_0^t\Vert\widetilde{\textbf{v}}(s)\Vert_{[W^{2,\infty}(\Omega)]^2}e^{\int_0^s\Vert\widetilde{\textbf{v}}(\tau)\Vert_{[W^{1,\infty}(\Omega)]^2}d\tau}ds) \right)^2\right)\\&\leq C\Vert\psi_{0}\Vert_{[W^{2,p}(\Omega)]^2}\left(1+
 C_\alpha e^{C_\alpha\int_0^t\Vert\widetilde{\textbf{u}}(\tau)\Vert}\left(\int_0^t\Vert\widetilde{\textbf{u}}(s)\Vert e^{C_\alpha\int_0^s\Vert\widetilde{\textbf{u}}(\tau)\Vert d\tau}ds\right)^2 \right)\nonumber\\&\leq C\Vert\psi_{0}\Vert_{[W^{2,p}(\Omega)]^2}\left(1+
 C_\alpha T_0^2 e^{C_\alpha T_0\Vert\widetilde{\textbf{u}}\Vert_{C([0,T_0];\textbf{H}_\sigma)}}\Vert\widetilde{\textbf{u}}\Vert^2_{C([0,T_0];\textbf{H}_\sigma)} \right)\nonumber\\&\leq  K_1:= C\Vert\psi_{0}\Vert_{[W^{2,p}(\Omega)]^2}\left(1+
 C_\alpha T_0^2 e^{C_\alpha T_0M}M^2\right).
    \label{psitilde}
    \end{align}
Here and in the sequel $C$ and $C_\alpha$ denote generic positive constants, the latter depending on the regularization parameter. Any other further dependence will be pointed out if necessary.
Note that we used the fact that $\vert \textbf{J}_G\vert \equiv1$ for any $t\geq0$, since $\text{div}\ \widetilde{\textbf{v}}\equiv0$. Also, we have exploited the embeddings $H^4(\Omega)\hookrightarrow W^{1,\infty}(\Omega)$, $H^4(\Omega)\hookrightarrow W^{2,\infty}(\Omega)$ and \eqref{calpha1}, together with the fact that $\psi_{0,\alpha}\rightarrow\psi_0$ in $W^{2,p}(\Omega)$ as $\alpha\rightarrow0$. In the last inequality, note that
$\Vert\textbf{u}\Vert_{C([0,T];\textbf{H}_\sigma)}\leq CM$ due to Aubin-Lions Lemma.

Let us now multiply \eqref{sys0}$_2$ by $\phi$ and integrate over $\Omega$. Integrating by parts, we get
\begin{equation}
\label{phitildeb}
\frac{d}{dt}\frac{1}{2}\Vert\widetilde{\phi}\Vert^2+\Vert\Delta\widetilde{\phi}\Vert^2+(\Psi^{\prime\prime}(\widetilde{\phi }),\vert\nabla\widetilde{\phi }\vert^2)-\frac{1}{2}(\lambda_\alpha^{\prime}(\widetilde{\phi})\vert\nabla\widetilde{\psi}\vert^2,\Delta\widetilde{\phi })=0.
\end{equation}
Recalling that $ \Psi^{\prime\prime}\geq -C$ we have
$$
(\Psi^{\prime\prime}(\widetilde{\phi }),\vert\nabla\widetilde{\phi }\vert^2) \geq - C \Vert\nabla\widetilde{\phi }\Vert^2.
$$
On the other hand, using Cauchy-Schwarz, Young's and Gagliardo-Nirenberg's inequalities, recalling $\lambda_\alpha^\prime$ bounded and \eqref{est}, we obtain
\begin{align*}
&\frac{1}{2}(\lambda_\alpha^{\prime}(\widetilde{\phi})\vert\nabla\widetilde{\psi}\vert^2,\Delta\widetilde{\phi })\leq \frac{1}{4}\Vert\Delta\widetilde{\phi}\Vert^2+C\Vert\nabla\widetilde{\psi}\Vert_{[L^4(\Omega)]^2}^4\\&\leq  \frac{1}{2}\Vert\Delta\widetilde{\phi}\Vert^2+C\Vert\widetilde{\psi}\Vert_{[H^2(\Omega)]^2}^2\Vert\widetilde{\psi }\Vert_{[L^\infty(\Omega)]^2}^2\leq \frac{1}{2}\Vert\Delta\widetilde{\phi}\Vert^2+C\Vert\widetilde{\psi}\Vert_{[H^2(\Omega)]^2}^2.
\end{align*}
Therefore, from \eqref{psitilde} we deduce
\begin{align*}
&\frac{1}{2}(\lambda_\alpha^{\prime}(\widetilde{\phi})\vert\nabla\widetilde{\psi}\vert^2,\Delta\widetilde{\phi })\\
&\leq \frac{1}{4}\Vert\Delta\widetilde{\phi}\Vert^2+C\Vert\psi_{0}\Vert_{[W^{2,p}(\Omega)]^2}^2\left(1+
C_\alpha T_0^4 e^{C_\alpha T_0\Vert\widetilde{\textbf{u}}\Vert_{C([0,T_0];\textbf{H}_\sigma)}}\Vert\widetilde{\textbf{u}}\Vert^4_{C([0,T_0];\textbf{H}_\sigma)}\right).
\end{align*}
Using also
\begin{equation}
\Vert\widetilde{\phi}\Vert_{V}^2\leq \Vert\Delta\widetilde{\phi}\Vert\Vert\widetilde{\phi}\Vert+\Vert\widetilde{\phi}\Vert^2,\label{V}
\end{equation}
from \eqref{physbound} and \eqref{phitildeb} we deduce the following inequality
\begin{equation}
\frac{d}{dt}\frac{1}{2}\Vert\widetilde{\phi}\Vert^2+\frac{1}{4}\Vert\Delta\widetilde{\phi}\Vert^2 \leq C
+ C\Vert\psi_{0}\Vert_{[W^{2,p}(\Omega)]^2}^2\left(1+
C_\alpha T_0^4 e^{C_\alpha T_0\Vert\widetilde{\textbf{u}}\Vert_{C([0,T_0];\textbf{H}_\sigma)}}\Vert\widetilde{\textbf{u}}\Vert^4_{C([0,T_0];\textbf{H}_\sigma)}\right).
\end{equation}
Thus, by Gronwall's Lemma, we obtain in the end, for any $t\in[0,T_0]$,
\begin{equation}
\sup_{t\in[0,T_0]}\Vert\widetilde{\phi}\Vert^2+\int_0^{T_0}\Vert\Delta\widetilde{\phi }\Vert^2ds
\leq K_2 := C(\Vert\phi_0\Vert^2 + 1)+C\Vert\psi_{0}\Vert_{[W^{2,p}(\Omega)]^2}^2\left(1+
C_\alpha T_0^4 e^{C_\alpha T_0M}M^4\right),
\label{phi}
\end{equation}
where we used again the fact that $\widetilde{\textbf{u}}\in X_T$.

In order to show that $S(\widetilde{\textbf{u}})\in X_T$, we multiply \eqref{sys}$_1$ by $\partial_t\textbf{u}$ and integrate over $\Omega$. This gives
\begin{align*}
&\Vert\partial_t\textbf{u}\Vert^2-(\text{div}(\nu(\widetilde{\varphi })D\textbf{u}),\partial_t\textbf{u})+((\widetilde{\textbf{v}}\cdot\nabla)\textbf{u},\partial_t\textbf{u})+(\nu(\widetilde{\phi })\frac{(1-\widetilde{\phi})}{2k(\widetilde{\phi })},\textbf{u}\cdot \partial_t\textbf{u})\\&-(\nabla\widetilde{\phi}\otimes\nabla\widetilde{\phi},\nabla(I+\alpha\textbf{A})^{-2}\partial_t\textbf{u})-(\lambda_\alpha(\widetilde{\phi})\nabla\widetilde{\psi}^T\nabla\widetilde{\psi},\nabla(I+\alpha\textbf{A})^{-2}\partial_t\textbf{u})=0.
\end{align*}
Observe that
\begin{align*}
-(\text{div}(\nu(\widetilde{\varphi })D\textbf{u}),\partial_t\textbf{u})=-(\nu^\prime(\widetilde{\varphi })D\textbf{u}\nabla\widetilde{\varphi },\partial_t\textbf{u})-\frac{1}{2}(\nu(\widetilde{\varphi })\Delta\textbf{u},\partial_t\textbf{u}).
\end{align*}
Recalling that $\nu\in W^{1,\infty}(\R)$, by \eqref{elliptic} (exploiting the embedding $H^3(\Omega)\hookrightarrow W^{1,\infty}(\Omega)$), \eqref{V} and \eqref{phi}, we deduce
\begin{align*}	
&(\nu^\prime(\widetilde{\varphi})D\textbf{u}\nabla\widetilde{\varphi },\partial_t\textbf{u})\leq C\Vert D\textbf{u}\Vert\Vert\nabla\widetilde{\varphi }\Vert_{[L^\infty(\Omega)]^2}\Vert\partial_t\textbf{u}\Vert\\&\leq C_\alpha\Vert D\textbf{u}\Vert\Vert\nabla\widetilde{\phi }\Vert\Vert\partial_t\textbf{u}\Vert\leq C_\alpha\Vert D\textbf{u}\Vert(\Vert\Delta\widetilde{\phi }\Vert+\Vert\widetilde{\phi }\Vert)\Vert\partial_t\textbf{u}\Vert\\&\leq \frac {1}{12}\Vert\partial_t\textbf{u}\Vert^2+C_\alpha(\Vert\Delta\widetilde{\phi }\Vert^2+\Vert\widetilde{\phi }\Vert^2)\Vert D\textbf{u}\Vert^{2}\\&\leq  \frac {1}{12}\Vert\partial_t\textbf{u}\Vert^2+C_\alpha(K_2+\Vert\Delta\widetilde{\phi }\Vert^2)\Vert D\textbf{u}\Vert^{2}.
\end{align*}
On the other hand, recalling \eqref{nuA}, we get
\begin{align*}
\frac{1}{2}(\nu(\widetilde{\varphi })\Delta\textbf{u},\partial_t\textbf{u})\leq C\Vert\textbf{A}\textbf{u}\Vert\Vert\partial_t\textbf{u}\Vert\leq C_1\Vert\textbf{Au}\Vert^2+\frac{1}{12}\Vert\partial_t\textbf{u}\Vert^2,
\end{align*}
for some $C_1>0$. Thus, by Young's and Poincar\'{e}'s inequalities, we infer
$$
(\nu(\widetilde{\phi })\frac{(1-\widetilde{\phi})}{2k(\widetilde{\phi })},\textbf{u}\cdot \partial_t\textbf{u})\leq C\Vert\nabla\textbf{u}\Vert^2+\frac{1}{12}\Vert\partial_t\textbf{u}\Vert^2.
$$
Moreover, by the embedding $[H^2(\Omega)]^2\hookrightarrow [L^\infty(\Omega)]^2$ and \eqref{calpha1}, recalling that $\widetilde{\textbf{u}}\in X_T$, we find
\begin{align*}
&((\widetilde{\textbf{v}}\cdot\nabla)\textbf{u},\partial_t\textbf{u})\leq C\Vert\widetilde{\textbf{v}}\Vert_{[L^\infty(\Omega)]^2}\Vert D\textbf{u}\Vert\Vert\partial_t\textbf{u}\Vert\\
&\leq C_\alpha M\Vert D\textbf{u}\Vert\Vert\partial_t\textbf{u}\Vert\leq 2C_\alpha^2M^2\Vert D\textbf{u}\Vert^2+\frac{1}{12}\Vert\partial_t\textbf{u}\Vert^2.
\end{align*}
Then, by H\"older's inequality, \eqref{calpha1}, \eqref{phi} and \eqref{V}, we get
\begin{align*}
&(\nabla\widetilde{\phi}\otimes\nabla\widetilde{\phi},\nabla(I+\alpha\textbf{A})^{-2}\partial_t\textbf{u})\leq C_\alpha\Vert\partial_t\textbf{u}\Vert\Vert\nabla\widetilde{\phi}\Vert^2\\
&\leq C_\alpha\Vert\partial_t\textbf{u}\Vert\left(\Vert\Delta\widetilde{\phi}\Vert\Vert\widetilde{\phi}\Vert+\Vert\widetilde{\phi}\Vert^2\right)\\
&\leq C_\alpha\Vert\partial_t\textbf{u}\Vert(K_2+K_2^{1/2}\Vert\Delta\widetilde{\phi}\Vert)\leq  \frac{1}{12}\Vert\partial_t\textbf{u}\Vert^2+C_\alpha(K_2^2+K_2\Vert\Delta\widetilde{\phi}\Vert^2).
\end{align*}
Using similar arguments, by Gagliardo-Nirenberg's inequality, recalling \eqref{est} and by \eqref{psitilde}, being $\lambda_\alpha$ bounded, we also have
\begin{align*}
&(\lambda_\alpha(\widetilde{\phi})\nabla\widetilde{\psi}^T\nabla\widetilde{\psi},\nabla(I+\alpha\textbf{A})^{-2}\partial_t\textbf{u})\leq C_\alpha\Vert\partial_t\textbf{u}\Vert\Vert\nabla\psi\Vert_{[L^4(\Omega)]^4}^2\\
&\leq C_\alpha\Vert\partial_t\textbf{u}\Vert\Vert\psi\Vert_{[L^\infty(\Omega)]^2}\Vert\psi\Vert_{[H^2(\Omega)]^2} \leq \frac{1}{12}\Vert\partial_t\textbf{u}\Vert^2+C_\alpha K_1^4.
\end{align*}
To sum up, taking the above estimates into account, we have
\begin{align}
\nonumber
\frac{1}{2}\Vert\partial_t\textbf{u}\Vert^2 &\leq C_2\Vert \textbf{A}\textbf{u}\Vert^2\\
&+C_\alpha(M^2+K_2+\Vert\Delta\widetilde{\phi }\Vert^2)\Vert \nabla\textbf{u}\Vert^2+C_\alpha(1+K_2^2+K_2\Vert\Delta\widetilde{\phi}\Vert^2)+C_\alpha K_1^4,
\label{dtn2}
\end{align}
We then multiply \eqref{sys}$_1$ by  $\textbf{A}\textbf{u}$ and integrate over $\Omega$. Recalling that $(I+\alpha^2\textbf{A})^{-2}$ is self-adjoint, we obtain
\begin{align}
&\nonumber(\partial_t\textbf{u},\textbf{Au})-(\text{div}(\nu(\widetilde{\varphi})D\textbf{u}),\textbf{A}\textbf{u})+((\textbf{v}\cdot \nabla)\textbf{u},\textbf{Au})+(\nu(\widetilde{\phi })\frac{(1-\widetilde{\phi})}{2k(\widetilde{\phi })},\textbf{Au})\\&=(\nabla\widetilde{\phi}\otimes\nabla\widetilde{\phi},\nabla(I+\alpha^2\textbf{A})^{-2}\textbf{Au})+(\lambda_\alpha(\widetilde{\phi})\nabla\widetilde{\psi}^T\nabla\widetilde{\psi},\nabla(I+\alpha^2\textbf{A})^{-2}\textbf{Au}).
\label{Au1}
\end{align}
We know that there exists $\tilde{p}\in L^2(0,T;V)$ such that
$$
-\Delta\textbf{u}+\nabla \tilde{p}=\textbf{A}\textbf{u},
$$
almost everywhere in $\Omega\times(0,T)$,
with
\begin{align}
\Vert \tilde{p}\Vert\leq C\Vert\nabla\textbf{u}\Vert^{1/2}\Vert\textbf{A}\textbf{u}\Vert^{1/2},\qquad \Vert \tilde{p}\Vert_V\leq C\Vert\textbf{A}\textbf{u}\Vert.
\label{pn1}
\end{align}
Arguing as in the proof of Theorem \ref{ns} (see Subsection \ref{proof2} below), we have
\begin{align*}
&-(\text{div}(\nu(\widetilde{\varphi})D\textbf{u}),\textbf{A}\textbf{u})\geq-(\nu^\prime(\widetilde{\varphi})D\textbf{u}\nabla\tilde{\varphi},\textbf{A}\textbf{u})+\frac{\nu_*}{28}\Vert\textbf{Au}\Vert^2+\frac{1}{2}(\nu^\prime(\widetilde{\varphi})\nabla\widetilde{\varphi} \tilde{p},\textbf{Au}).
\end{align*}
Then, by \eqref{elliptic} and \eqref{phi}, being $\nu\in W^{1,\infty}(\R)$,
\begin{align*}
&(\nu^\prime(\widetilde{\varphi})D\textbf{u}\nabla\widetilde{\varphi},\textbf{A}\textbf{u})\leq  C\Vert D\textbf{u}\Vert\Vert\nabla\widetilde{\varphi }\Vert_{[L^\infty(\Omega)]^2}\Vert\textbf{A}\textbf{u}\Vert\\&\leq C_\alpha\Vert D\textbf{u}\Vert\Vert\nabla\widetilde{\phi }\Vert\Vert\textbf{A}\textbf{u}\Vert\leq C_\alpha\Vert D\textbf{u}\Vert(\Vert\Delta\widetilde{\phi }\Vert+\Vert\widetilde{\phi }\Vert)\Vert\textbf{A}\textbf{u}\Vert\\&\leq \frac {\nu_*}{28}\Vert\textbf{A}\textbf{u}\Vert^2+C_\alpha(\Vert\Delta\widetilde{\phi }\Vert^2+\Vert\widetilde{\phi }\Vert^2)\Vert D\textbf{u}\Vert^{2}\\&\leq  \frac {\nu_*}{28}\Vert\textbf{A}\textbf{u}\Vert^2+C_\alpha(K_2+\Vert\Delta\widetilde{\phi }\Vert^2)\Vert D\textbf{u}\Vert^{2},		
\end{align*}
so that, thanks to \eqref{V}, \eqref{phi} and \eqref{pn1}, we infer
\begin{align*}
&\frac{1}{2}(\nu^\prime(\widetilde{\varphi})\nabla\widetilde{\varphi}\tilde{p},\textbf{Au})\leq C\Vert\nabla\widetilde{\varphi}\Vert_{[L^\infty(\Omega)]^2}\Vert \tilde{p}\Vert\Vert\textbf{Au}\Vert\\&\leq C_\alpha\Vert\nabla\widetilde{\phi}\Vert\Vert \tilde{p}\Vert\Vert\textbf{Au}\Vert\leq C_\alpha\Vert\nabla\widetilde{\phi}\Vert\Vert \nabla\textbf{u}\Vert^{1/2}\Vert \textbf{Au}\Vert^{3/2} \leq \frac{\nu_*}{28}\Vert\textbf{A}\textbf{u}\Vert^{2}+C\Vert\nabla\widetilde{\phi}\Vert^4\Vert D\textbf{u}\Vert^2\\&\leq \frac{\nu_*}{28}\Vert\textbf{A}\textbf{u}\Vert^{2}+C(\Vert\widetilde{\phi}\Vert^2\Vert\Delta\widetilde{\phi}\Vert^2+\Vert\widetilde{\phi}\Vert^4)\Vert D\textbf{u}\Vert^2\leq \frac{\nu_*}{28}\Vert\textbf{A}\textbf{u}\Vert^{2}+C(K_2\Vert\Delta\widetilde{\phi}\Vert^2+K_2^2)\Vert D\textbf{u}\Vert^2.
\end{align*}
Observe now that
$$
((\widetilde{\textbf{v}}\cdot \nabla)\textbf{u},\textbf{Au})\leq C\Vert\widetilde{\textbf{v}}\Vert_{L^\infty(\Omega)}\Vert D\textbf{u}\Vert\Vert\textbf{Au}\Vert\leq \frac{\nu_*}{28}\Vert\textbf{Au}\Vert^2+C\Vert\widetilde{\textbf{v}}\Vert_{L^\infty(\Omega)}^2\Vert D\textbf{u}\Vert^2.
$$
Moreover, by the regularity properties of the operator $(I+\alpha^2\textbf{A})^{-2}$ (see \eqref{calpha1}) and \eqref{phi}, we find
\begin{align*}
&(\nabla\widetilde{\phi}\otimes\nabla\widetilde{\phi},\nabla(I+\alpha^2\textbf{A})^{-2}\textbf{Au})\leq \Vert\nabla\widetilde{\phi}\Vert^2\Vert\nabla(I+\alpha^2\textbf{A})^{-2}\textbf{Au}\Vert_{[L^\infty(\Omega)]^4}\\&\leq C_\alpha\Vert\nabla\widetilde{\phi}\Vert^2\Vert\textbf{Au}\Vert\leq C_\alpha(\Vert\widetilde{\phi}\Vert\Vert\Delta\widetilde{\phi}\Vert+\Vert\widetilde{\phi}\Vert^2)\Vert\textbf{Au}\Vert\\&\leq  \frac{\nu_*}{28}\Vert\textbf{Au}\Vert^2+C_\alpha(K_2\Vert\Delta\widetilde{\phi}\Vert^2+K_2^2).
\end{align*}
Arguing similarly and using \eqref{psitilde}, we deduce
\begin{align*}
&(\lambda_\alpha(\widetilde{\phi})\nabla\widetilde{\psi}^T\nabla\widetilde{\psi},\nabla(I+\alpha^2\textbf{A})^{-2}\textbf{Au})\\
&\leq \Vert\lambda_\alpha(\widetilde{\phi})\nabla\widetilde{\psi}^T\nabla\widetilde{\psi}\Vert_{[L^1(\Omega)]^4}\Vert\nabla(I+\alpha^2\textbf{A})^{-2}\textbf{Au}\Vert_{[L^\infty(\Omega)]^4}\\&\leq C_\alpha\Vert\nabla\widetilde{\psi}\Vert^2\Vert\textbf{Au}\Vert\leq C_\alpha K_1^2\Vert\textbf{Au}\Vert\leq  \frac{\nu_*}{28}\Vert\textbf{Au}\Vert^2+C_\alpha K_1^4.
\end{align*}
Finally, Poincar\'{e}'s inequality yields
$$
(\nu(\widetilde{\phi })\frac{(1-\widetilde{\phi})}{2k(\widetilde{\phi })},\textbf{u}\cdot \textbf{A}\textbf{u})\leq C\Vert\nabla\textbf{u}\Vert^2+\frac{\nu_*}{28}\Vert\textbf{A}\textbf{u}\Vert^2,
$$

Summing up, recalling \eqref{calpha1} and $\widetilde{\textbf{u}}\in X_T$, on account of the above estimates, we infer from \eqref{Au1} that
\begin{align*}
&\frac{1}{2}\frac{d}{dt}\Vert\nabla\textbf{u}\Vert^2+\frac{\nu_*}{4}\Vert\textbf{Au}\Vert^2\\&\leq C_\alpha\left(1+K_2+(1+K_2)\Vert\Delta\widetilde{\phi}\Vert^2+K_2^2+\Vert\widetilde{\textbf{v}}\Vert_{L^\infty(\Omega)}^2\right)\Vert\nabla\textbf{u}\Vert^2+C_\alpha( K_1^4+K_2\Vert\Delta\widetilde{\phi}\Vert^2+K_2^2)\\&\leq
C_\alpha\left(1+K_2+(1+K_2)\Vert\Delta\widetilde{\phi}\Vert^2+K_2^2+M^2\right)\Vert\nabla\textbf{u}\Vert^2+C_\alpha( K_1^4+K_2\Vert\Delta\widetilde{\phi}\Vert^2+K_2^2) ,
\end{align*}
and adding this estimate and \eqref{dtn2} multiplied by $\omega=\frac{\nu_*}{8C_2}$ together give
\begin{align}
\nonumber&\frac{1}{2}\frac{d}{dt}\Vert\nabla\textbf{u}\Vert^2+\frac{\omega}{2}\Vert\partial_t\textbf{u}\Vert^2+\frac{\nu_*}{8}\Vert\textbf{Au}\Vert^2\\&\leq C_\alpha\left(1+K_2+(1+K_2)\Vert\Delta\widetilde{\phi}\Vert^2+K_2^2+M^2\right)\Vert\nabla\textbf{u}\Vert^2\nonumber\\&+C_\alpha(K_2\Vert\Delta\widetilde{\phi}\Vert^2+K_2^2+K_1^4).
\label{G}
\end{align}
Thus, exploiting \eqref{phi} and recalling that $\textbf{u}\in C([0,T];\textbf{V}_\sigma)$, an application of Gronwall's Lemma yields
\begin{align}
\nonumber\sup_{t\in[0,T_0]}\Vert\nabla\textbf{u}\Vert^2&\leq \Vert\nabla\textbf{u}_{0,\alpha}\Vert^2e^{C_\alpha T_0\left(1+K_2+K_2^2+M^2\right)+C_\alpha(1+K_2)K_2}\\&+e^{C_\alpha T_0\left(1+K_2+K_2^2+M^2\right)+(1+K_2)K_2}\left(C_\alpha T_0(K_2^2+K_1^4)+C_\alpha K_2^2\right).
\label{Gr}
\end{align}
We now observe that $K_1$ and $K_2$ can be rewritten as
$$
K_1=C\Vert\psi_0\Vert_{[W^{2,p}(\Omega)]^2}+\mathcal{Q}_{\alpha,1}(T_0)T_0,\quad K_2=C(\Vert\phi_0\Vert^2 + \Vert\psi_0\Vert^2_{[W^{2,p}(\Omega)]^2} + 1) +\mathcal{Q}_{\alpha,2}(T_0)T_0,
$$
 where from now on we consider $\mathcal{Q}_{\alpha,i}(T_0)$ as a generic function depending on the initial data, $\alpha$ and $M$, nonincreasing in $T_0$, i.e., such that $\lim_{T_0\to0}\mathcal{Q}_{\alpha,i}(T_0)T_0=0$. Therefore, \eqref{Gr} can be rewritten as
$$
\sup_{t\in[0,T_0]}\Vert\nabla\textbf{u}\Vert^2 \leq
e^{C_\alpha (\phi_0,\psi_0) +\mathcal{Q}_{\alpha,3}(T_0)T_0}\left( \Vert\nabla\textbf{u}_{0,\alpha}\Vert^2+ C_\alpha(\phi_0,\psi_0){\color{blue}+} \mathcal{Q}_{\alpha,4}(T_0)T_0\right),
$$
where $C_\alpha(\phi_0,\psi_0) := C_\alpha(1 + \Vert\phi_0\Vert^4 + \Vert\psi_0\Vert^4_{[W^{2,p}(\Omega)]^2})$.
 Note that we have used $\Vert\nabla\textbf{u}_{0,\alpha}\Vert^2=R^2$ and that all the terms depending on $M$ are encapsulated in some $\mathcal{Q}_{\alpha,i}(T_0)$.
We now set $T_0=T_0(M)$ such that
\begin{align}
\begin{cases}
\mathcal{Q}_{\alpha,3}(T_0)T_0\leq C_\alpha(\phi_0,\psi_0)  ,\\
\mathcal{Q}_{\alpha,4}(T_0)T_0 \leq C_\alpha(\phi_0,\psi_0) ,
\end{cases}
\label{T0}
\end{align}
so that we obtain
$$
\sup_{t\in[0,T_0]}\Vert\nabla\textbf{u}\Vert^2 \leq
e^{2C_\alpha(\Vert\phi_0\Vert^2+\Vert\phi_0\Vert^4)}\left( \Vert\nabla\textbf{u}_{0,\alpha}\Vert^2+ 2C_\alpha(\phi_0,\psi_0)\right):=\widetilde{\mathcal{Q}}_\alpha,
$$
where $\widetilde{\mathcal{Q}}_\alpha$ does not depend neither on $M$ nor on $T_0$.
Then, by integrating \eqref{G} over $(0,T_0)$, we also obtain
\begin{align*}
&\int_0^{T_0}\Vert\partial_t\textbf{u}\Vert^2 dt\\
&\leq \max\left\{\frac{2}{\omega},1\right\}R^2+\left(C_\alpha T_0\left(1+K_2+K_2^2+M^2\right)+C_\alpha(1+K_2)K_2\right)\sup_{t\in[0,T_0]}\Vert\nabla\textbf{u}\Vert^2\\&+C_\alpha T_0(K_2^2+K_1^4)+C_\alpha K_2^2\leq \max\left\{\frac{2}{\omega},1\right\}R^2+\mathcal{Q}_{\alpha,5}(T_0)T_0+C_\alpha K_2(1+K_2)\widetilde{\mathcal{Q}}_\alpha+C_\alpha K_2^2\\&\leq
\max\left\{\frac{2}{\omega},1\right\}R^2+C_\alpha(\phi_0,\psi_0)(1 + \widetilde{\mathcal{Q}}_\alpha)+\mathcal{Q}_{\alpha,6}(T_0)T_0.
\end{align*}
Therefore in the end we get
\begin{align*}
\sup_{t\in[0,T_0]}\Vert\nabla\textbf{u}\Vert^2+\int_0^{T_0}\Vert\partial_t\textbf{u}\Vert^2 dt
\leq \max\left\{\frac{2}{\omega},1\right\}R^2+ C_\alpha(\phi_0,\psi_0)(1 + \widetilde{\mathcal{Q}}_\alpha) +\mathcal{Q}_{\alpha,6}(T_0)T_0.
\end{align*}
 If we now fix $T_0=T_0(M)$ sufficiently small such that \eqref{T0} and
\begin{align}
\mathcal{Q}_{\alpha,6}(T_0)T_0\leq C_\alpha(\phi_0,\psi_0)
\end{align}
hold, then we get
\begin{align*}
\sup_{t\in[0,T_0]}\Vert\nabla\textbf{u}\Vert^2+\int_0^{T_0}\Vert\partial_t\textbf{u}\Vert^2 dt
\leq \max\left\{\frac{2}{\omega},1\right\}R^2+2C_\alpha(\phi_0,\psi_0) +  C_\alpha(\phi_0,\psi_0)\widetilde{\mathcal{Q}}_\alpha:=M^2.
\end{align*}
Thus, for this choice of $M$ and, consequently, for the choice of $T_0=T_0(M)$, we have
\begin{equation}
\Vert\partial_t\textbf{u}\Vert_{L^2(0,T_0;\textbf{H}_\sigma)}^2+\Vert \textbf{u}\Vert_{L^\infty(0,T_0;\textbf{V}_\sigma)}^2\leq M^2.
\label{est3}
\end{equation}
This implies that $S:X_{T_0}\to X_{T_0}$ is well defined.

We are left to show that the map $S$ is continuous from $X_{T_0}$ to itself with respect to the topology of $C([0,T_0];\textbf{H}_\sigma)$.
Let us define $\widetilde{\textbf{u}}=\widetilde{\textbf{u}}_1-\widetilde{\textbf{u}}_2$ and $\textbf{u}=S(\widetilde{\textbf{u}}_1)-S(\widetilde{\textbf{u}}_2)$. As usual, we start from $\widetilde{\psi}=\widetilde{\psi}_1-\widetilde{\psi}_2$ which solves the following
$$
\partial_t\widetilde{\psi}+\widetilde{\textbf{v}}_1\cdot \nabla\widetilde{\psi}+\widetilde{\textbf{v}}\cdot \nabla\widetilde{\psi}_2=0.
$$
Recalling \eqref{psiregtilde}, which means in particular that $\widetilde{\psi}\in AC([0,T];[W^{1,p}(\Omega)]^2)$,
we take the gradient of the above equation, and then multiply it (double dot product) by $\vert\nabla\widetilde{\psi}\vert^{p-2}\nabla\widetilde{\psi}$. This gives
\begin{align}
\frac{1}{p}\frac{d}{dt}\Vert\nabla\widetilde{\psi}\Vert_{L^p(\Omega)}^p+( \nabla(\widetilde{\textbf{v}}_1\cdot \nabla\widetilde{\psi}), \vert\nabla\widetilde{\psi}\vert^{p-2}\nabla\widetilde{\psi})+( \nabla(\widetilde{\textbf{v}}\cdot \nabla\widetilde{\psi}_2), \vert\nabla\widetilde{\psi}\vert^{p-2}\nabla\widetilde{\psi})=0.
\label{diff}
\end{align}
For the first term, using \eqref{calpha1} and recalling that $\widetilde{\textbf{v}}_1$ is divergence free, we have
\begin{align*}
&( \nabla(\widetilde{\textbf{v}}_1\cdot \nabla\widetilde{\psi}), \vert\nabla\widetilde{\psi}\vert^{p-2}\nabla\widetilde{\psi})=(  \nabla\widetilde{\psi}\nabla\widetilde{\textbf{v}}_1, \vert\nabla\widetilde{\psi}\vert^{p-2}\nabla\widetilde{\psi})\\&\leq \Vert\nabla\widetilde{\textbf{v}}_1\Vert_{[W^{1,\infty}(\Omega)]^2}
\Vert\nabla\widetilde{\psi}\Vert^p_{L^p(\Omega)}\leq C_\alpha M\Vert\nabla\widetilde{\psi}\Vert^{p}_{L^p(\Omega)},
\end{align*}
where $C_\alpha$ also depends on $T_0,R,\psi_0,M$.
Concerning the second one, we can argue in a similar way, thanks to H\"older's inequality and \eqref{calpha1}. This gives
\begin{align*}
&( \nabla(\widetilde{\textbf{v}}\cdot \nabla\widetilde{\psi}_2), \vert\nabla\widetilde{\psi}\vert^{p-2}\nabla\widetilde{\psi})\\&\leq \Vert\nabla\widetilde{\textbf{v}}\Vert_{[W^{1,\infty}(\Omega)]^2}
\Vert\widetilde{\psi}_2\Vert_{[W^{1,p}(\Omega)]^2}\Vert\nabla\widetilde{\psi}\Vert_{L^p(\Omega)}^{p-1}+\Vert\widetilde{\textbf{v}}\Vert_{[L^\infty(\overline{\Omega})]^2}\Vert\widetilde{\psi}_2\Vert_{[W^{2,p}(\Omega)]^2}\Vert\nabla\widetilde{\psi}\Vert_{L^p(\Omega)}^{p-1}\\&\leq C_\alpha\Vert\widetilde{\textbf{u}}\Vert\Vert\nabla\widetilde{\psi}\Vert_{L^p(\Omega)}^{p-1}.
\end{align*}
Here we have exploited the $W^{2,p}$-control on $\widetilde{\psi}_2$ (see \eqref{psitilde}). Therefore $C_\alpha$ also depends on $T_0,R,\psi_0,M$.

Applying Gronwall's Lemma to \eqref{diff} yields
 \begin{align}
\sup_{t\in[0,T_0]}\Vert\nabla\widetilde{\psi}(t)\Vert_{[L^p(\Omega)]^4}\leq C_\alpha e^{C_\alpha MT_0}\int_0^{T_0}\Vert\widetilde{\textbf{u}}(s)\Vert ds
 \leq C_\alpha\int_0^{T_0}\Vert\widetilde{\textbf{u}}(s)\Vert ds.\label{Gronwallpsi2}
 \end{align}
We can now consider \eqref{sys0}$_2$ and \eqref{sys0}$_3$. Setting $\widetilde{\phi}=\widetilde{\phi }_1-\widetilde{\phi }_2$ and $\widetilde{\mu}=\widetilde{\mu}_1-\widetilde{\mu}_2$ we have
$$
\frac{\partial\widetilde{\phi}}{\partial t}+\widetilde{\textbf{v}}_1\cdot \nabla\widetilde{\phi}+\widetilde{\textbf{v}}\cdot \nabla\widetilde{\phi}_2-\Delta\widetilde{\mu}=0.
$$
Multiplying it by $\widetilde{\phi }$ and integrating over $\Omega$, we get
\begin{align*}
&\frac{d}{dt}\frac{1}{2}\Vert\widetilde{\phi}\Vert^2+\Vert\Delta\widetilde{\phi }\Vert^2 +(\widetilde{\textbf{v}}\cdot \nabla\widetilde{\phi}_2,\widetilde{\phi})-(\Psi^\prime(\widetilde{\phi}_1)-\Psi^\prime(\widetilde{\phi }_2),\Delta\widetilde{\phi }) \\&-\left(\lambda_\alpha^\prime(\widetilde{\phi }_1)-\lambda_\alpha^\prime(\widetilde{\phi }_2),\frac{1}{2}\vert\nabla\widetilde{\psi}_1\vert^2\Delta\widetilde{\phi }\right)-\left(\frac{\lambda_\alpha^\prime(\widetilde{\phi }_2)}{2}\left(\vert\nabla\widetilde{\psi}_1\vert^2-\vert\nabla\widetilde{\psi}_2\vert^2\right),\Delta\widetilde{\phi}\right)=0.
\end{align*}	
We have, by Young's inequality together with \eqref{calpha1}, \eqref{V} and \eqref{phi},
\begin{align*}
(\widetilde{\textbf{v}}\cdot \nabla\widetilde{\phi}_2,\widetilde{\phi})&\leq \Vert\widetilde{\textbf{v}}\Vert_{[L^\infty(\Omega)]^2}\Vert\nabla\widetilde{\phi}_2\Vert\Vert\widetilde{\phi}\Vert\leq C_\alpha\Vert\widetilde{\textbf{u}}\Vert^2+C_\alpha\left(\Vert\Delta\widetilde{\phi}_2\Vert\Vert\widetilde{\phi}_2\Vert+\Vert\widetilde{\phi}_2\Vert^2\right)\Vert\widetilde{\phi}\Vert^2\\&\leq C_\alpha\Vert\widetilde{\phi }\Vert^2\left(\Vert\Delta\widetilde{\phi}_2\Vert+1\right)+C_\alpha\Vert\widetilde{\textbf{u}}\Vert^2.
\end{align*}
Then recalling that $F^{\prime\prime}$ is convex and (see \eqref{potbound}) $\Vert F^{\prime\prime}(\widetilde{\phi_i})\Vert_{L^\infty(0,T_0;L^p(\Omega))}\leq C$, for any $p\in[2,\infty)$, we obtain
(see also \eqref{V})
\begin{align*}
&(\Psi^\prime(\widetilde{\phi}_1)-\Psi^\prime(\widetilde{\phi }_2),\Delta\widetilde{\phi })=\left(\int_0^1\Psi^{\prime\prime}(s\widetilde{\phi }_1+(1-s)\widetilde{\phi}_2)\widetilde{\phi }ds,\Delta\widetilde{\phi }\right)\\&\leq \left(\int_0^1 (s\Psi^{\prime\prime}(\widetilde{\phi }_1)+(1-s)\Psi^{\prime\prime}(\widetilde{\phi}_2))\widetilde{\phi }ds,\Delta\widetilde{\phi }\right)\\&\leq C\left(\Vert\Psi^{\prime\prime}(\widetilde{\phi }_1)\Vert_{L^4(\Omega)}+\Vert\Psi^{\prime\prime}(\widetilde{\phi }_2)\Vert_{L^4(\Omega)}\right)\Vert\widetilde{\phi }\Vert_{L^4(\Omega)}\Vert\Delta\widetilde{\phi}\Vert\\&\leq C\Vert\nabla\widetilde{\phi }\Vert\Vert\Delta\widetilde{\phi}\Vert\leq C(\Vert\widetilde{\phi}\Vert^{1/2}\Vert\Delta\widetilde{\phi}\Vert^{1/2}+\Vert\widetilde{\phi}\Vert)\Vert\Delta\widetilde{\phi}\Vert\\&\leq \frac{1}{6}\Vert\Delta\widetilde{\phi }\Vert^2+C_\alpha\Vert\widetilde{\phi }\Vert^2.
\end{align*}
On the other hand, since $\lambda_\alpha^\prime$ is globally Lipschitz,	using Young's inequality, \eqref{psitilde} and Agmon's inequality, we deduce that
\begin{align*}
&\left(\lambda_\alpha^\prime(\widetilde{\phi }_1)-\lambda_\alpha^\prime(\widetilde{\phi }_2),\frac{1}{2}\vert\nabla\widetilde{\psi}_1\vert^2\Delta\widetilde{\phi }\right)= \left(\int_0^1\lambda_\alpha^{\prime\prime}(s\widetilde{\psi }_1+(1-s)\widetilde{\psi}_2))\widetilde{\phi }ds,\frac{1}{2}\vert\nabla\widetilde{\psi}_1\vert^2\Delta\widetilde{\phi }\right)\\&\leq C\Vert\widetilde{\phi }\Vert_{L^\infty(\Omega)}\Vert\nabla\widetilde{\psi}_1\Vert^2_{L^4(\Omega)}\Vert\Delta\widetilde{\phi }\Vert\leq C_\alpha\Vert\widetilde{\phi}\Vert^{1/2}\left(\Vert\widetilde{\phi}\Vert^{1/2}+\Vert\Delta\widetilde{\phi }\Vert^{1/2}\right)\Vert\Delta\widetilde{\phi }\Vert\\
&\leq \frac{1}{6}\Vert\Delta\widetilde{\phi }\Vert^2+C_\alpha\Vert\widetilde{\phi }\Vert^2,
\end{align*}
where we used the norm $\Vert\cdot\Vert^2_{H^2(\Omega)}=\Vert\cdot\Vert^2+\Vert\Delta\cdot\Vert^2$.
Note that $C_\alpha$ also depends on $T_0,R,\psi_0,\phi_0,M$. In conclusion, being $\lambda_\alpha^\prime$ bounded, thanks to Young's inequality, \eqref{psitilde} and \eqref{Gronwallpsi2}, we infer
\begin{align*}
&\left(\frac{\lambda_\alpha^\prime(\widetilde{\phi }_2)}{2}\left(\vert\nabla\widetilde{\psi}_1\vert^2-\vert\nabla\widetilde{\psi}_2\vert^2\right),\Delta\widetilde{\phi}\right)=\left(\frac{\lambda_\alpha^\prime(\widetilde{\phi }_2)}{2}\left(\vert\nabla\widetilde{\psi}_1\vert-\vert\nabla\widetilde{\psi}_2\vert\right)\left(\vert\nabla\widetilde{\psi}_1\vert+\vert\nabla\widetilde{\psi}_2\vert\right),\Delta\widetilde{\phi}\right)\\&\leq C_\alpha\Vert\nabla\widetilde{\psi }\Vert_{[L^4(\Omega)]^4}\left(\Vert\nabla\widetilde{\psi }_1\Vert_{[L^4(\Omega)]^4}+\Vert\nabla\widetilde{\psi }_2\Vert_{[L^4(\Omega)]^4}\right)\Vert\Delta\widetilde{\phi}\Vert\leq C_\alpha\Vert\nabla\widetilde{\psi }\Vert_{[L^4(\Omega)]^4}^2+\frac{1}{6}\Vert\Delta\widetilde{\phi}\Vert^2\\&\leq C_\alpha\left(\int_0^{T_0}\Vert\widetilde{\textbf{u}}(s)\Vert ds\right)^2+\frac{1}{6}\Vert\Delta\widetilde{\phi}\Vert^2\leq C_\alpha T_0^2\Vert\widetilde{\textbf{u}}\Vert^2_{C([0,T_0];\textbf{H}_\sigma)}+\frac{1}{6}\Vert\Delta\widetilde{\phi}\Vert^2.
\end{align*}
Again $C_\alpha$ also depends on $T_0,R,\psi_0,\phi_0,M$.

Summing up, we have
$$
\frac{d}{dt}\frac{1}{2}\Vert\widetilde{\phi}\Vert^2+\frac {1}{2}\Vert\Delta\widetilde{\phi }\Vert^2\leq C_\alpha\left(\Vert\Delta\widetilde{\phi}_2\Vert+1\right)\Vert\widetilde{\phi }\Vert^2+ C_\alpha T_0^2\Vert\widetilde{\textbf{u}}\Vert^2_{C([0,T_0];\textbf{H}_\sigma)}.
$$
Thus, recalling that $\widetilde{\phi}(0)=0$ and, by \eqref{phi}, that $\int_0^{T_0}\Vert\Delta\widetilde{\phi }_2(s)\Vert ds\leq C_\alpha$, Gronwall's Lemma entails
\begin{align}
\sup_{t\in[0,T_0]}\Vert\widetilde{\phi}(t)\Vert^2+\int_0^{T_0}\Vert\Delta\widetilde{\phi }(s)\Vert^2ds\leq C_\alpha T_0^3e^{CT_0}\Vert\widetilde{\textbf{u}}\Vert_{C([0,T_0];\textbf{H}_\sigma)}^2.
\label{phidiff}
\end{align}	
Let us now consider \eqref{sys}$_1$-\eqref{sys}$_2$. Observe that ${\textbf{u}}:=\textbf{u}_1-\textbf{u}_2$ solves
\begin{align*}
&\partial_t\textbf{u}+(\widetilde{\textbf{v}}_1\cdot\nabla)\textbf{u}+(\widetilde{\textbf{v}}\cdot\nabla)\textbf{u}_2+\nabla p-\text{div}(\nu(\widetilde{\varphi}_1) D\textbf{u})-\text{div}((\nu(\widetilde{\varphi}_1)-\nu(\widetilde{\varphi}_2)) D\textbf{u}_2)\\&+(I+\alpha \textbf{A})^{-2}\textbf{P}\text{div}(\lambda_\alpha(\widetilde{\phi}_1)\nabla\widetilde{\psi}_1^T\nabla\widetilde{\psi})+(I+\alpha \textbf{A})^{-2}\textbf{P}\text{div}((\lambda_\alpha(\widetilde{\phi}_1)-\lambda_\alpha(\widetilde{\phi}_2))\nabla\widetilde{\psi}_1^T\nabla\widetilde{\psi}_2)\\&+(I+\alpha \textbf{A})^{-2}\textbf{P}\text{div}(\lambda_\alpha(\widetilde{\phi}_2)\nabla\widetilde{\psi}^T\nabla\widetilde{\psi}_2)+(I+\alpha \textbf{A})^{-2}\textbf{P}\text{div}(\nabla\widetilde{\phi}_1\otimes\nabla\widetilde{\phi})\\&+(I+\alpha \textbf{A})^{-2}\textbf{P}\text{div}(\nabla\widetilde{\phi}\otimes\nabla\widetilde{\phi}_2)-\nu(\widetilde{\phi }_2)\frac{1-\widetilde{\phi}_2}{2k(\widetilde{\phi}_1)}\textbf{u}-\nu(\widetilde{\phi}_1)\frac{\left(\widetilde{\phi}_1-\widetilde{\phi}_2\right)\textbf{u}_1}{2k(\widetilde{\phi }_1)}\\&-\left(\nu(\widetilde{\phi }_2)-\nu(\widetilde{\phi }_2)\right)\frac{1-\widetilde{\phi}_2}{2k(\widetilde{\phi}_1)}\textbf{u}_1-\nu(\widetilde{\phi }_2)(1-\widetilde{\phi}_2)\frac{k(\widetilde{\phi }_2)-k(\widetilde{\phi }_1)}{2k(\widetilde{\phi}_1)k(\widetilde{\phi}_2)}\textbf{u}_2=\textbf{0}.
\end{align*}
We now multiply by $\textbf{u}$ and integrate over $\Omega$, obtaining, setting $\textbf{v}:=(I+\alpha \textbf{A})^{-2}\textbf{u}$,
\begin{align*}
&\frac{d}{dt}\frac{1}{2}\Vert\textbf{u}\Vert^2+\left((\widetilde{\textbf{v}}\cdot\nabla)\textbf{u}_2,\textbf{u}\right)+\left(\nu(\widetilde{\varphi}_1) D\textbf{u},D\textbf{u}\right)+((\nu(\widetilde{\varphi}_1)-\nu(\widetilde{\varphi}_2)) D\textbf{u}_2,D\textbf{u})\\&-\left(\lambda_\alpha(\widetilde{\phi}_1)\nabla\widetilde{\psi}_1^T\nabla\widetilde{\psi},\nabla\textbf{v}\right)-\left((\lambda_\alpha(\widetilde{\phi}_1)-\lambda_\alpha(\widetilde{\phi}_2))\nabla\widetilde{\psi}_1^T\nabla\widetilde{\psi}_2,\nabla\textbf{v}\right)\\&-\left(\lambda_\alpha(\widetilde{\phi}_2)\nabla\widetilde{\psi}^T\nabla\widetilde{\psi}_2,\nabla\textbf{v}\right)-\left(\nabla\widetilde{\phi}_1\otimes\nabla\widetilde{\phi},\nabla\textbf{v}\right)\\&-\left(\nabla\widetilde{\phi}\otimes\nabla\widetilde{\phi}_2,\nabla\textbf{v}\right)-\left(\nu(\widetilde{\phi }_2)\frac{1-\widetilde{\phi}_2}{2k(\widetilde{\phi}_1)}\textbf{u},\textbf{u}\right)-\left(\nu(\widetilde{\phi}_1)\frac{\left(\widetilde{\phi}_1-\widetilde{\phi}_2\right)\textbf{u}_1}{2k(\widetilde{\phi }_1)},\textbf{u}\right)\\&-\left(\left(\nu(\widetilde{\phi }_2)-\nu(\widetilde{\phi }_2)\right)\frac{1-\widetilde{\phi}_2}{2k(\widetilde{\phi}_1)}\textbf{u}_1,\textbf{u}\right)-\left(\nu(\widetilde{\phi }_2)(1-\widetilde{\phi}_2)\frac{k(\widetilde{\phi }_2)-k(\widetilde{\phi }_1)}{2k(\widetilde{\phi}_1)k(\widetilde{\phi}_2)}\textbf{u}_2,\textbf{u}\right)=\textbf{0}.
\end{align*}
Notice that using H\"older's inequality  and \eqref{elliptic}, since $\nu\in W^{1,\infty}(\R)$ and $\textbf{u}_2\in X_{T_0}$, we find
\begin{align*}
&\vert((\nu(\widetilde{\varphi}_1)-\nu(\widetilde{\varphi}_2)) D\textbf{u}_2,D\textbf{u})\vert\leq C\Vert \widetilde{\varphi}\Vert_{L^\infty(\Omega)}\Vert D\textbf{u}_2\Vert\Vert D\textbf{u}\Vert \leq C\Vert \widetilde{\varphi}\Vert_{H^2(\Omega)}\Vert D\textbf{u}_2\Vert\Vert D\textbf{u}\Vert\\&\leq C_\alpha\Vert \widetilde{\phi}\Vert\Vert D\textbf{u}_2\Vert\Vert D\textbf{u}\Vert\leq \frac{\nu_*}{4}\Vert D\textbf{u}\Vert^2+C_\alpha\Vert\widetilde{\phi }\Vert^2\leq \frac{\nu_*}{10}\Vert D\textbf{u}\Vert^2+C_\alpha T_0^3e^{CT_0}\Vert\widetilde{\textbf{u}}\Vert_{C([0,T_0];\textbf{H}_\sigma)}^2,
\end{align*}
where we exploited \eqref{phidiff}. By Young's inequality and \eqref{calpha1}, being $\textbf{u}_2=S(\widetilde{\textbf{u}}_2)\in X_{T_0}$, we deduce
$$
\vert\left((\widetilde{\textbf{v}}\cdot\nabla)\textbf{u}_2,\textbf{u}\right)\vert\leq C\Vert\widetilde{\textbf{v}}\Vert_{[L^\infty(\Omega)]^2}\Vert \textbf{u}_2\Vert\Vert D\textbf{u}\Vert\leq C_\alpha\Vert\widetilde{\textbf{u}}\Vert\Vert D\textbf{u}\Vert\leq \frac{\nu_*}{10}\Vert D\textbf{u}\Vert^2+C_\alpha\Vert\widetilde{\textbf{u}}\Vert^2.
$$
Furthermore, being $\lambda_\alpha$ bounded and using \eqref{calpha1}, \eqref{psitilde}, \eqref{Gronwallpsi2}, we get
\begin{align*}
&\left(\lambda_\alpha(\widetilde{\phi}_1)\nabla\widetilde{\psi}_1^T\nabla\widetilde{\psi},\textbf{v}\right)\leq C\Vert \nabla\widetilde{\psi }_1\Vert\Vert\nabla\widetilde{\psi }\Vert\Vert\nabla\textbf{v}\Vert_{L^\infty(\Omega)}\\&\leq C_\alpha \Vert \nabla\widetilde{\psi }_1\Vert\Vert\nabla\widetilde{\psi }\Vert\Vert\textbf{u}\Vert\leq C_\alpha\int_0^{T_0}\Vert\widetilde{\textbf{u}}(s)\Vert ds\Vert\textbf{u}\Vert\leq C_\alpha T_0 \Vert\widetilde{\textbf{u}}\Vert_{C([0,T_0];\textbf{H}_\sigma)}\Vert\textbf{u}\Vert\\&\leq C_\alpha\Vert\textbf{u}\Vert^2+C_\alpha T_0^2 \Vert\widetilde{\textbf{u}}\Vert_{C([0,T_0];\textbf{H}_\sigma)}^2,
\end{align*}
and, arguing similarly, we obtain
$$
\left(\lambda_\alpha(\widetilde{\phi}_2)\nabla\widetilde{\psi}^T\nabla\widetilde{\psi}_2,\nabla\textbf{v}\right)\leq C_\alpha\Vert\textbf{u}\Vert^2+C_\alpha T_0^2 \Vert\widetilde{\textbf{u}}\Vert_{C([0,T_0];\textbf{H}_\sigma)}^2.
$$
Recalling now that $\lambda_\alpha$ globally Lipschitz, owing to \eqref{calpha1}, \eqref{psitilde} and \eqref{phidiff}, we infer
\begin{align*}
&\left((\lambda_\alpha(\widetilde{\phi}_1)-\lambda_\alpha(\widetilde{\phi}_2))\nabla\widetilde{\psi}_1^T\nabla\widetilde{\psi}_2,\nabla\textbf{v}\right)\leq C\Vert\widetilde{\phi }\Vert\Vert\nabla\widetilde{\psi}_1\Vert_{[L^4(\Omega)]^4}\Vert\nabla\widetilde{\psi}_2\Vert_{[L^4(\Omega)]^4}\Vert\nabla\textbf{v}\Vert_{L^\infty(\Omega)}\\&\leq
C_\alpha\Vert\widetilde{\phi }\Vert\Vert\textbf{u}\Vert\leq C_\alpha\Vert\textbf{u}\Vert^2+C_\alpha\Vert\widetilde{\phi }\Vert^2\leq C_\alpha\Vert\textbf{u}\Vert^2+C_\alpha T_0^3e^{CT_0}\Vert\widetilde{\textbf{u}}\Vert_{C([0,T_0];\textbf{H}_\sigma)}^2.
\end{align*}
On the other hand, by \eqref{calpha1}, \eqref{phi}, \eqref{V} and \eqref{phidiff}, and by Young's inequality, we have
\begin{align*}
&\left(\nabla\widetilde{\phi}_1\otimes\nabla\widetilde{\phi},\nabla\textbf{v}\right)\leq \Vert\nabla\widetilde{\phi}\Vert\Vert\nabla\widetilde{\phi}_1\Vert\Vert\nabla\textbf{v}\Vert_{L^\infty(\Omega)}\\&\leq C_\alpha\Vert\textbf{u}\Vert^2+C_\alpha\left(1+\Vert\Delta\widetilde{\phi}_1\Vert\right)\left(\Vert\widetilde{\phi}\Vert^2+\Vert\widetilde{\phi}\Vert\Vert\Delta\widetilde{\phi }\Vert\right)\\&\leq C_\alpha\Vert\textbf{u}\Vert^2+C_\alpha\left(1+\Vert\Delta\widetilde{\phi}_1\Vert^2\right)\Vert\widetilde{\phi}\Vert^2+C_\alpha\Vert\Delta\widetilde{\phi}\Vert^2\\&\leq C_\alpha\Vert\textbf{u}\Vert^2+C_\alpha\left(1+\Vert\Delta\widetilde{\phi}_1\Vert^2\right) T_0^3e^{CT_0}\Vert\widetilde{\textbf{u}}\Vert_{C([0,T_0];\textbf{H}_\sigma)}^2+C_\alpha\Vert\Delta\widetilde{\phi}\Vert^2
\end{align*}
and, similarly, we find
\begin{align*}
\left(\nabla\widetilde{\phi}\otimes\nabla\widetilde{\phi}_2,\nabla\textbf{v}\right)\leq C_\alpha\Vert\textbf{u}\Vert^2+C_\alpha\left(1+\Vert\Delta\widetilde{\phi}_2\Vert^2\right) T_0^3e^{CT_0}\Vert\widetilde{\textbf{u}}\Vert_{C([0,T_0];\textbf{H}_\sigma)}^2+C_\alpha\Vert\Delta\widetilde{\phi}\Vert^2.
\end{align*}
Using Sobolev embeddings and Poincar\'{e}'s inequality, being $\textbf{u}_i=S(\widetilde{\textbf{u}}_i)\in X_{T_0}$, $i=1,2$, and by \eqref{V} and \eqref{phidiff}, we deduce
\begin{align*}
\left\vert\left(\nu(\widetilde{\phi }_2)\frac{1-\widetilde{\phi}_2}{2k(\widetilde{\phi}_1)}\textbf{u},\textbf{u}\right)\right\vert\leq C\Vert\textbf{u}\Vert^2
\end{align*}
and
\begin{align*}
&\left\vert\left(\nu(\widetilde{\phi}_1)\frac{\left(\widetilde{\phi}_1-\widetilde{\phi}_2\right)\textbf{u}_1}{2k(\widetilde{\phi }_1)},\textbf{u}\right)\right\vert\leq C\Vert\widetilde{\phi }\Vert_{L^4(\Omega)}\Vert\textbf{u}_1\Vert\Vert\textbf{u}\Vert_{[L^4(\Omega)]^2}\\&\leq  C_\alpha\Vert\nabla\widetilde{\phi }\Vert^2+\frac{\nu}{4}\Vert D\textbf{u}\Vert^2\leq
C_\alpha\left(\Vert\Delta\widetilde{\phi}\Vert\Vert\widetilde{\phi}\Vert+\Vert\widetilde{\phi}\Vert^2\right)+\frac{\nu}{4}\Vert D\textbf{u}\Vert^2\\
&\leq C_\alpha\Vert\widetilde{\phi}\Vert^2+C_\alpha\Vert\Delta\widetilde{\phi}\Vert^2+\frac{\nu_*}{10}\Vert D\textbf{u}\Vert^2\\
&\leq C_\alpha T_0^3e^{CT_0}\Vert\widetilde{\textbf{u}}\Vert_{C([0,T_0];\textbf{H}_\sigma)}^2+C\Vert\Delta\widetilde{\phi}\Vert^2+\frac{\nu_*}{10}\Vert D\textbf{u}\Vert^2.
\end{align*}
On account of $\nu\in W^{1,\infty}(\Omega)$, a similar argument gives
\begin{align*}
&\left\vert\left(\left(\nu(\widetilde{\phi }_2)-\nu(\widetilde{\phi }_2)\right)\frac{1-\widetilde{\phi}_2}{2k(\widetilde{\phi}_1)}\textbf{u}_1,\textbf{u}\right)\right\vert\leq C\Vert\widetilde{\phi }\Vert_{L^4(\Omega)}\Vert\textbf{u}_1\Vert\Vert\textbf{u}\Vert_{[L^4(\Omega)]^2}\\&\leq C_\alpha T_0^3e^{CT_0}\Vert\widetilde{\textbf{u}}\Vert_{C([0,T_0];\textbf{H}_\sigma)}^2+C\Vert\Delta\widetilde{\phi}\Vert^2+\frac{\nu_*}{10}\Vert D\textbf{u}\Vert^2.
\end{align*}
In conclusion, arguing as above and using that fact that $k\in W^{1,\infty}(\R)$, we obtain
\begin{align*}
&\left\vert\left(\nu(\widetilde{\phi }_2)(1-\widetilde{\phi}_2)\frac{k(\widetilde{\phi }_2)-k(\widetilde{\phi }_1)}{2k(\widetilde{\phi}_1)k(\widetilde{\phi}_2)}\textbf{u}_2,\textbf{u}\right)\right\vert\leq C\Vert\widetilde{\phi }\Vert_{L^4(\Omega)}\Vert\textbf{u}_2\Vert\Vert\textbf{u}\Vert_{[L^4(\Omega)]^2}\\&\leq C_\alpha T_0^3e^{CT_0}\Vert\widetilde{\textbf{u}}\Vert_{C([0,T_0];\textbf{H}_\sigma)}^2+C\Vert\Delta\widetilde{\phi}\Vert^2+\frac{\nu_*}{10}\Vert D\textbf{u}\Vert^2.
\end{align*}
Therefore, recalling \eqref{nuA}, we have, for any $t\in[0,T_0]$,
$$
\frac{d}{dt}\frac{1}{2}\Vert\textbf{u}\Vert^2+\frac{\nu_*}{2}\Vert D\textbf{u}\Vert^2\leq C_\alpha\Vert\textbf{u}\Vert^2+C_\alpha(T_0)\Vert\widetilde{\textbf{u}}\Vert_{C([0,T_0];\textbf{H}_\sigma)}^2+C_\alpha(T_0)\Vert\Delta\widetilde{\phi}\Vert^2,
$$
which entails, by Gronwall's Lemma (see also \eqref{phidiff}) the continuity of the map $S: X_{T_0}\to X_{T_0}$ with respect to the topology of $C([0,T_0];\textbf{H}_\sigma)$.
Thus Schauder's fixed point theorem yields that  $S$ has a fixed point in $X_{T_0}$, that is, a solution to \eqref{problem}-\eqref{problem_BIC}.

More precisely, we have found a triplet $(\textbf{u}_{\alpha},\phi_{\alpha},\psi_{\alpha})$ which solves \eqref{problem}-\eqref{problem_BIC} and such that, for any $q\in[2,\infty)$,
\begin{align}
\begin{cases}
\textbf{u}_{\alpha}\in L^\infty(0,T_0;\textbf{V}_\sigma)\cap L^2(0,T_0;\textbf{W}_\sigma)\cap H^1(0,T_0;\textbf{H}_\sigma),\\
\Vert\phi_\alpha\Vert_{L^\infty(\Omega\times(0,T))}\leq 1,\\\phi_{\alpha}\in L^\infty(0,T_0;W^{2,q}(\Omega))\cap H^1(0,T_0;V),\\
{\mu}_{\alpha}\in L^\infty(0,T_0;V),\\
F^{\prime}(\phi_\alpha),\ F^{\prime\prime}(\phi_\alpha)\in L^\infty(0,T_0;L^q(\Omega)),\\
{\psi}_{\alpha}\in C([0,T_0];C^2(\overline{\Omega}))\cap C^1([0,T_0];C^1(\overline{\Omega})).
\end{cases}
\label{regularities}
\end{align}
In particular, observe that  $\textbf{u}_{\alpha}\in C([0,T_0];\textbf{V}_\sigma)$ implies
$$
\textbf{v}_{\alpha}\in C([0,T_0];[H^5(\Omega)]^2)\hookrightarrow C([0,T_0];C^3(\overline{\Omega})),
$$
so that (see Lemma \ref{transport})
\begin{align}
\psi_{\alpha}\in  C^1([0,T_0];C^2(\overline{\Omega})).
\label{psir}
\end{align}
\subsection{Higher-order regularity properties of the approximated solution}
\label{high}
Here we obtain some useful regularity properties which will be used next.

Let us define, for a generic function $g$ with values in a (real) Banach space,
$$
\partial_t^hg(t)=\dfrac{g(t+h)-g(t)}{h}, \quad h>0,\, t+h\in(0,T_0).
$$
By the linearity of the operator $(I+\alpha\textbf{A})^{-2}$ and by \eqref{lp}, we infer
$$
\Vert\partial_t^h\textbf{v}_{\alpha}\Vert\leq C\Vert\partial_t^h\textbf{u}_{\alpha}\Vert,
$$
and also (see \eqref{regularities})
$$
\Vert\partial_t^h\textbf{v}_{\alpha}\Vert_{L^2(0,T_0-h,\textbf{H}_\sigma)}\leq C\Vert\partial_t^h\textbf{u}_{\alpha}\Vert_{L^2(0,T_0-h,\textbf{H}_\sigma)}\leq C\Vert\partial_t\textbf{u}_{\alpha}\Vert_{L^2(0,T_0,\textbf{H}_\sigma)}\leq C_\alpha,
$$
uniformly in $h$. Thus $\partial_t\textbf{v}_{\alpha}\in L^2(0,T_0;\textbf{H}_\sigma)$ and $\partial_t^h\textbf{v}_{\alpha}\rightarrow \partial_t\textbf{v}_{\alpha}$ in $L^2(0,T_0; \textbf{H}_\sigma)$. Therefore, passing to the limit up to subsequences, for almost any $t\in (0,T_0)$, we get
\begin{align}
\Vert\partial_t\textbf{v}_{\alpha}\Vert\leq C\Vert\partial_t\textbf{u}_{\alpha}\Vert,
\label{dtu}
\end{align}
independently of $\alpha$.
Then we observe that $\mu_\alpha$ solves, in distributional sense, the following elliptic problem
$$
\begin{cases}
-\Delta\mu_\alpha=-\partial_t\phi_{\alpha}-\textbf{v}_{\alpha}\cdot \nabla\phi_\alpha\quad\text{ in }\Omega\times(0,T_0),\\
\partial_\textbf{n}\mu_{\alpha}=0\quad\text{ on }\partial\Omega\times(0,T_0).
\end{cases}
$$
However, from  \eqref{calpha1} and \eqref{regularities}, we have
\begin{align*}
\Vert\nabla(\textbf{v}_{\alpha}\cdot \nabla\phi_\alpha)\Vert&\leq C\Vert\textbf{v}_{\alpha}\Vert_{[W^{1,\infty}(\Omega)]^2}\Vert\phi_\alpha\Vert_V+C\Vert\textbf{v}_{\alpha}\Vert_{[L^{\infty}(\Omega)]^2}\Vert\phi_\alpha\Vert_{H^2(\Omega)}\\&\leq C_\alpha \Vert\textbf{u}_{\alpha}\Vert+C_\alpha\Vert\textbf{u}_{\alpha}\Vert\Vert\phi_\alpha\Vert_{H^2(\Omega)},
\end{align*}
which means, being $\textbf{u}_{\alpha}\in C([0,T_0];\textbf{H}_\sigma)$ and $\phi_\alpha\in L^2(0,T_0;H^2(\Omega))$, $\textbf{v}_{\alpha}\cdot \nabla\phi_\alpha\in L^2(0,T_0;V)$.
This fact, together with $\partial_t\phi_{\alpha}\in L^2(0,T_0;V)$, entails, by elliptic regularity,
\begin{align}
\mu_{\alpha}\in L^2(0,T_0;H^3(\Omega)).
\label{mu3}
\end{align}
Furthermore, we set
$$
\partial_t^h\mu_{\alpha}(t)=\dfrac{\mu_{\alpha}(t+h)-\mu_{\alpha}(t)}{h}, \quad h>0,\, t+h\in(0,T_0)
$$
and notice that, for any $v\in V$, integrating by parts,
\begin{align*}
&(\partial_t^h\mu_{\alpha}(t),v)=(\nabla\partial_t^h\phi_\alpha(t),\nabla v)+(\partial_t^h\Psi^\prime(\phi_\alpha(t)),v)\\&+\left(\frac{\lambda_\alpha^\prime(\phi_\alpha(t+h))}{2}\partial_t^h\vert\nabla\psi_\alpha(t\vert^2,v\right)
+\left(\frac{1}{2}\vert\nabla\psi_\alpha(t)\vert^2\partial_t^h\lambda^{\prime}_\alpha(\phi_\alpha(t)),v\right).
\end{align*}
Note that, by \eqref{regularities} and being $\Psi^{\prime\prime}$ convex, it holds
\begin{align*}
&(\partial_t^h\Psi^\prime(\phi_\alpha(t)),v)=\left(\int_0^1\Psi^{\prime\prime}(s\phi_\alpha(t+h)+(1-s)\phi_\alpha(t))\partial_t^h\phi_\alpha ds,v\right)\\&\leq C_\alpha(\Vert\Psi^{\prime\prime}(\phi_\alpha(t+h))\Vert_{L^4(\Omega)}+\Vert\Psi^{\prime\prime}(\phi_\alpha(t))\Vert_{L^4(\Omega)})\Vert\partial_t^h\phi_\alpha\Vert\Vert v\Vert_{L^4(\Omega)}\leq C\Vert\partial_t^h\phi_\alpha\Vert\Vert v\Vert_V.
\end{align*}
Moreover, on account of the boundedness of $\lambda_\alpha^\prime, \lambda_\alpha^{\prime\prime}$, we get
$$
\left(\frac{\lambda_\alpha^\prime(\phi_\alpha(t+h))}{2}\partial_t^h\vert\nabla\psi_\alpha(t)\vert^2,v\right)\leq C\Vert \partial_t^h\vert\nabla\psi_\alpha(t)\vert^2\Vert\Vert v\Vert,
$$
and, by H\"older's inequality,
\begin{align*}
&\left(\frac{1}{2}\vert\nabla\psi_\alpha(t)\vert^2\partial_t^h\lambda^{\prime}_\alpha(\phi_\alpha(t)),v\right)\\&=\left(\frac{1}{2}\vert\nabla\psi_\alpha(t)\vert^2\int_0^1\lambda^{\prime\prime}_\alpha(s\phi_\alpha(t+h)+(1-s)\phi_\alpha(t))\partial_t^h\phi_\alpha ds,v\right)\\&\leq C\Vert \nabla\psi_\alpha\Vert^2_{L^4(\Omega)}\Vert\partial_t^h\phi_\alpha \Vert_{L^4(\Omega)}\Vert v\Vert_{L^4(\Omega)}\leq C\Vert\psi_\alpha\Vert_{[W^{1,4}(\Omega)]^2}^2\Vert\partial_t^h\phi_\alpha\Vert_V\Vert v\Vert_V
\end{align*}
Therefore, recalling \eqref{regularities} we deduce that $\Vert\partial_t^h\mu_\alpha\Vert_{L^2(0,T_0-h,V^\prime)}\leq C$ for any $h>0$, which implies
$\partial_t\mu_\alpha\in L^2(0,T_0; V^\prime)$.
Thus, together with \eqref{mu3}, we obtain
$$
\frac{1}{2}\frac{d}{dt}\Vert\nabla\mu\Vert^2=<\partial_t\mu, -\Delta\mu>\in L^1(0,T_0).
$$
Hence (see also \cite{Giorgini} and \cite{HeWu})
\begin{align}
\Vert\nabla\mu_\alpha\Vert^2\in AC([0,T_0]).
\label{muAC}
\end{align}

In order to make the rest of the computations rigorous, we need to find higher-order regularity on $\textbf{u}_\alpha$.
In particular, given $h>0$, we consider the equation for $\textbf{u}_\alpha$ evaluated at $t+h$ and $t$ ($0<t<T_0-h$), take the difference and divide by $h$. We get, for any $\textbf{w}\in \textbf{V}_\sigma$,
\begin{align}
\nonumber&(\partial_t\partial_t^h \textbf{u}_\alpha(t),\textbf{w})+(\nu(\varphi_\alpha(t))D\partial^h_t\textbf{u}_\alpha(t), D\textbf{w}_\alpha)\\&+(\partial_t^h\nu(\varphi_\alpha(t))D\textbf{u}_\alpha(t+h),\textbf{w})\nonumber
+((\partial_t^h\textbf{v}_\alpha(t)\cdot\nabla)\textbf{u}_\alpha(t+h),\textbf{w})\\&+ ((\textbf{v}_\alpha(t)\cdot \nabla)\partial_t^h\textbf{u}_\alpha(t),\textbf{w})-(\partial_t^h\lambda(\phi_\alpha(t))\nabla^T\psi_\alpha(t+h)\nabla\psi_\alpha(t+h),\nabla (I+\alpha\textbf{A})^{-2}\textbf{w})\nonumber\\&-(\lambda(\phi_\alpha(t))\partial_t^h\nabla^T{\psi_\alpha}(t)\nabla\psi_\alpha(t+h),\nabla (I+\alpha\textbf{A})^{-2}\textbf{w})\nonumber\\&-(\lambda(\phi_\alpha(t))\nabla^T\psi_\alpha(t)\partial_t^h\nabla\psi_\alpha(t),\nabla (I+\alpha\textbf{A})^{-2}\textbf{w})\nonumber\\&-(\nabla\partial_t^h\phi(t)\otimes\nabla\phi_\alpha(t+h),\nabla (I+\alpha\textbf{A})^{-2}\textbf{w})-(\nabla\phi_\alpha(t)\otimes\nabla\partial_t^h\phi_\alpha(t),\nabla (I+\alpha\textbf{A})^{-2}\textbf{w})\nonumber\\&+\left(\frac{\partial_t^h\nu(\phi_\alpha(t))}{2k(\phi_\alpha(t+h))}(1-\phi_\alpha(t+h))\textbf{u}_\alpha(t+h),\textbf{w}\right)\nonumber\\&-\left(\frac{\nu(\phi_\alpha(t))\partial_t^hk(\phi_\alpha(t))}{2k(\phi_\alpha(t))k(\phi_\alpha(t+h))}(1-\phi_\alpha(t+h))\textbf{u}_\alpha(t+h),\textbf{w}\right)\nonumber\\&
-\left(\frac{\nu(\phi_\alpha(t))}{2k(\phi_\alpha(t))}\partial_t^h\phi_\alpha(t)\textbf{u}_\alpha(t+h),\textbf{w}\right)+\left(\frac{\nu(\phi_\alpha(t))}{2k(\phi_\alpha(t))}(1-\phi_\alpha(t))\partial_t^h\textbf{u}_\alpha(t),\textbf{w}\right)=0.
\label{partialta}
\end{align}
If we now set $\textbf{w}=\partial_t^h \textbf{u}_\alpha(t)$ we get
\begin{align}
\nonumber&\frac{1}{2}\frac{d}{dt}\Vert\partial_t^h \textbf{u}_\alpha(t)\Vert^2+(\nu(\varphi_\alpha(t))D\partial^h_t\textbf{u}_\alpha(t), D\partial_t^h \textbf{u}_\alpha(t))\\&
=-(\partial_t^h\nu(\varphi_\alpha(t))D\textbf{u}_\alpha(t+h),\partial_t^h \textbf{u}_\alpha(t))\nonumber - ((\partial_t^h\textbf{v}_\alpha(t)\cdot\nabla)\textbf{u}_\alpha(t+h),\partial_t^h \textbf{u}_\alpha(t))\\&+(\partial_t^h\lambda(\phi_\alpha(t))\nabla^T\psi_\alpha(t+h)\nabla\psi_\alpha(t+h),\nabla (I+\alpha\textbf{A})^{-2}\partial_t^h \textbf{u}_\alpha(t))\nonumber\\&+(\lambda(\phi_\alpha(t))\partial_t^h\nabla^T{\psi_\alpha}(t)\nabla\psi_\alpha(t+h),\nabla (I+\alpha\textbf{A})^{-2}\partial_t^h \textbf{u}_\alpha(t))\nonumber\\&+(\lambda(\phi_\alpha(t))\nabla^T\psi_\alpha(t)\partial_t^h\nabla\psi_\alpha(t),\nabla (I+\alpha\textbf{A})^{-2}\partial_t^h \textbf{u}_\alpha(t))\nonumber\\&+(\nabla\partial_t^h\phi(t)\otimes\nabla\phi_\alpha(t+h),\nabla (I+\alpha\textbf{A})^{-2}\partial_t^h \textbf{u}_\alpha(t))\nonumber\\&\nonumber+(\nabla\phi_\alpha(t)\otimes\nabla\partial_t^h\phi_\alpha(t),\nabla (I+\alpha\textbf{A})^{-2}\partial_t^h \textbf{u}_\alpha(t))\nonumber\\&-\left(\frac{\partial_t^h\nu(\phi_\alpha(t))}{2k(\phi_\alpha(t+h))}(1-\phi(t+h))\textbf{u}_\alpha(t+h),\partial_t^h \textbf{u}_\alpha(t)\right)\nonumber\\&+\left(\frac{\nu(\phi_\alpha(t))\partial_t^hk(\phi_\alpha(t))}{2k(\phi_\alpha(t))k(\phi_\alpha(t+h))}(1-\phi_\alpha(t+h))\textbf{u}_\alpha(t+h),\partial_t^h \textbf{u}_\alpha(t)\right)\nonumber\\&
+\left(\frac{\nu(\phi_\alpha(t))}{2k(\phi_\alpha(t))}\partial_t^h\phi_\alpha(t)\textbf{u}_\alpha(t+h),\partial_t^h \textbf{u}_\alpha(t)\right)\nonumber\\&
-\left(\frac{\nu(\phi_\alpha(t))}{2k(\phi_\alpha(t))}(1-\phi_\alpha(t))\partial_t^h\textbf{u}_\alpha(t),\partial_t^h \textbf{u}_\alpha(t)\right)=0,
\label{partialta1}
\end{align}
which we rewrite, recalling \eqref{nuA}, in the following form
\begin{align}
\frac{1}{2}\frac{d}{dt}\Vert\partial_t^h \textbf{u}_\alpha(t)\Vert^2+\nu_*\Vert D\partial_t^h \textbf{u}_\alpha(t)\Vert^2 \leq \sum_{i=1}^{11}\mathcal{K}_i,
\label{vv}
\end{align}
where $\mathcal{K}_i$ is the $i$-th term on the right-hand side of \eqref{partialta1}.
We first observe that, arguing as we did to obtain \eqref{dtu}, we get, by \eqref{calpha1} and \eqref{elliptic},
\begin{align}
\Vert\partial_t^h\textbf{v}_\alpha\Vert_{[H^2(\Omega)]^2}\leq C_\alpha \Vert\partial_t^h\textbf{u}_\alpha\Vert,\quad\Vert\partial_t^h\varphi_\alpha\Vert_{H^2(\Omega)}\leq C_\alpha \Vert\partial_t^h\phi_\alpha\Vert,
\label{dtt}
\end{align}
therefore, being $\nu\in W^{1,\infty}(\R)$ and on account of \eqref{regularities}, we get
\begin{align*}
\vert\mathcal{K}_1\vert&\leq C\Vert\partial_t^h\varphi_\alpha(t)\Vert_{L^\infty(\Omega)}\Vert D\textbf{u}_\alpha(t+h)\Vert\Vert\partial_t^h\textbf{u}_\alpha(t)\Vert\\&\leq C\Vert\partial_t^h\varphi_\alpha(t)\Vert_{H^2(\Omega)}\Vert D\textbf{u}_\alpha(t+h)\Vert\Vert\partial_t^h\textbf{u}_\alpha\Vert\leq C_\alpha\Vert\partial_t^h\nabla\phi_\alpha(t)\Vert\Vert D\textbf{u}_\alpha(t+h)\Vert\Vert\partial_t^h\textbf{u}_\alpha(t)\Vert\\&\leq C_\alpha\Vert\partial_t^h\textbf{u}_\alpha(t)\Vert^2+C_\alpha\Vert\partial_t^h\nabla\phi_\alpha(t)\Vert^2.
\end{align*}
Then, by \eqref{regularities} and \eqref{dtt}, we have
\begin{align*}
\vert\mathcal{K}_2\vert\leq\Vert\partial_t^h\textbf{v}_\alpha(t)\Vert_{[L^\infty(\Omega)]^2}\Vert\nabla\textbf{u}(t+h)\Vert\Vert\partial_t^h\textbf{u}(t)\Vert\leq C_\alpha\Vert\partial_t^h\textbf{u}(t)\Vert^2.
\end{align*}
By \eqref{calpha1}, being $\lambda_\alpha\in W^{1,\infty}(\R)$, we deduce
\begin{align*}
\vert\mathcal{K}_3\vert&\leq C\Vert\partial_t^h\phi_\alpha(t)\Vert\Vert\nabla\psi_\alpha(t+h)\Vert_{[L^4(\Omega)]^4}^2\Vert\nabla(\textbf{I}+\alpha\textbf{A})^{-2}\partial_t^h\textbf{u}(t)\Vert_{[L^\infty(\Omega)]^4}\\&\leq C_\alpha\Vert\partial_t^h\nabla\phi_\alpha(t)\Vert\Vert\partial_t^h\textbf{u}(t)\Vert\leq C_\alpha\Vert\partial_t^h\textbf{u}(t)\Vert^2+C_\alpha\Vert\partial_t^h\nabla\phi_\alpha(t)\Vert^2.
\end{align*}
Similarly, thanks to \eqref{regularities}, we infer
\begin{align*}
\vert\mathcal{K}_4\vert+\vert\mathcal{K}_5\vert&\leq C\Vert\partial_t^h\nabla\psi_\alpha(t)\Vert(\Vert\nabla\psi_\alpha(t)\Vert+\Vert\nabla\psi_\alpha(t+h)\Vert)\Vert\nabla(\textbf{I}+\alpha\textbf{A})^{-2}\partial_t^h\textbf{u}(t)\Vert_{[L^\infty(\Omega)]^4}\\&\leq C_\alpha\Vert\partial_t^h\textbf{u}(t)\Vert^2+C_\alpha\Vert\partial_t^h\nabla\psi_\alpha(t)\Vert^2,
\end{align*}
and
\begin{align*}
\vert\mathcal{K}_6\vert+\vert\mathcal{K}_7\vert&\leq C\Vert\partial_t^h\nabla\phi_\alpha(t)\Vert(\Vert\nabla\phi_\alpha(t)\Vert+\Vert\nabla\phi_\alpha(t+h)\Vert)\Vert\nabla(\textbf{I}+\alpha\textbf{A})^{-2}\partial_t^h\textbf{u}(t)\Vert_{[L^\infty(\Omega)]^4}\\&\leq C_\alpha\Vert\partial_t^h\textbf{u}(t)\Vert^2+C_\alpha\Vert\partial_t^h\nabla\phi_\alpha(t)\Vert^2.
\end{align*}
Then, on account of \eqref{regularities}, being $\nu,k\in W^{1,\infty}(\R)$ and $\Vert\phi_\alpha\Vert_{L^\infty(\Omega\times(0,T_0))}\leq 1$, we have
\begin{align*}
\vert\mathcal{K}_8\vert+\vert\mathcal{K}_9\vert+\vert\mathcal{K}_{10}\vert&\leq C\Vert\partial_t^h\phi_\alpha(t)\Vert_{L^4(\Omega)}\Vert\textbf{u}_\alpha(t+h)\Vert\Vert\partial_t^h\textbf{u}(t)\Vert_{[L^4(\Omega)]^2}\\&\leq C\Vert\partial_t^h\nabla\phi_\alpha(t)\Vert\Vert\partial_t^h\textbf{u}(t)\Vert^{1/2}\Vert D\partial_t^h\textbf{u}(t)\Vert^{1/2}\\&\leq \frac{\nu_*}{4}\Vert D\partial_t^h\textbf{u}(t)\Vert^2+C\Vert\partial_t^h\textbf{u}(t)\Vert^2+C\Vert\partial_t^h\nabla\phi_\alpha(t)\Vert^2,\\
\vert\mathcal{K}_{11}\vert&\leq C\Vert\partial_t^h\textbf{u}(t)\Vert^2.
\end{align*}
Take now $\textbf{w}=\textbf{A}^{-1}\partial_t^h\partial_t\textbf{u}(t)\in D(\textbf{A})$ for almost any $t\in(0,T_0)$
in \eqref{partialta}. This gives
\begin{align}
\nonumber&\Vert\partial_t^h\partial_t \textbf{u}_\alpha(t)\Vert^2_{\textbf{V}_\sigma^{-1}}+(\nu(\varphi_\alpha(t))D\partial^h_t\textbf{u}_\alpha(t), D\textbf{A}^{-1}\partial_t^h\partial_t \textbf{u}_\alpha(t))\\&(\partial_t^h\nu(\varphi_\alpha(t))D\textbf{u}_\alpha(t+h),\textbf{A}^{-1}\partial_t^h\partial_t \textbf{u}_\alpha(t))\nonumber+((\partial_t^h\textbf{v}_\alpha(t)\cdot\nabla)\textbf{u}_\alpha(t+h),\textbf{A}^{-1}\partial_t^h\partial_t \textbf{u}_\alpha(t))\\&+((\textbf{v}_\alpha(t)\cdot \nabla)\partial_t^h\textbf{u}_\alpha(t),\textbf{A}^{-1}\partial_t^h\partial_t \textbf{u}_\alpha(t))\nonumber\\&-(\partial_t^h\lambda(\phi_\alpha(t))\nabla^T\psi_\alpha(t+h)\nabla\psi_\alpha(t+h),\nabla (I+\alpha\textbf{A})^{-2}\textbf{A}^{-1}\partial_t^h\partial_t \textbf{u}_\alpha(t))\nonumber\\&-(\lambda(\phi_\alpha(t))\partial_t^h\nabla^T{\psi_\alpha}(t)\nabla\psi_\alpha(t+h),\nabla (I+\alpha\textbf{A})^{-2}\textbf{A}^{-1}\partial_t^h\partial_t \textbf{u}_\alpha(t))\nonumber\\&-(\lambda(\phi_\alpha(t))\nabla^T\psi_\alpha(t)\partial_t^h\nabla\psi_\alpha(t),\nabla (I+\alpha\textbf{A})^{-2}\textbf{A}^{-1}\partial_t^h\partial_t \textbf{u}_\alpha(t))\nonumber\\&-(\nabla\partial_t^h\phi(t)\otimes\nabla\phi_\alpha(t+h),\nabla (I+\alpha\textbf{A})^{-2}\textbf{A}^{-1}\partial_t^h\partial_t \textbf{u}_\alpha(t))\nonumber\\&-(\nabla\phi_\alpha(t)\otimes\nabla\partial_t^h\phi_\alpha(t),\nabla (I+\alpha\textbf{A})^{-2}\textbf{A}^{-1}\partial_t^h\partial_t \textbf{u}_\alpha(t))\nonumber\\&+\left(\frac{\partial_t^h\nu(\phi_\alpha(t))}{2k(\phi_\alpha(t+h))}(1-\phi_\alpha(t+h))\textbf{u}_\alpha(t+h),\textbf{A}^{-1}\partial_t^h\partial_t \textbf{u}_\alpha(t)\right)\nonumber\\&-\left(\frac{\nu(\phi_\alpha(t))\partial_t^hk(\phi_\alpha(t))}{2k(\phi_\alpha(t))k(\phi_\alpha(t+h))}(1-\phi_\alpha(t+h))\textbf{u}_\alpha(t+h),\textbf{A}^{-1}\partial_t^h\partial_t \textbf{u}_\alpha(t)\right)\nonumber\\&
-\left(\frac{\nu(\phi_\alpha(t))}{2k(\phi_\alpha(t))}\partial_t^h\phi_\alpha(t)\textbf{u}_\alpha(t+h),\textbf{A}^{-1}\partial_t^h\partial_t \textbf{u}_\alpha(t)\right)\nonumber\\&+\left(\frac{\nu(\phi_\alpha(t))}{2k(\phi_\alpha(t))}(1-\phi_\alpha(t))\partial_t^h\textbf{u}_\alpha(t),\textbf{A}^{-1}\partial_t^h\partial_t \textbf{u}_\alpha(t)\right)=0.
\nonumber\\&
\label{partialta2}
\end{align}
By standard estimates, recalling \eqref{calpha1}, \eqref{regularities}, \eqref{dtt} and the embedding $H^2(\Omega)\hookrightarrow L^\infty(\Omega)$, we find
\begin{align}
\nonumber&\Vert\partial_t^h\partial_t \textbf{u}_\alpha(t)\Vert^2_{\textbf{V}_\sigma^{-1}}\\&\nonumber\leq C\Vert D\partial^h_t\textbf{u}_\alpha(t)\Vert \Vert D\textbf{A}^{-1}\partial_t^h\partial_t \textbf{u}_\alpha(t)\Vert\\&\nonumber +C\Vert\partial_t^h\varphi_\alpha(t)\Vert_{L^\infty(\Omega)}\Vert D\textbf{u}_\alpha(t+h)\Vert\Vert\textbf{A}^{-1}\partial_t^h\partial_t \textbf{u}_\alpha(t)\Vert\\&\nonumber+\Vert\partial_t^h\textbf{v}_\alpha(t)\Vert_{[L^4(\Omega)]^2}\Vert\textbf{u}_\alpha(t+h)\Vert_{[L^4(\Omega)]^2}\Vert\nabla\textbf{A}^{-1}\partial_t^h\partial_t \textbf{u}_\alpha(t)\Vert\\&+\Vert\textbf{v}_\alpha(t)\Vert_{[L^\infty(\Omega)]^2}\Vert\partial_t^h\textbf{u}_\alpha(t)\Vert\Vert\textbf{A}^{-1}\partial_t^h\partial_t \textbf{u}_\alpha(t)\Vert\nonumber\\&+\Vert\partial_t^h\phi_\alpha(t)\Vert\Vert\nabla\psi_\alpha(t+h)\Vert_{[L^4(\Omega)]^4}^2\Vert\nabla (I+\alpha\textbf{A})^{-2}\textbf{A}^{-1}\partial_t^h\partial_t \textbf{u}_\alpha(t)\Vert_{[L^\infty(\Omega)]^4}\nonumber\\&+C\Vert\partial_t^h\nabla{\psi_\alpha}(t)\Vert\Vert\nabla\psi_\alpha(t+h)\Vert\nabla (I+\alpha\textbf{A})^{-2}\textbf{A}^{-1}\partial_t^h\partial_t \textbf{u}_\alpha(t)\Vert_{[L^\infty(\Omega)]^4}\nonumber\\&+C\Vert\partial_t^h\nabla{\psi_\alpha}(t)\Vert\Vert\nabla\psi_\alpha(t)\Vert\nabla (I+\alpha\textbf{A})^{-2}\textbf{A}^{-1}\partial_t^h\partial_t \textbf{u}_\alpha(t)\Vert_{[L^\infty(\Omega)]^4}\nonumber\\&+\Vert\nabla\partial_t^h\phi(t)\Vert\Vert\nabla\phi_\alpha(t+h)\Vert\Vert\nabla (I+\alpha\textbf{A})^{-2}\textbf{A}^{-1}\partial_t^h\partial_t \textbf{u}_\alpha(t)\Vert_{[L^\infty(\Omega)]^4}\nonumber\\&+\Vert\nabla\partial_t^h\phi(t)\Vert\Vert\nabla\phi_\alpha(t)\Vert\Vert\nabla (I+\alpha\textbf{A})^{-2}\textbf{A}^{-1}\partial_t^h\partial_t \textbf{u}_\alpha(t)\Vert_{[L^\infty(\Omega)]^4}\nonumber\\&+C\Vert\partial_t^h\phi_\alpha(t)\Vert\Vert\textbf{u}_\alpha(t+h)\Vert_{[L^4(\Omega)]^2}\Vert\textbf{A}^{-1}\partial_t^h\partial_t\textbf{u}\Vert_{[L^4(\Omega)]^2}+C\Vert \partial_t^h\textbf{u}_\alpha(t)\Vert \Vert\textbf{A}^{-1}\partial_t^h\partial_t \textbf{u}_\alpha(t)\Vert\nonumber\\&\leq \frac{1}{2}\Vert \partial_t^h\partial_t \textbf{u}_\alpha(t)\Vert_{\textbf{V}_\sigma^{-1}}^2 + C_0\Vert D\partial^h_t\textbf{u}_\alpha(t)\Vert^2+C_\alpha \left(\Vert\partial_t^h\textbf{u}_\alpha(t)\Vert^2+ \Vert\partial_t^h\nabla\phi_\alpha(t)\Vert^2\right).
\label{partialta3}\end{align}
Collecting the above estimates, multiply \eqref{partialta3} by $\omega=\frac{\nu_*}{4 C_0}$, and add them together with \eqref{partialta1}, we obtain in the end, recalling \eqref{nuA},
\begin{align*}
&\frac{1}{2}\frac{d}{dt}\Vert\partial_t^h \textbf{u}_\alpha(t)\Vert^2+\frac{\nu_*}{2}\Vert D\partial^h_t\textbf{u}_\alpha(t)\Vert^2+\frac{\omega}{2}\Vert\partial_t^h\partial_t \textbf{u}_\alpha(t)\Vert^2_{\textbf{V}_\sigma^{-1}}\\&\leq C_\alpha \left(\Vert\partial_t^h\textbf{u}_\alpha(t)\Vert^2+ \Vert\partial_t^h\nabla\phi_\alpha(t)\Vert^2+\Vert\partial_t^h\nabla\psi_\alpha(t)\Vert^2\right).
\end{align*}
We then integrate in time over $(0,t)$ to get
\begin{align}
&\nonumber\frac{1}{2}\Vert\partial_t^h \textbf{u}_\alpha(t)\Vert^2+\int_0^t\left(\frac{\nu_*}{2}\Vert D\partial^h_t\textbf{u}_\alpha(s)\Vert^2+\frac{\omega}{2}\Vert\partial_t^h\partial_t \textbf{u}_\alpha(s)\Vert^2_{\textbf{V}_\sigma^{-1}}\right)ds\\&\leq \frac{1}{2}\Vert\partial_t^h \textbf{u}_\alpha(0)\Vert^2+C_\alpha \int_0^t\left(\Vert\partial_t^h\textbf{u}_\alpha(s)\Vert^2+ \Vert\partial_t^h\nabla\phi_\alpha(s)\Vert^2+\Vert\partial_t^h\nabla\psi_\alpha(s)\Vert^2\right)ds.\label{gronwal}
\end{align}
The second term on the right-hand side is bounded uniformly in $h$ since, by \eqref{regularities} and \eqref{psir}, we have $\partial_t\textbf{u}_\alpha\in L^2(0,T;\textbf{H}_\sigma)$, $\partial_t\phi_\alpha\in L^2(0,T;V)$ and $\partial_t\nabla\psi_\alpha\in C([0,T_0];C^1(\overline{\Omega}))$. Therefore, we only need to find an $h$-uniform estimate for the initial datum. To this purpose, observe that $\textbf{u}_\alpha(t)-\textbf{u}_{0,\alpha}$ solves
\begin{align*}
&\frac{1}{2}\frac{d}{dt}\Vert\textbf{u}_\alpha-\textbf{u}_{0,\alpha}\Vert^2+((\textbf{v}_\alpha\cdot\nabla)\textbf{u}_{0,\alpha},\textbf{u}_\alpha-\textbf{u}_{0,\alpha})+(\nu(\varphi_\alpha) D(\textbf{u}_\alpha-\textbf{u}_{0,\alpha}),D(\textbf{u}_\alpha-\textbf{u}_{0,\alpha}))\\&-(\nu^\prime(\varphi_\alpha) D\textbf{u}_{0,\alpha}\nabla\varphi_\alpha,\textbf{u}_\alpha-\textbf{u}_{0,\alpha})-\frac{1}{2}(\nu(\varphi_\alpha) \Delta\textbf{u}_{0,\alpha},\textbf{u}_\alpha-\textbf{u}_{0,\alpha})\\&+(\text{div}(\lambda_\alpha(\phi)\nabla\psi^T\nabla\psi), (I+\alpha \textbf{A})^{-2}(\textbf{u}_\alpha-\textbf{u}_{0,\alpha}))\\&+(\text{div}(\nabla\phi\otimes\nabla\phi), (I+\alpha \textbf{A})^{-2}(\textbf{u}_\alpha-\textbf{u}_{0,\alpha}))+\left(\nu(\phi)\frac{(1-\phi)\textbf{u}_\alpha}{2k(\phi)}, \textbf{u}_\alpha-\textbf{u}_{0,\alpha}\right)=0,
\end{align*}
which implies (see \eqref{nuA}, \eqref{calpha1}, \eqref{elliptic} and \eqref{regularities})
\begin{align*}
&\frac{1}{2}\frac{d}{dt}\Vert\textbf{u}_\alpha-\textbf{u}_{0,\alpha}\Vert^2\\
&\leq \Vert\textbf{v}_\alpha\Vert_{[L^\infty(\Omega)]^2}\Vert \nabla\textbf{u}_{0,\alpha}\Vert\Vert \textbf{u}_\alpha-\textbf{u}_{0,\alpha}\Vert +C\Vert D\textbf{u}_{0,\alpha}\Vert\Vert\nabla\varphi_\alpha\Vert_{[L^\infty(\Omega)]^2}\Vert\textbf{u}_\alpha-\textbf{u}_{0,\alpha}\Vert\\&+C\Vert\Delta\textbf{u}_{0,\alpha}\Vert\Vert\textbf{u}_\alpha-\textbf{u}_{0,\alpha}\Vert
+C\Vert\nabla\psi_\alpha\Vert_{[L^4(\Omega)]^4}^2\Vert\nabla(I+\alpha \textbf{A})^{-2}(\textbf{u}_\alpha-\textbf{u}_{0,\alpha})\Vert_{[L^\infty(\Omega)]^4}\\&+C\Vert\nabla\phi_\alpha\Vert_{[L^4(\Omega)]^2}^2\Vert\nabla(I+\alpha \textbf{A})^{-2}(\textbf{u}_\alpha-\textbf{u}_{0,\alpha})\Vert_{[L^\infty(\Omega)]^4}+C\Vert\textbf{u}_\alpha\Vert\Vert \textbf{u}_\alpha-\textbf{u}_{0,\alpha}\Vert\\
&\leq C_\alpha\Vert \textbf{u}_\alpha-\textbf{u}_{0,\alpha}\Vert,
\end{align*}
where we exploited the fact that $\textbf{u}_{0,\alpha}\in \textbf{W}_\sigma$. Therefore, by a well-known version of Gronwall's Lemma (see \cite[Lemma A.5]{Brezis}), we obtain
$$
\Vert\textbf{u}_\alpha-\textbf{u}_{0,\alpha}\Vert\leq C_\alpha t,\qquad \forall t\in[0,T_0],
$$
and choosing $t=h$, we find that,
$$
\Vert\partial_t^h\textbf{u}_\alpha(0)\Vert\leq C_\alpha,
$$
uniformly with respect to $h$.
Then, from \eqref{gronwal}, passing to the limit as $h\to0$, we deduce
\begin{align}
\Vert\partial_t\textbf{u}_\alpha\Vert_{L^\infty(0,T_0;\textbf{H}_\sigma)}+\Vert\partial_t\textbf{u}_\alpha\Vert_{L^2(0,T_0;\textbf{V}_\sigma)}+\Vert\partial_{tt}\textbf{u}_\alpha\Vert_{L^2(0,T_0;\textbf{V}_\sigma^{-1})}\leq C_\alpha.
\label{newreg}
\end{align}
Therefore, we have $\partial_t\textbf{u}_\alpha\in C([0,T_0];\textbf{H}_\sigma)$ so that
\begin{align}
\textbf{u}_\alpha\in C^1([0,T_0];\textbf{H}_\sigma).
\label{C1}
\end{align}
The regularities obtained so far are enough to make the next steps rigorous and obtain higher-order estimates independent of the regularizing parameter $\alpha$.

\subsection{Estimates independent of $\alpha$}

In this subsection we obtain some estimates on the approximating solution which are independent of $\alpha$.
We will omit the subscripts, in order to simplify the notation. Also, $C$ or $C_i$ will stand for a generic positive constant independent of $t,\alpha$, which may vary from line to line.
\label{ests}
\subsubsection{Energy estimate}
We look first for an energy identity. Thanks to the regularity of the solution, we can multiply \eqref{problem}$_1$ by $\textbf{u}$, \eqref{problem}$_4$ by $\mu$ and in the end we take the gradient of \eqref{problem}$_2$ and multiply (double dot product) by $\lambda_\alpha(\phi)\nabla\psi$. Integrating the resulting identities over $\Omega$, we obtain, being $(I+\alpha\textbf{A})^{-2}\textbf{P}$ self-adjoint and recalling that $\textbf{v}=(I+\alpha\textbf{A})^{-2}\textbf{u}$ and $\varphi=(I+\alpha{A}_1)^{-1}{\phi}$,
\begin{align*}
\begin{cases}
\frac{d}{dt}\frac{1}{2}\Vert\textbf{u}\Vert^2+\int_\Omega\nu(\varphi)\vert D\textbf{u}\vert^2dx-\left(\mu\nabla\phi,\textbf{v}\right)\\+\left(\nabla\psi^T\text{div}(\lambda_\alpha(\phi)\nabla\psi),\textbf{v}\right)+\frac{\nu(\phi)}{2k(\phi)}\left(1-\phi,\vert\textbf{u}\vert^2\right)=0,\\
\frac{d}{dt}\left(\frac{1}{2}\Vert\nabla\phi\Vert^2+\int_\Omega\Psi(\phi)dx\right)+\frac{1}{2}\left(\lambda^\prime_\alpha(\phi)\partial_t\phi,\vert\nabla\psi\vert^2\right)+\Vert\nabla\mu\Vert^2+(\textbf{v}\cdot\nabla\phi,\mu)=0,\\
\left(\partial_t\nabla\psi,\lambda_\alpha(\phi)\nabla\psi\right)-\left(\text{div}(\lambda_\alpha(\phi)\nabla\psi),\textbf{v}\cdot\nabla\psi\right)=0,
\end{cases}
\end{align*}
where we have exploited the following equivalence and the fact that $\text{div}\,\textbf{v}=0$ in $\Omega\times(0,T_0)$:
\begin{align*}
&\text{div}(\lambda_\alpha(\phi)\nabla\psi^T\nabla\psi)+\text{div}(\nabla\phi\otimes\nabla\phi)\\&=\nabla\psi^T\text{div}(\lambda_\alpha(\phi)\nabla\psi)+\frac{\lambda_\alpha(\phi)}{2}\nabla\vert\nabla\psi\vert^2-\mu\nabla\phi+\nabla\left(-\Psi(\phi)+\frac{1}{2}\vert\nabla\phi\vert^2\right)+\frac{\lambda_\alpha^\prime(\phi)}{2}\vert\nabla\psi\vert^2\nabla\phi\\&=-\mu\nabla\phi+\nabla\psi^T\text{div}(\lambda_\alpha(\phi)\nabla\psi)+\nabla\left(-\Psi(\phi)+\frac{1}{2}\vert\nabla\phi\vert^2+\frac{\lambda_\alpha(\phi)}{2}\vert\nabla\psi\vert^2\right).
\end{align*}
 Moreover, in the last estimate we used $\textbf{v}=0$ on $\partial\Omega\times(0,T)$ in the integration by parts.
Therefore, recalling that
$$
\left(\text{div}(\lambda_\alpha(\phi)\nabla\psi),\textbf{v}\cdot\nabla\psi\right)=\left(\nabla\psi^T\text{div}(\lambda_\alpha(\phi)\nabla\psi),\textbf{v}\right),
$$
we can add the three identities together. This gives
\begin{align}
&\nonumber	\frac{d}{dt}\left(\frac{1}{2}\Vert\textbf{u}\Vert^2+\int_\Omega \Psi(\phi)dx+\frac{1}{2}\Vert\nabla\phi\Vert^2+\frac{1}{2}\int_\Omega\lambda_\alpha(\phi)\vert\nabla\psi\vert^2dx\right)\\&+\Vert\nabla\mu\Vert^2+\int_\Omega\nu(\varphi)\vert D\textbf{u}\vert^2dx+\frac{\nu(\phi)}{2k(\phi)}\left(1-\phi,\vert\textbf{u}\vert^2\right)=0.
\end{align}
Then, recalling \eqref{lambdamin}, Gronwall's Lemma entails
\begin{equation}
0\leq \frac{1}{2}\Vert\textbf{u}\Vert^2+\int_\Omega \Psi(\phi)dx+\hat{C}+\frac{1}{2}\Vert\nabla\phi\Vert^2+\frac{1}{2}\int_\Omega\tilde{\lambda}_*\vert\nabla\psi\vert^2dx\leq C, \quad \forall\,t\in(0,T_0).
\label{energy3}
\end{equation}
where $C$ depends on the initial data energy, but neither on $\alpha$ nor on $t$. Note that the above estimate holds because (see \eqref{regularities})
$$
\frac{\nu(\phi)}{2k(\phi)}\left(1-\phi,\vert\textbf{u}\vert^2\right)\geq 0.
$$
Moreover, we also obtain
\begin{align}
\int_0^{T_0 }\Vert\nabla\mu\Vert^2ds+\int_0^{T_0} \nu_*\Vert D\textbf{u}\Vert^2ds \leq C(T_0).
\label{L2}
\end{align}

\subsubsection{Higher-order estimates}
In order to pass to the limit with respect to $\alpha$, the energy estimate does not suffice, but we need some higher-order bounds.

\textbf{First estimate.} Let us multiply \eqref{problem}$_1$ by $\partial_t\textbf{u}$ and integrate over $\Omega$. Note that this can be done, since $\textbf{u}\in H^1(0,T_0;\textbf{H}_\sigma)\cap L^2(0,T_0;\textbf{W}_\sigma)\hookrightarrow AC([0,T_0];\textbf{V}_\sigma)$. We obtain

\begin{align}
\nonumber&\frac{1}{2}\frac{d}{dt}\int_\Omega\nu(\varphi)\vert D\textbf{u}\vert^2dt+\underbrace{((\textbf{v}\cdot\nabla)\textbf{u},\partial_t\textbf{u})}_{I_1}+\Vert\partial_t\textbf{u}\Vert^2+\underbrace{(\text{div}(\lambda_\alpha(\phi)\nabla\psi^T\nabla\psi),(I+\alpha\textbf{A})^{-2}\partial_t\textbf{u})}_{I_2}\\&+\underbrace{\left(\frac{\nu(\phi)}{2k(\phi)}(1-\phi)\textbf{u},\partial_t\textbf{u}\right)}_{I_3}+\underbrace{(\text{div}(\nabla\phi\otimes\nabla\phi),(I+\alpha\textbf{A})^{-2}\partial_t\textbf{u})}_{I_4}=\underbrace{(\nu^\prime(\phi)\partial_t\phi,\vert D\textbf{u}\vert^2)}_{I_5}.
\label{est1}
\end{align}
Considering the corresponding Stokes problem, we get an $H^2$-estimate for the velocity by Theorem \ref{stokes}. In particular, setting $p=s=2$ ad $r=\infty$, we have
\begin{align}
&\nonumber\Vert\textbf{u}\Vert_{[H^2(\Omega)]^2}\leq \Vert\partial_t\textbf{u}\Vert+\Vert(\textbf{v}\cdot \nabla)\textbf{u}\Vert+\Vert(I+\alpha\textbf{A})^{-2}\textbf{P}\text{div}(\nabla\phi\otimes\nabla\phi)\Vert\\&+\Vert(I+\alpha\textbf{A})^{-2}\textbf{P}\text{div}(\lambda_\alpha(\phi)\nabla\psi^T\nabla\psi)\Vert+\left\Vert\frac{\nu(\phi)}{2k(\phi)}(1-\phi)\textbf{u}\right\Vert+C\Vert\nabla\varphi\Vert_{[L^\infty(\Omega)]^2}\Vert D\textbf{u}\Vert.
\label{H^2}
\end{align}
Observe now that, for $i=1,2$,
\begin{equation}
\text{div}(\lambda_\alpha(\phi)\nabla\psi^T\nabla\psi)_i=\lambda_\alpha^\prime\left(\phi)(\nabla\psi^T\nabla\psi\nabla\phi\right)_i
+\lambda_\alpha(\phi)\left(\sum_{k=1}^2\Delta\psi_k(\nabla\psi)_{ki}+\frac{1}{2}(\nabla\vert\nabla\psi\vert^2)_i\right).
\label{psi3}
\end{equation}
Therefore, being $\vert\lambda^\prime_\alpha(\phi)\vert\leq C$ almost everywhere in $\Omega\times(0,T_0)$ (see \eqref{lambdaest}), independently of $\alpha$, using \eqref{lp}, \eqref{energy3} and H\"older and Gagliardo-Nirenberg's inequalities, we deduce
\begin{align}
&\Vert(I+\alpha\textbf{A})^{-2}\textbf{P}\text{div}(\lambda_\alpha(\phi)\nabla\psi^T\nabla\psi)\Vert\leq C\Vert\text{div}(\lambda_\alpha(\phi)\nabla\psi^T\nabla\psi)\Vert\nonumber\\&\leq C\Vert\nabla\psi\Vert^2_{[L^4(\Omega)]^4}\Vert\nabla\phi\Vert+C\Vert\nabla\psi\Vert_{[L^\infty(\Omega)]^4}\Vert\psi\Vert_{[H^2(\Omega)]^2}\nonumber\\&\leq C\Vert\psi\Vert_{[H^2(\Omega)]^2}(1+\Vert\nabla\psi\Vert_{[L^\infty(\Omega)]^4}).
\label{divpsi}
\end{align}
Moreover, by Agmon's inequality, \eqref{lp}, \eqref{w2pv} and \eqref{energy3}, we get
$$
\Vert(\textbf{v}\cdot \nabla)\textbf{u}\Vert\leq C\Vert\textbf{u}\Vert_{[H^2(\Omega)]^2}^{1/2}\Vert D\textbf{u}\Vert\leq C\Vert D\textbf{u}\Vert^2+\frac{1}{2}\Vert\textbf{u}\Vert_{[H^2(\Omega)]^2}.
$$
Then, for $p\in(2,\infty)$ and $1/q+1/p=1/2$, \eqref{energy3}, the embedding $W^{2,p}(\Omega)\hookrightarrow H^2(\Omega)$ and \eqref{lp}, we have
\begin{align}
\nonumber
&\Vert(I+\alpha\textbf{A})^{-2}\textbf{P}\text{div}(\nabla\phi\otimes\nabla\phi)\Vert\leq C \Vert\text{div}(\nabla\phi\otimes\nabla\phi)\Vert\leq C\Vert\nabla\phi\Vert_{L^{q}(\Omega)}\Vert D^2\phi\Vert_{L^{p}(\Omega)}\\&\leq C\Vert\nabla\phi\Vert^{2/q}\Vert\phi\Vert^{1-2/q}_{H^2(\Omega)}\Vert\phi\Vert_{W^{2,p}(\Omega)}\leq C\Vert\phi\Vert_{W^{2,p}(\Omega)}^{2-2/q}=C\Vert\phi\Vert_{W^{2,p}(\Omega)}^{1+2/p}+\Vert\phi\Vert_{W^{2,p}(\Omega)}^2.
\label{phi5}
\end{align}
On the other hand, using the embedding $W^{2,p}(\Omega)\hookrightarrow W^{1,\infty}(\Omega)$ (being $p>2$) and recalling \eqref{lpphi2}, we deduce
$$
C\Vert\nabla\varphi\Vert_{[L^\infty(\Omega)]^2}\Vert D\textbf{u}\Vert\leq C\Vert\varphi\Vert_{W^{2,p}(\Omega)}\Vert D\textbf{u}\Vert
\leq C\left(\Vert\phi\Vert^2_{W^{2,p}(\Omega)}+\Vert D\textbf{u} \Vert^2\right).
$$
Finally, owing to \eqref{regularities} and \eqref{energy3}, we have
$$
\left\Vert\frac{\nu(\phi)}{2k(\phi)}(1-\phi)\textbf{u}\right\Vert\leq C.
$$
Summing up, rearranging the terms, we get
\begin{align}
&\Vert\textbf{u}\Vert_{[H^2(\Omega)]^2} \nonumber\\
&\leq C\left( \Vert\partial_t\textbf{u}\Vert+\Vert\psi\Vert_{[H^2(\Omega)]^2}(1+\Vert\nabla\psi\Vert_{[L^\infty(\Omega)]^4})+\Vert D\textbf{u}\Vert^2+\Vert\phi\Vert_{W^{2,p}(\Omega)}^{1+2/p}+\Vert\phi\Vert_{W^{2,p}(\Omega)}^2\right). \label{H23}
\end{align}
Consider now \eqref{est1}. Recalling again Agmon's inequality, \eqref{lp}, \eqref{w2pv} and \eqref{energy3}, we have
\begin{align*}
&\vert I_1\vert\leq \Vert\textbf{v}\Vert_{[L^\infty(\Omega)]^2}\Vert\nabla\textbf{u}\Vert\Vert\partial_t\textbf{u}\Vert\leq C\Vert\textbf{v}\Vert^{1/2}\Vert\textbf{v}\Vert^{1/2}_{[H^2(\Omega)]^2}\Vert\nabla\textbf{u}\Vert\Vert\partial_t\textbf{u}\Vert\\&\leq C\Vert\textbf{u}\Vert^{1/2}\Vert\textbf{u}\Vert^{1/2}_{[H^2(\Omega)]^2}\Vert\nabla\textbf{u}\Vert\Vert\partial_t\textbf{u}\Vert \leq \delta_1\Vert\textbf{u}\Vert_{[H^2(\Omega)]^2}^2+C(\delta,\delta_1)\Vert D\textbf{u}\Vert^4+\delta \Vert\partial_t\textbf{u}\Vert^2,
\end{align*}
for some $\delta,\delta_1>0$ suitably chosen.
On account of $\vert\lambda_\alpha(\phi)\vert\leq C$ almost everywhere in $\Omega\times(0,T_0)$ (see \eqref{lambdaest}), and by \eqref{lp}, \eqref{energy3} and \eqref{divpsi},
we obtain
\begin{align*}
\vert I_2\vert&\leq \Vert(I+\alpha\textbf{A})^{-2}\textbf{P}\text{div}(\lambda(\phi)\nabla\psi^T\nabla\psi)\Vert\Vert\partial_t\textbf{u}\Vert\leq C\Vert\psi\Vert_{[H^2(\Omega)]^2}(1+\Vert\nabla\psi\Vert_{[L^\infty(\Omega)]^4})\Vert\partial_t\textbf{u}\Vert\\&\leq \delta \Vert\partial_t\textbf{u}\Vert^2+C(\delta)\Vert\psi\Vert_{[H^2(\Omega)]^2}^2(1+\Vert\nabla\psi\Vert_{[L^\infty(\Omega)]^4}^2).
\end{align*}
For $I_3$, being $\vert1-\phi\vert\leq 1$ almost everywhere on $\Omega\times(0,T_0)$ and using \eqref{energy3}, we get
$$\vert I_3\vert\leq C\Vert\textbf{u}\Vert\Vert\partial_t\textbf{u}\Vert\leq C+\delta\Vert\partial_t\textbf{u}\Vert^2.$$
On the other hand, by \eqref{lp}, \eqref{energy3} and \eqref{phi5}, we deduce
$$
\vert I_4\vert\leq \Vert(I+\alpha\textbf{A})^{-2}\textbf{P}\text{div}(\nabla\phi\otimes\nabla\phi)\Vert\Vert\partial_t\textbf{u}\Vert\leq C\Vert\phi\Vert_{W^{2,p}(\Omega)}^{1+2/p}\Vert\partial_t\textbf{u}\Vert\leq  C\Vert\phi\Vert_{W^{2,p}(\Omega)}^{2+4/p}+\delta\Vert\partial_t\textbf{u}\Vert^2.
$$
Concerning $I_5$, Sobolev embeddings yield
\begin{align*}
\vert I_5\vert&\leq C\Vert\partial_t\phi\Vert_{L^4(\Omega)}\Vert D\textbf{u}\Vert_{[L^4(\Omega)]^4}\Vert D\textbf{u}\Vert\leq C\Vert\partial_t\nabla\phi\Vert\Vert\textbf{u}\Vert^{1/2}_{[H^2(\Omega)]^2}\Vert D\textbf{u}\Vert^{3/2}\\&\leq \delta_2\Vert\partial_t\nabla\phi\Vert^2+\delta_1\Vert\textbf{u}\Vert_{[H^2(\Omega)]^2}^2+C(\delta_1,\delta_2)\Vert D\textbf{u}\Vert^6,
\end{align*}
for $\delta_1,\delta_2>0$ sufficiently small.

Taking the above estimates into account, combining \eqref{est1} with \eqref{H23} multiplied by $\omega>0$, we find
\begin{align*}
&\frac{1}{2}\frac{d}{dt}\int_\Omega \nu(\varphi)\vert D\textbf{u}\vert^2dx+\Vert\partial_t\textbf{u}\Vert^2+\omega\Vert\textbf{u}\Vert_{[H^2(\Omega)]^2}^2\\&\leq C(\omega,\delta,\delta_1)\left(1+\Vert D\textbf{u}\Vert^4+\Vert\psi\Vert_{[H^2(\Omega)]^2}^2(1+\Vert\nabla\psi\Vert_{[L^\infty(\Omega)]^4}^2)+\Vert\phi\Vert_{W^{2,p}(\Omega)}^{2+4/p}+\Vert\phi\Vert_{W^{2,p}(\Omega)}^2+\Vert D\textbf{u}\Vert^6\right)\\&+(\omega C+4\delta)\Vert\partial_t\textbf{u}\Vert^2+\delta_1\Vert\textbf{u}\Vert^2_{[H^2(\Omega)]^2}+\delta_2\Vert\partial_t\nabla\phi\Vert^2.
\end{align*}
Therefore, choosing $\omega=\frac{1}{4C}$, $\delta=\frac{1}{16},\ \delta_1=\frac{1}{8C}$, $\delta_2=\frac{1}{8}$, we immediately obtain
 \begin{align}
 \nonumber&\frac{1}{2}\frac{d}{dt}\int_\Omega \nu(\varphi)\vert D\textbf{u}\vert^2dx+\frac{1}{2}\Vert\partial_t\textbf{u}\Vert^2+\frac{1}{8C}\Vert\textbf{u}\Vert_{[H^2(\Omega)]^2}^2\\&\leq C\left(1+\Vert D\textbf{u}\Vert^4+\Vert\psi\Vert_{[H^2(\Omega)]^2}^2(1+\Vert\nabla\psi\Vert_{[L^\infty(\Omega)]^4}^2)+\Vert\phi\Vert_{W^{2,p}(\Omega)}^{2+4/p}+\Vert\phi\Vert_{W^{2,p}(\Omega)}^2+\Vert D\textbf{u}\Vert^6\right)\nonumber\\&+\frac{1}{8}\Vert\partial_t\nabla\phi\Vert^2.
 \label{estvelox}
 \end{align}

\textbf{Second estimate.} We need to estimate $\Vert\psi\Vert_{[W^{2,p}(\Omega)]^2}$. In this case, we can exploit Lemma \ref{transport} and, in particular, \eqref{psiw2p}, for $p\in(2,\infty)$, to get
(see also \eqref{w2pv})
\begin{equation}
\frac{1}{2}\frac{d}{dt}\Vert \psi\Vert_{[W^{2,p}(\Omega)]^2}^2\leq C\Vert\textbf{v}\Vert_{[W^{2,p}(\Omega)]^2}\Vert\psi\Vert_{[W^{2,p}(\Omega)]^2}^2\leq C\Vert\textbf{u}\Vert_{[W^{2,p}(\Omega)]^2}\Vert\psi\Vert_{[W^{2,p}(\Omega)]^2}^2.
\label{psiw2}
\end{equation}

\textbf{Third estimate.} The next inequality we need is the one coming from the Cahn-Hilliard type equation. Firstly, we need some preliminary observations. From \eqref{energy3} and the conservation of mass, we immediately deduce that
\begin{equation}
\Vert\phi\Vert_{V}\leq C, \quad \text{ a.e. in } (0,T_0).
\label{phi323}
\end{equation}
Moreover, we observe that, from \eqref{problem}$_4$, thanks to \eqref{h1v}, \eqref{energy3} and $\Vert\phi\Vert_{L^\infty(\Omega\times(0,T_0))}\leq 1$, we have
\begin{equation}
\Vert \partial_t\phi \Vert_{V^\prime}\leq C\left(\Vert\nabla\mu\Vert+\Vert \textbf{u}\Vert\Vert\phi\Vert_{L^\infty(\Omega)}\right)\leq C\left(\Vert\nabla\mu\Vert+1\right).
\label{phistar223}
\end{equation}
In conclusion, by well-established results (see, e.g., \cite{Giorginitemam}), \eqref{energy3}, Sobolev embeddings and Gagliardo-Nirenberg's inequalities, for any $p\in(2,\infty)$, we find
\begin{align}
&\nonumber\Vert\phi\Vert_{W^{2,p}(\Omega)}+\Vert F^\prime(\phi)\Vert_{L^p(\Omega)}\\&\leq C\left(\Vert\mu\Vert_{L^p(\Omega)}+\Vert\phi\Vert_{L^p(\Omega)}+\Vert\nabla\psi\Vert_{[L^{2p}(\Omega)]^4}^2\right)\leq C\left(1+\Vert\mu\Vert_{V}+\Vert\psi\Vert^{2-2/p}_{[W^{2,p}(\Omega)]^2}\right).
\label{pot}
\end{align}
Observe now that
\begin{equation*}
(\mu,\phi-\overline{\phi})=\Vert   \nabla\phi\Vert   ^2+(F'(\phi),\phi-\overline{\phi})-\alpha_0(\phi,\phi-\overline{\phi})+\left(\lambda^\prime_\alpha(\phi)\frac{\vert\nabla\psi\vert^2}{2},\phi-\overline{\phi}\right).
\end{equation*}
Using $(\overline{\mu},\phi-\overline{\phi})=0$ and applying standard inequalities together with Gagliardo-Nirenberg's inequality, from \eqref{lambdaest} we infer
\begin{align*}
(F'(\phi),\phi-\overline{\phi})=&(\mu-\overline{\mu},\phi-\overline{\phi})-\Vert   \nabla\phi\Vert   ^2+\alpha_0(\phi,\phi-\overline{\phi})-\left(\lambda^\prime_\alpha(\phi)\frac{\vert\nabla\psi\vert^2}{2},\phi-\overline{\phi}\right)\\&\leq C_0^2\left(\Vert   \nabla\mu\Vert   \ \Vert   \nabla\phi\Vert  \right )-\alpha\Vert   \nabla\phi\Vert   ^2\\&+2\alpha_0\left(\Vert   \phi\Vert   ^2+\Vert   \phi-\overline{\phi}\Vert   ^2\right)+\Vert\nabla\psi\Vert^2_{[L^4(\Omega)]^4}\Vert   \phi-\overline{\phi}\Vert\\
&\leq C (1+\Vert   \nabla\mu\Vert+\Vert\psi\Vert_{[H^2(\Omega)]^2}   ).
\end{align*}
Notice that we have used the well-known property (see, e.g., \cite{Miranville}) \eqref{fprimo}, which is also valid for the logarithmic potential and not only for its regular approximation.
Therefore, we have (see again \eqref{lambdaest} and \eqref{energy3})
\begin{align*}
\vert \overline{\mu} \vert &\leq \frac{1}{ \vert \Omega \vert }\left(\int_\Omega \vert \Psi'(\phi) \vert dx+\int_\Omega \frac{\lambda_\alpha^\prime(\phi)}{2}\vert\nabla\psi\vert^2dx\right) \leq\frac{1}{ \vert \Omega \vert }\left(\int_\Omega \vert F'(\phi) \vert +\alpha_0\int_\Omega \vert \phi \vert +C\Vert\nabla\psi\Vert^2\right)\\&\leq\frac{C}{ \vert \Omega \vert }\left(\int_\Omega \vert F'(\phi) \vert +\alpha_0\sqrt{ \vert \Omega \vert }\ \Vert   \phi\Vert  +1\right )\leq\frac{C}{ \vert \Omega \vert }\left( \left\vert \int_\Omega F'(\phi)(\phi-\overline{\phi}) \right\vert + 1 \right)\\&\leq C\left(1+\Vert   \nabla\mu\Vert+\Vert\psi\Vert_{[H^2(\Omega)]^2}   \right).
\end{align*}
Therefore, by Poincar\'{e}'s inequality we have
\begin{align}
\Vert\mu\Vert_{V}\leq C\left(1+\Vert\nabla\mu\Vert+\Vert\psi\Vert_{[H^2(\Omega)]^2}\right).
\label{muH1}
\end{align}
Thus, we get (see \eqref{pot})
\begin{align}
\Vert\phi\Vert_{W^{2,p}(\Omega)} + \Vert F^\prime(\phi)\Vert_{L^p(\Omega)} \leq C\left(1+\Vert\nabla\mu\Vert+\Vert\psi\Vert_{[H^2(\Omega)]^2}+\Vert\psi\Vert^{2-2/p}_{[W^{2,p}(\Omega)]^2}\right).
\label{w2p2}
\end{align}
Multiplying \eqref{problem}$_4$ by $\partial_t\mu$ and integrating over $\Omega$, we obtain (see also \eqref{muAC})
\begin{equation}
\label{mut}
\frac{1}{2}\frac{d}{dt}\Vert\nabla\mu\Vert^2+<\partial_t\mu,\partial_t\phi>+\underbrace{<\partial_t\mu,\textbf{v}\cdot\nabla\phi>}_{K_1}=0.
\end{equation}
Observe that (see \eqref{problem}$_5$)
$$
<\partial_t\mu,\partial_t\phi>=\Vert\nabla\partial_t\phi\Vert^2+\underbrace{(\Psi^{\prime\prime}(\phi)\partial_t\phi,\partial_t\phi)}_{K_2}+\underbrace{(\nabla\psi,\nabla\partial_t\psi\lambda_\alpha^\prime(\phi)\partial_t\phi)}_{K_3}+\underbrace{(\frac{1}{2}\vert\nabla\psi\vert^2\lambda_\alpha^{\prime\prime}(\phi)\partial_t\phi,\partial_t\phi)}_{K_4},
$$
note that, due to the regularity of the functions involved, the dualities are actually scalar products in $H$.
Let us now rewrite $K_1$ as follows
$$
K_1=\frac{d}{dt}(\textbf{v}\cdot \nabla\phi,\mu)-\underbrace{(\partial_t\textbf{v}\cdot\nabla\phi,\mu)}_{K_5}-\underbrace{(\textbf{v}\cdot\nabla\partial_t\phi,\mu)}_{K_6}.
$$
Observe now that, by \eqref{lp}, \eqref{h1v} and \eqref{energy3} and by Gagliardo-Nirenberg's inequality, we have
\begin{align}
\nonumber\Vert\mu\Vert_{H^2(\Omega)}&\leq C\left(\Vert\partial_t\phi\Vert+\Vert\textbf{v}\cdot \nabla\phi\Vert\right)\leq C(\Vert\partial_t\phi\Vert+\Vert\textbf{u}\Vert_{[L^4(\Omega)]^2}\Vert\nabla\phi\Vert_{[L^4(\Omega)]^4})\\&\leq C(\Vert\partial_t\phi\Vert+\Vert D\textbf{u}\Vert^{1/2}\Vert\phi\Vert_{H^2(\Omega)}^{1/2})\leq C(\Vert\partial_t\phi\Vert+\Vert D\textbf{u}\Vert+\Vert\phi\Vert_{H^2(\Omega)}).
\label{muH2}
\end{align}
Hence, by \eqref{dtu}, \eqref{energy3}, Sobolev embeddings, Agmon's inequality, \eqref{muH1}, \eqref{w2p2} and \eqref{muH2}, we deduce
\begin{align*}
K_5&\leq \Vert\partial_t\textbf{v}\Vert\Vert\nabla\phi\Vert\Vert\mu\Vert_{L^\infty(\Omega)}\leq C\Vert\partial_t\textbf{u}\Vert\Vert\mu\Vert^{1/2}\Vert\mu\Vert^{1/2}_{H^2(\Omega)}\\
&\leq \delta \Vert\partial_t\textbf{u}\Vert^2+C(\delta)\Vert\mu\Vert_V^2+C(\delta)\Vert\mu\Vert^2_{H^2(\Omega)}\\&\leq
\delta \Vert\partial_t\textbf{u}\Vert^2+C(\delta)\left(1+\Vert\nabla\mu\Vert^2+\Vert\psi\Vert_{[H^2(\Omega)
]^2}^2\right)+C(\delta)(\Vert\partial_t\phi\Vert^2+\Vert D\textbf{u}\Vert^2+\Vert\phi\Vert_{H^2(\Omega)}^2)\\&\leq \delta \Vert\partial_t\textbf{u}\Vert^2+C(\delta)\left(1+\Vert\nabla\mu\Vert^2+\Vert\psi\Vert_{[H^2(\Omega)]^2}^2+\Vert D\textbf{u}\Vert^2\right)+C(\delta)\Vert\partial_t\phi\Vert^2,
\end{align*}
but from \eqref{phistar223} we have
\begin{align}
C(\delta)\Vert\partial_t\phi\Vert^2\leq \frac{1}{8}\Vert\partial_t\nabla\phi\Vert^2+C(\delta)\Vert\partial_t\phi\Vert_{V^\prime}^2\leq \frac{1}{8}\Vert\partial_t\nabla\phi\Vert^2
+C(\delta)\left(\Vert\nabla\mu\Vert^2+\Vert D\textbf{u}\Vert^2\right).
\label{dtphi}
\end{align}
Therefore, we get
\begin{align*}
&K_5\leq \delta \Vert\partial_t\textbf{u}\Vert^2+C(\delta)\left(1+\Vert\nabla\mu\Vert^2+\Vert\psi\Vert_{[H^2(\Omega)]^2}^2+\Vert D\textbf{u}\Vert^2\right)+\frac{1}{8}\Vert\partial_t\nabla\phi\Vert^2.
\end{align*}
Then, from \eqref{lp}, \eqref{energy3}, \eqref{muH1}, we deduce
\begin{align*}
K_6&\leq \Vert\mu\Vert_{L^4(\Omega)}\Vert\textbf{v}\Vert_{[L^4(\Omega)]^2}\Vert\nabla\partial_t\phi\Vert\leq C\Vert\mu\Vert_V\Vert D\textbf{u}\Vert^{1/2}\Vert\nabla\partial_t\phi\Vert\\&\leq \frac{1}{8}\Vert\nabla\partial_t\phi\Vert^2+C\Vert\mu\Vert_V^2\Vert D\textbf{u}\Vert\leq \frac{1}{8}\Vert\nabla\partial_t\phi\Vert^2+C\Vert D\textbf{u}\Vert\left(1+\Vert\nabla\mu\Vert^2+\Vert\psi\Vert_{[H^2(\Omega)]^2}^2\right)
\end{align*}
Observe now that (see \eqref{quasiconv})
$$
K_2\geq -C\Vert\partial_t\phi\Vert^2,
$$
which entails (see \eqref{phistar223} and \eqref{dtphi})
$$
K_2\geq -\frac{1}{8}\Vert\partial_t\nabla\phi\Vert^2-C\Vert\partial_t\phi\Vert_{V^\prime}^2\geq-\frac{1}{8}\Vert\partial_t\nabla\phi\Vert^2
- C\left(\Vert\nabla\mu\Vert^2+\Vert D\textbf{u}\Vert^2\right).
$$
Consider now \eqref{problem}$_3$. Taking the gradient, multiplying by $\partial_t\nabla\psi$ and integrating over $\Omega$, by \eqref{h1v}, \eqref{energy3} and
the Sobolev embeddings $\textbf{V}_\sigma\hookrightarrow [L^q(\Omega)]^2$, $2< q< \infty$, $[W^{2,p}(\Omega)]^2\hookrightarrow[W^{1,\infty}(\Omega)]^2$, $2<p< \infty$, we get
\begin{align}
\nonumber\Vert\nabla\partial_t\psi\Vert&\leq C\Vert\nabla(\textbf{v}\cdot\nabla\psi)\Vert \\&\leq C\Vert D\textbf{u}\Vert\left(\Vert\nabla\psi\Vert_{[L^\infty(\Omega)]^4}+\Vert\psi\Vert_{[W^{2,p}(\Omega)]^2}\right)\leq C\Vert D\textbf{u}\Vert\Vert\psi\Vert_{[W^{2,p}(\Omega)]^2}.
\label{psit3}
\end{align}
Then, again by Sobolev embeddings, Gagliardo-Nirenberg's inequality, \eqref{lambdaest}, \eqref{energy3} and \eqref{psit3}, we obtain
\begin{align*}
K_3&\leq \Vert\nabla\psi\Vert_{L^4(\Omega)}\Vert\nabla\partial_t\psi\Vert\Vert\partial_t\phi\Vert_{L^4(\Omega)}\leq C\Vert\psi\Vert_{[H^2(\Omega)]^2}^{1/2}\Vert\nabla\partial_t\psi\Vert\Vert\nabla\partial_t\phi\Vert \\&\leq
C\Vert D\textbf{u}\Vert\Vert\psi\Vert_{[W^{2,p}(\Omega)]^2}^{3/2}\Vert\nabla\partial_t\phi\Vert\\&
\leq \frac{1}{8}\Vert\nabla\partial_t\phi\Vert^2+C\Vert D\textbf{u}\Vert^2\Vert\psi\Vert^3_{[W^{2,p}(\Omega)]^2},
\end{align*}
whereas for $K_4$ it holds (see \eqref{lambdamin})
\begin{align*}
K_4=\left(\frac{1}{2}\vert\nabla\psi\vert^2\lambda_\alpha^{\prime\prime}(\phi),\vert\partial_t\phi\vert^2\right)\geq 0,\quad \text{ a.e. in  } (0,T_0).
\end{align*}
Taking the above estimates into account, we infer from \eqref{mut} that
\begin{align}
\nonumber&\frac{d}{dt}\left(\frac{1}{2}\Vert\nabla\mu\Vert^2+(\textbf{v}\cdot \nabla\phi,\mu)\right)+\frac{1}{2}\Vert\nabla\partial_t\phi\Vert^2\\&\leq \delta \Vert\partial_t\textbf{u}\Vert^2+C(\delta)\left(1+\Vert\nabla\mu\Vert^2+\Vert\nonumber\psi\Vert_{[H^2(\Omega)]^2}^2+\Vert D\textbf{u}\Vert^2\right)\nonumber\\&+C\Vert D\textbf{u}\Vert\left(1+\Vert\nabla\mu\Vert^2+\Vert\psi\Vert_{[H^2(\Omega)]^2}^2\right)
+C\Vert D\textbf{u}\Vert^2\Vert\psi\Vert^3_{[W^{2,p}(\Omega)]^2}.
\label{estmu}
\end{align}
On the other hand, due to Sobolev embeddings, \eqref{lp} and \eqref{energy3}, we have
\begin{align*}
\vert(\textbf{v}\cdot \nabla\phi,\mu)\vert=\vert(\textbf{v}\cdot\nabla\mu,\phi)\vert &\leq C \Vert\textbf{u}\Vert_{[L^4(\Omega)]^2}\Vert\nabla\mu\Vert\Vert\phi\Vert_{L^4(\Omega)}\leq C\Vert D\textbf{u}\Vert^{1/2}\Vert\nabla\mu\Vert\\&\leq  C+\frac{1}{4}\Vert D\textbf{u}\Vert^2+\frac{1}{4}\Vert\nabla\mu\Vert^2,
\end{align*}
so that
\begin{align}
(\textbf{v}\cdot \nabla\phi,\mu)\geq -\vert(\textbf{v}\cdot \nabla\phi,\mu)\vert\geq - C-\frac{1}{4}\Vert D\textbf{u}\Vert^2-\frac{1}{4}\Vert\nabla\mu\Vert^2.
\label{vcontrol}
\end{align}
\textbf{Fourth estimate.}
Observe that (see \eqref{regularities} and \eqref{newreg})
\begin{align}
\Vert \textbf{u}\Vert_{L^\infty(0,T_0;\textbf{V}_\sigma^2)}\leq C_\alpha.
\label{reg}
\end{align}
Therefore, by \eqref{regularities}, \eqref{newreg}, \eqref{C1} and \eqref{reg}, we have that (see, e.g., \cite[Lemma II.5.9]{Boyerlibro}),
\begin{align}
\textbf{u}\in C^1([0,T_0], \textbf{H}_\sigma)\cap C_w([0,T_0];\textbf{V}_\sigma^2).
\label{ureg}
\end{align}
This regularity allows us to differentiate with respect to time the weak formulation of \eqref{problem}$_1$ and test it with $\textbf{A}^{-1/2-\varepsilon}\partial_t\textbf{u}(t)$, for $\varepsilon\in(0,1/6)$ (see \ref{Hmedium}). To be precise, this step should be performed by means of the increments $\partial_t^h\textbf{u}$ as in \eqref{partialta}, to which we refer for the rigorous details. Indeed, notice that also for the initial conditions, thanks to \eqref{ureg}, we have that $\partial_t^h\textbf{u}(0)\rightarrow \partial_t\textbf{u}(0)$ as $h\to0$ in $\textbf{H}_\sigma$. Therefore, for the sake of brevity, we can argue formally without considering the increments. We obtain
\begin{align*}
&\frac{1}{2}\frac{d}{dt}\Vert \partial_t\textbf{u}\Vert^2_{\textbf{V}_\sigma^{-1/2-\varepsilon}}+\underbrace{((\partial_t\textbf{v}\cdot\nabla)\textbf{u},\textbf{A}^{-1/2-\varepsilon}\partial_t\textbf{u})}_{J_1}-\underbrace{(\lambda_\alpha^\prime(\phi)\partial_t\phi\nabla\psi^T\nabla\psi,\nabla(I+\alpha\textbf{A})^{-2}\textbf{A}^{-1/2-\varepsilon}\partial_t\textbf{u})}_{J_2}\\&-\underbrace{\left((\lambda_\alpha(\phi)\nabla\partial_t\psi^T\nabla\psi,\nabla(I+\alpha\textbf{A})^{-2}\textbf{A}^{-1/2-\varepsilon}\partial_t\textbf{u})-(\lambda_\alpha(\phi)\nabla\psi^T\nabla\partial_t\psi,\nabla(I+\alpha\textbf{A})^{-2}\textbf{A}^{-1/2-\varepsilon}\partial_t\textbf{u})\right)}_{J_3}\\&-\underbrace{\left(\frac{\nu^\prime(\phi)}{2k(\phi)}\partial_t\phi\textbf{u},\textbf{A}^{-1/2-\varepsilon}\partial_t\textbf{u}\right)}_{J_4}+\underbrace{\left(\frac{\nu(\phi)}{2k(\phi)}(1-\phi)\partial_t\textbf{u},\textbf{A}^{-1/2-\varepsilon}\partial_t\textbf{u}\right)}_{J_5}\\&+\underbrace{\left(\frac{\nu^\prime(\phi)}{2k(\phi)}(1-\phi)\partial_t\phi\textbf{u},\textbf{A}^{-1/2-\varepsilon}\partial_t\textbf{u}\right)}_{J_6}-\underbrace{\left(\frac{\nu(\phi)k^\prime(\phi)}{2k^2(\phi)}\partial_t\phi(1-\phi)\textbf{u},\textbf{A}^{-1/2-\varepsilon}\partial_t\textbf{u}\right)}_{J_7}\\&-\underbrace{(\nabla\partial_t\phi\otimes\nabla\phi,\nabla\textbf{A}^{-1/2-\varepsilon}\partial_t\textbf{u})}_{J_8}-\underbrace{(\nabla\phi\otimes\nabla\partial_t\phi,\nabla\textbf{A}^{-1/2-\varepsilon}\partial_t\textbf{u})}_{J_9}\\&-\underbrace{(\text{div}(\nu(\varphi) D\partial_t\textbf{u}),\textbf{A}^{-1/2-\varepsilon}\partial_t\textbf{u})}_{J_{10}}-\underbrace{(\text{div}(\nu^\prime(\varphi)\partial_t\varphi D\textbf{u}),\textbf{A}^{-1/2-\varepsilon}\partial_t\textbf{u})}_{J_{11}}=0.
\end{align*}
Consider $J_1$. Using H\"older's inequality and Sobolev embeddings, \eqref{lp}, \eqref{interp3} and \eqref{interp4}, we find
\begin{align*}
\vert J_1\vert &\leq C\Vert\partial_t\textbf{u}\Vert_{[L^{\frac{4}{1-2\varepsilon}}(\Omega)]^2}\Vert\nabla\textbf{u}\Vert\Vert\textbf{A}^{-1/2-\varepsilon}\partial_t\textbf{u}\Vert_{[L^\frac{4}{1+2\varepsilon}(\Omega)]^2}\\
&\leq C\Vert\partial_t\textbf{u}\Vert_{\textbf{V}_\sigma^{1/2-\varepsilon}}\Vert\nabla\textbf{u}\Vert\Vert\partial_t\textbf{u}\Vert_{\textbf{V}_\sigma^{-1/2-3\varepsilon}}\\&\leq {\delta_3}\Vert\partial_t\textbf{u}\Vert_{\textbf{V}_\sigma^{1/2-\varepsilon}}^2+C(\delta_1)\Vert\nabla\textbf{u}\Vert^2\Vert\partial_t\textbf{u}\Vert_{\textbf{V}_\sigma^{-1/2-\varepsilon}}^2,
\end{align*}
for $\delta_3>0$.
Concerning $J_2$, since $\overline{\partial_t\phi}\equiv0$, by interpolation we get
\begin{equation}
\Vert\partial_t\phi\Vert\leq C\Vert\nabla\partial_t\phi\Vert^{1/2}\Vert\partial_t\phi\Vert_{V^\prime}^{1/2}.
\label{phistar}
\end{equation}
Therefore, being $\lambda_\alpha^\prime$ bounded on $[0,1]$, recalling the conservation of mass, \eqref{h1v} and \eqref{phistar223}, we infer
\begin{align*}
\vert J_2\vert&\leq C\Vert \partial_t\phi\Vert\Vert\nabla\psi\Vert_{[L^\infty(\Omega)]^2}^2\Vert (I+\alpha\textbf{A})^{-2}\textbf{A}^{-1/2-\varepsilon}\partial_t\textbf{u}\Vert_{\textbf{V}_\sigma}\\&\leq C\Vert \partial_t\phi\Vert\Vert\nabla\psi\Vert_{[L^\infty(\Omega)]^4}^2\Vert \textbf{A}^{-1/2-\varepsilon}\partial_t\textbf{u}\Vert_{\textbf{V}_\sigma}\\&\leq C\Vert\nabla\partial_t\phi\Vert^{1/2}\Vert\partial_t\phi\Vert_{V^\prime}^{1/2}\Vert\psi\Vert_{[W^{2,p}(\Omega)]^2}^2\Vert \partial_t\textbf{u}\Vert\leq \delta_2\Vert\nabla \partial_t\phi\Vert^2\\&+\delta_1\Vert \partial_t\textbf{u}\Vert^2+C(\delta_1,\delta_2)\Vert\partial_t\phi\Vert_{V^\prime}^2\Vert\psi\Vert_{[W^{2,p}(\Omega)]^2}^8\\&\leq \delta_2\Vert\nabla \partial_t\phi\Vert^2+\delta_1\Vert \partial_t\textbf{u}\Vert^2+C(\delta_1,\delta_2)(1+\Vert\nabla\mu\Vert^2)\Vert\psi\Vert_{[W^{2,p}(\Omega)]^2}^8,
\end{align*}
for $\delta_2>0$.
On the other hand, by \eqref{h1v} together with \eqref{psit3}, we have
\begin{align*}
\vert J_3\vert&\leq C\Vert\nabla\partial_t\psi\Vert\Vert\nabla\psi\Vert_{[L^\infty(\Omega)]^4}\Vert (I+\alpha\textbf{A})^{-2}\textbf{A}^{-1/2-\varepsilon}\partial_t\textbf{u}\Vert_{\textbf{V}_\sigma}\\&\leq C(\delta_1)\Vert\nabla\partial_t\psi\Vert^2\Vert\nabla\psi\Vert_{[L^\infty(\Omega)]^4}^2+\delta \Vert \partial_t\textbf{u}\Vert^2\\&\leq C(\delta_1)\Vert D\textbf{u}\Vert^2\Vert\psi\Vert_{[W^{2,p}(\Omega)]^2}^4+\delta_1 \Vert \partial_t\textbf{u}\Vert^2.
\end{align*}
For $J_4$, by \eqref{lp}, \eqref{interp1}, \eqref{interp4} and \eqref{energy3}, we deduce
\begin{align*}
\vert J_4\vert&\leq C\Vert\textbf{u}\Vert\Vert\partial_t\phi\Vert_{L^{\frac{4}{1-2\varepsilon}}(\Omega)}\Vert\textbf{A}^{-1/2-\varepsilon}\partial_t\textbf{u}\Vert_{[L^{\frac{4}{1+2\varepsilon}}(\Omega)]^2}\\&\leq C\Vert\nabla\partial_t\phi\Vert\Vert \partial_t\textbf{u}\Vert_{\textbf{V}_\sigma^{-1/2-\varepsilon}}\leq
\delta_2\Vert\nabla\partial_t\phi\Vert^2+C(\delta_2)\Vert\partial_t\textbf{u}\Vert^{2}_{\textbf{V}_\sigma^{-1/2-\varepsilon}}.
\end{align*}
By the embedding $\textbf{V}_\sigma^{-1-2\varepsilon}\hookrightarrow \textbf{V}_\sigma^{-1/2-\varepsilon}$, we obtain
$$
\vert J_5\vert \leq C\Vert\partial_t\textbf{u}\Vert\Vert\textbf{A}^{-1/2-\varepsilon}\partial_t\textbf{u}\Vert\leq {\delta}\Vert\partial_t\textbf{u}\Vert^2+C(\delta)\Vert\partial_t\textbf{u}\Vert_{\textbf{V}_\sigma^{-1/2-\varepsilon}}^2,
$$
and, similarly, by \eqref{energy3} and Sobolev embeddings, being $\nu,k\in W^{1,\infty}(\R)$,
\begin{align*}
\vert J_6\vert +\vert J_7\vert&\leq C\Vert\partial_t\phi\Vert_{L^4(\Omega)}\Vert\textbf{u}\Vert_{[L^4(\Omega)]^2}\Vert\textbf{A}^{-1/2-\varepsilon}\partial_t\textbf{u}\Vert\\&\leq C\Vert\nabla\partial_t\phi\Vert\Vert\nabla\textbf{u}\Vert^{1/2}\Vert\partial_t\textbf{u}\Vert_{\textbf{V}_\sigma^{-1/2-\varepsilon}}\leq \delta_2\Vert\nabla\partial_t\phi\Vert^2+C(\delta_2)\Vert\nabla\textbf{u}\Vert\Vert\partial_t\textbf{u}\Vert_{\textbf{V}_\sigma^{-1/2-\varepsilon}}^2.
\end{align*}
Then, thanks to \eqref{energy3}, Sobolev-Gagliardo-Nirenberg's and Young's inequalities and \eqref{l3}, we deduce
\begin{align*}
\vert J_8\vert +\vert J_9\vert&\leq \Vert\nabla\phi\Vert_{[L^6(\Omega)]^2}\Vert\nabla\partial_t\phi\Vert\Vert \nabla(I+\alpha\textbf{A})^{-2}\textbf{A}^{-1/2}\partial_t\textbf{u}\Vert_{[L^3(\Omega)]^4}\\&\leq C\Vert\nabla\phi\Vert_{[L^6(\Omega)]^2}\Vert\nabla\partial_t\phi\Vert \Vert\partial_t\textbf{u}\Vert_{\textbf{V}_\sigma^{1/2-\varepsilon}}^{\frac{5-12\varepsilon^2-4\varepsilon}{2(3-6\varepsilon)}}\Vert\partial_t\textbf{u}\Vert_{\textbf{V}_\sigma^{-1/2-\varepsilon}}^{\frac{1+12\varepsilon^2-8\varepsilon}{2(3-6\varepsilon)}}\\&\leq C\Vert\nabla\phi\Vert^{1/3}\Vert\phi\Vert_{H^2(\Omega)}^{2/3}\Vert\nabla\partial_t\phi\Vert\Vert\partial_t\textbf{u}\Vert_{\textbf{V}_\sigma^{1/2-\varepsilon}}^{\frac{5-12\varepsilon^2-4\varepsilon}{2(3-6\varepsilon)}}\Vert\partial_t\textbf{u}\Vert_{\textbf{V}_\sigma^{-1/2-\varepsilon}}^{\frac{1+12\varepsilon^2-8\varepsilon}{2(3-6\varepsilon)}}\\&\leq \delta_2\Vert\nabla\partial_t\phi\Vert^2+\delta_1\Vert\partial_t\textbf{u}\Vert_{\textbf{V}_\sigma^{1/2-\varepsilon}}^2+C(\delta_1,\delta_2)\Vert\phi\Vert_{H^2(\Omega)}^{\frac{8}{1-6\varepsilon}}\Vert\partial_t\textbf{u}\Vert^{2}_{\textbf{V}_\sigma^{-1/2-\varepsilon}}.
\end{align*}
We now rewrite $J_{10}$ as follows
\begin{align*}
J_{10}&=(\nu(\varphi)D\partial_t\textbf{u}, \nabla\textbf{A}^{-1/2-\varepsilon}\partial_t\textbf{u})=(\nu(\varphi)\nabla\partial_t\textbf{u}, D\textbf{A}^{-1/2-\varepsilon}\partial_t\textbf{u})\\&=-(\partial_t\textbf{u},\text{div}(\nu(\varphi)D\textbf{A}^{-1/2-\varepsilon}\partial_t\textbf{u}))\\&= -(\partial_t\textbf{u},\nu^\prime(\varphi)D\textbf{A}^{-1/2-\varepsilon}\partial_t\textbf{u}\nabla\varphi)-\frac{1}{2}(\partial_t\textbf{u},\nu(\varphi)\Delta\textbf{A}^{-1/2-\varepsilon}\partial_t\textbf{u})\\&=-(\partial_t\textbf{u},\nu^\prime(\varphi)D\textbf{A}^{-1/2-\varepsilon}\partial_t\textbf{u}\nabla\varphi)+\frac{1}{2}(\partial_t\textbf{u},\nu(\varphi)\textbf{A}^{1/2-\varepsilon}\partial_t\textbf{u})-\frac{1}{2}(\partial_t\textbf{u},\nu(\varphi)\nabla \hat{\pi})\\&=-\underbrace{(\partial_t\textbf{u},\nu^\prime(\varphi)D\textbf{A}^{-1/2-\varepsilon}\partial_t\textbf{u}\nabla\varphi)}_{J_{10,1}}+{\frac{1}{2}(\partial_t\textbf{u},\nu(\varphi)\textbf{A}^{1/2-\varepsilon}\partial_t\textbf{u})}+\underbrace{\frac{1}{2}(\nu^\prime(\varphi)\partial_t\textbf{u},\hat{\pi}\nabla\varphi )}_{J_{10,2}},
\end{align*}
where $\hat{\pi}$ is given by $-\Delta\textbf{A}^{-1/2-\varepsilon}\partial_t\textbf{u}+\nabla \pi=\textbf{A}^{1/2-\varepsilon}\partial_t\textbf{u}$, therefore, by Lemma \ref{press},
\begin{align}
\nonumber
\Vert \hat{\pi}\Vert &\leq C\left(\Vert\textbf{A}^{1/2-\varepsilon}\partial_t\textbf{u}\Vert_{\textbf{V}_\sigma^{-1}}+\Vert\textbf{A}^{1/2-\varepsilon}\partial_t\textbf{u}\Vert_{\textbf{V}_\sigma^{-1/2+\varepsilon}}\right)\\
&\leq C\left(\Vert\partial_t\textbf{u}\Vert_{\textbf{V}_\sigma^{-2\varepsilon}}+\Vert\partial_t\textbf{u}\Vert_{\textbf{V}_\sigma^{1/2-\varepsilon}}\right).
\label{p}
\end{align}
Moreover, notice that, by \eqref{nu1}, we infer
$$
\frac{1}{2}(\partial_t\textbf{u},\nu(\varphi)\textbf{A}^{1/2-\varepsilon}\partial_t\textbf{u})=\frac{\nu_1+\nu_2}{4}\Vert\partial_t\textbf{u}\Vert^2_{\textbf{V}_\sigma^{1/2-\varepsilon}}+\underbrace{\frac{\nu_1-\nu_2}{4}(\partial_t\textbf{u},\varphi\textbf{A}^{1/2-\varepsilon}\partial_t\textbf{u})}_{J_{10,3}}.
$$
where we used the fact that $\Vert \varphi \Vert_{L^\infty(\Omega\times (0,T_0))}\leq 1$ (see \eqref{limit} and \eqref{regularities}).

Consider now $J_{10,1}$. In this case,  by \eqref{interp1}, being $\Vert D\textbf{A}^{-1/2-\varepsilon}\partial_t\textbf{u}\Vert\leq C\Vert\partial_t\textbf{u}\Vert$, and by \eqref{lpphi2}, we have
\begin{align*}
\vert J_{10,1}\vert&\leq C\Vert\partial_t\textbf{u}\Vert\Vert D\textbf{A}^{-1/2-\varepsilon}\partial_t\textbf{u}\Vert\Vert\nabla\varphi\Vert_{[L^\infty(\Omega)]^2}\leq C\Vert\partial_t\textbf{u}\Vert^2\Vert\nabla\varphi\Vert_{[L^\infty(\Omega)]^2}\\&\leq C\Vert\partial_t\textbf{u}\Vert^{1-2\varepsilon}_{\textbf{V}_\sigma^{-1/2-\varepsilon}}  \Vert\partial_t\textbf{u}\Vert^{1+2\varepsilon}_{\textbf{V}_\sigma^{1/2-\varepsilon}}\Vert\nabla\varphi\Vert_{[L^\infty(\Omega)]^2}\\&\leq {\delta_3}\Vert\partial_t\textbf{u}\Vert_{\textbf{V}_\sigma^{1/2-\varepsilon}}^2+C(\delta_3)\Vert\partial_t\textbf{u}\Vert_{\textbf{V}_\sigma^{-1/2-\varepsilon}}^2\Vert\nabla\varphi\Vert_{[L^\infty(\Omega)]^2}^{\frac{2}{1-2\varepsilon}}\\&\leq {\delta_3}\Vert\partial_t\textbf{u}\Vert_{\textbf{V}_\sigma^{1/2-\varepsilon}}^2+C(\delta_3)\Vert\partial_t\textbf{u}\Vert_{\textbf{V}_\sigma^{-1/2-\varepsilon}}^2\Vert\varphi\Vert_{W^{2,p}(\Omega)}^{\frac{2}{1-2\varepsilon}}\\&\leq {\delta_3}\Vert\partial_t\textbf{u}\Vert_{\textbf{V}_\sigma^{1/2-\varepsilon}}^2+C(\delta_3)\Vert\partial_t\textbf{u}\Vert_{\textbf{V}_\sigma^{-1/2-\varepsilon}}^2\Vert\phi\Vert_{W^{2,p}(\Omega)}^{\frac{2}{1-2\varepsilon}},
\end{align*}
for $p>2$. Then, by the same inequalities, by \eqref{interp5} and \eqref{p}, we infer, for $p>2$,
\begin{align*}
\vert J_{10,2}\vert&\leq C\Vert\partial_t\textbf{u}\Vert \Vert \hat{\pi}\Vert \Vert \nabla\varphi\Vert_{[L^\infty(\Omega)]^2}\leq C\Vert\partial_t\textbf{u}\Vert\left(\Vert\partial_t\textbf{u}\Vert_{\textbf{V}_\sigma^{-2\varepsilon}}+\Vert\partial_t\textbf{u}\Vert_{\textbf{V}_\sigma^{1/2-\varepsilon}}\right)\Vert \nabla\varphi\Vert_{[L^\infty(\Omega)]^2}\\&\leq C\Vert\partial_t\textbf{u}\Vert_{\textbf{V}_\sigma^{-1/2-\varepsilon}}^{1/2-\varepsilon}\Vert\partial_t\textbf{u}\Vert_{\textbf{V}_\sigma^{1/2-\varepsilon}}^{1/2+\varepsilon}\left(\Vert\partial_t\textbf{u}\Vert_{\textbf{V}_\sigma^{-1/2-\varepsilon}}^{1/2+\varepsilon}\Vert\partial_t\textbf{u}\Vert_{\textbf{V}_\sigma^{1/2-\varepsilon}}^{1/2-\varepsilon}+\Vert\partial_t\textbf{u}\Vert_{\textbf{V}_\sigma^{1/2-\varepsilon}}\right)\Vert \varphi\Vert_{W^{2,p}(\Omega)}\\&\leq \delta_3\Vert\partial_t\textbf{u}\Vert_{\textbf{V}_\sigma^{1/2-\varepsilon}}^{2}+C(\delta_3)\left(\Vert \phi\Vert_{W^{2,p}(\Omega)}^{2}+\Vert \phi\Vert_{W^{2,p}(\Omega)}^{\frac{4}{1-2\varepsilon}}\right)\Vert\partial_t\textbf{u}\Vert_{\textbf{V}_\sigma^{-1/2-\varepsilon}}^{2}.
\end{align*}
We now recall the following essential estimate  (see \cite[Lemma 4.11]{Maxwell}):
\begin{align}
\Vert \partial_t\textbf{u}\varphi\Vert_{[H^{1/2-\varepsilon}(\Omega)]^2}\leq {C}_1\left(\Vert \varphi\Vert_{L^\infty(\Omega)}\Vert \partial_t\textbf{u}\Vert_{[H^{1/2-\varepsilon}(\Omega)]^2}+\Vert \varphi\Vert_{H^2(\Omega)}\Vert\partial_t\textbf{u}\Vert\right),
\label{Maxwell}
\end{align}
for some ${C}_1>0$. Moreover, by the regularity results on the projector $\textbf{P}$ (see, e.g., \cite[Ch.1]{temam}) we have that, for any $\textbf{w}\in [H^{1/2-\varepsilon}(\Omega)]^2$,
\begin{align}
\Vert\textbf{P}\textbf{w}\Vert_{[H^{1/2-\varepsilon}(\Omega)]^2}\leq C_2\Vert\textbf{w}\Vert_{[H^{1/2-\varepsilon}(\Omega)]^2},
\end{align}
for some ${C}_2>0$.
Observe now that $\textbf{P}(\partial_t\textbf{u}\varphi)\in D(\textbf{A}^{1/4-\varepsilon/2})=\textbf{V}_\sigma^{1/2-\varepsilon}$ for almost any $t$. Indeed, by \cite[Corollary 2.1]{salgado}, we know that, being $\varepsilon\in(0,1/6)$, $\textbf{V}_\sigma^{1/2-\varepsilon}=[H^{1/2-\varepsilon}(\Omega)]^2\cap \textbf{H}_\sigma$. Since, by \eqref{regularities} and \eqref{newreg}, $\partial_t\textbf{u}\varphi\in [H^1(\Omega)]^2$ for almost any $t$, by \cite[Ch.1]{temam}, $\textbf{P}(\partial_t\textbf{u}\varphi)\in [H^1(\Omega)]^2\cap \textbf{H}_\sigma$, which implies that $\textbf{P}(\partial_t\textbf{u}\varphi)\in \textbf{V}_\sigma^{1/2-\varepsilon}$. Note that the case $\varepsilon=0$ is not allowed, since in general we have only $\textbf{V}_\sigma^{1/2}\subsetneqq[H^{1/2}(\Omega)]^2\cap \textbf{H}_\sigma$, being $\textbf{V}_\sigma^{1/2}=[H^{1/2}_{00}(\Omega)]^2\cap \textbf{H}_\sigma$ (see, e.g., \cite{LM}), and thus we are not able to ensure the result.

On the other hand, recalling as before that $\Vert \varphi \Vert_{L^\infty(\Omega\times (0,T_0))}\leq 1$, by the equivalence of the ${H}^{1/2-\varepsilon}$ norms over $D(\textbf{A}^{1/4-\varepsilon/2})$, by \eqref{lpphi2}, \eqref{interp1} and \eqref{Maxwell}, we find
\begin{align*}
\vert J_{10,3}\vert&=\frac{\nu_1-\nu_2}{4}(\partial_t\textbf{u}\varphi,\textbf{A}^{1/2-\varepsilon}\partial_t\textbf{u})=\frac{\nu_1-\nu_2}{4}(\textbf{P}(\partial_t\textbf{u}\varphi),\textbf{A}^{1/2-\varepsilon}\partial_t\textbf{u})\\
&= \frac{\nu_1-\nu_2}{4}(\textbf{A}^{1/4-\varepsilon/2}\textbf{P}(\partial_t\textbf{u}\varphi),\textbf{A}^{1/4-\varepsilon/2}\partial_t\textbf{u})\\
&\leq \frac{\vert \nu_1-\nu_2\vert }{4}\Vert \textbf{A}^{1/4-\varepsilon/2}\textbf{P}(\partial_t\textbf{u}\varphi)\Vert\Vert\textbf{A}^{1/4-\varepsilon/2}\partial_t\textbf{u}\Vert\\&\leq \frac{C_2\vert \nu_1-\nu_2\vert }{4}\Vert\partial_t\textbf{u}\varphi\Vert_{[{H}^{1/2-\varepsilon(\Omega)]^2}(\Omega)}\Vert\textbf{A}^{1/4-\varepsilon/2}\partial_t\textbf{u}\Vert\\&\leq \frac{C_1C_2\vert \nu_1-\nu_2\vert }{4}\left(\Vert \varphi\Vert_{L^\infty(\Omega)}\Vert \partial_t\textbf{u}\Vert_{[H^{1/2-\varepsilon}(\Omega)]^2}+\Vert \varphi\Vert_{H^2(\Omega)}\Vert \partial_t\textbf{u}\Vert\right)\Vert \partial_t\textbf{u}\Vert_{\textbf{V}^{1/2-\varepsilon}_\sigma}\\&\leq  \frac{C_3\vert \nu_1-\nu_2\vert }{4}\left(\Vert \partial_t\textbf{u}\Vert_{\textbf{V}^{1/2-\varepsilon}_\sigma(\Omega)}+\Vert \phi\Vert_{H^2(\Omega)}\Vert \partial_t\textbf{u}\Vert\right)\Vert \partial_t\textbf{u}\Vert_{\textbf{V}^{1/2-\varepsilon}_\sigma}\\&\leq
\frac{C_3\vert \nu_1-\nu_2\vert }{4}\Vert \partial_t\textbf{u}\Vert_{\textbf{V}^{1/2-\varepsilon}_\sigma}^2+C\Vert \phi\Vert_{H^2(\Omega)}\Vert \partial_t\textbf{u}\Vert_{\textbf{V}^{-1/2-\varepsilon}_\sigma}^{1/2-\varepsilon}\Vert \partial_t\textbf{u}\Vert_{\textbf{V}^{1/2-\varepsilon}_\sigma}^{3/2+\varepsilon}\\&\leq
\frac{C_3\vert \nu_1-\nu_2\vert }{4}\Vert \partial_t\textbf{u}\Vert_{\textbf{V}^{1/2-\varepsilon}_\sigma}^2+\delta_3\Vert \partial_t\textbf{u}\Vert_{\textbf{V}^{1/2-\varepsilon}_\sigma}^2+C \Vert \partial_t\textbf{u}\Vert_{\textbf{V}^{-1/2-\varepsilon}_\sigma}^2\Vert \phi\Vert_{H^2(\Omega)}^\frac{4}{1-2\varepsilon},
\end{align*}
for some $C_3>0$ depending only on $\Omega$.
Consider now $J_{11}$. Recalling \eqref{lpphi} and using interpolation, we get
\begin{align*}
\Vert\partial_t\varphi\Vert_{L^{\frac{4}{1-2\varepsilon}}(\Omega)}&\leq C\Vert\partial_t\phi\Vert_{L^{\frac{4}{1-2\varepsilon}}(\Omega)}\\
&\leq C\Vert\partial_t\phi\Vert^{1/2-\varepsilon}\Vert\nabla\partial_t\phi\Vert^{1/2+\varepsilon}\leq C\Vert\partial_t\phi\Vert^{1/4-\varepsilon/2}_{V^\prime}\Vert\nabla\partial_t\phi\Vert^{3/4+\varepsilon/2}.
\end{align*}
This inequality, \eqref{interp4bis3} and $\nu\in W^{1,\infty}(\R)$ allow us to deduce (see also \eqref{phistar223})
\begin{align*}
\vert J_{11}\vert&=(\nu^\prime(\varphi)\partial_t\varphi D\textbf{u},\nabla\textbf{A}^{-1/2-\varepsilon}\partial_t\textbf{u})\leq C\Vert\partial_t\varphi\Vert_{L^{\frac{4}{1-2\varepsilon}}(\Omega)}\Vert D\textbf{u}\Vert\Vert\nabla\textbf{A}^{-1/2-\varepsilon}\partial_t\textbf{u}\Vert_{L^{\frac{4}{1+2\varepsilon}}(\Omega)}\\&\leq C\Vert\partial_t\phi\Vert_{V^\prime}^{1/4-\varepsilon/2}\Vert\nabla\partial_t\phi\Vert^{3/4+\varepsilon/2}\Vert D\textbf{u}\Vert\Vert \partial_t\textbf{u}\Vert_{\textbf{V}^{1/2-\varepsilon}_\sigma}\\&\leq \delta_3 \Vert \partial_t\textbf{u}\Vert_{\textbf{V}^{1/2-\varepsilon}_\sigma}^2+\delta_2\Vert\nabla\partial_t\phi\Vert^2+C(\delta_2,\delta_3)\Vert\partial_t\phi\Vert_{V^\prime}^2\Vert D\textbf{u}\Vert^\frac{8}{1-2\varepsilon}\\&\leq \delta_3 \Vert \partial_t\textbf{u}\Vert_{\textbf{V}^{1/2-\varepsilon}_\sigma}^2+\delta_2\Vert\nabla\partial_t\phi\Vert^2+C(\delta_2,\delta_3)(1+\Vert\nabla\mu\Vert^2)\Vert D\textbf{u}\Vert^\frac{8}{1-2\varepsilon}.
\end{align*}
Therefore, choosing $\delta_1,\delta_2,\delta_3$ sufficiently small, we get in the end, recalling that $\varepsilon\in(0,1/6)$,
\begin{align}
\nonumber&\frac{1}{2}\frac{d}{dt}\Vert \partial_t\textbf{u}\Vert^2_{\textbf{V}_\sigma^{-1/2-\varepsilon}}+\Vert \partial_t\textbf{u}\Vert_{\textbf{V}^{1/2-\varepsilon}_\sigma}^2\left(\frac{\nu_1+\nu_2}{6}-\frac{C_3\vert \nu_1-\nu_2\vert }{4}\right)\\&\leq C\left(1+\Vert\nabla\textbf{u}\Vert^2+\Vert\phi\Vert_{H^2(\Omega)}^{\frac{8}{1-6\varepsilon}}+\Vert\phi\Vert_{W^{2,p}(\Omega)}^{\frac{2}{1-2\varepsilon}}+\Vert\phi\Vert_{W^{2,p}(\Omega)}^{\frac{4}{1-2\varepsilon}}+\Vert\phi\Vert_{W^{2,p}(\Omega)}^{\frac{4}{2-3\varepsilon}}\right)\Vert \partial_t\textbf{u}\Vert^2_{\textbf{V}_\sigma^{-1/2-\varepsilon}}\nonumber\\&
+C\left((1+\Vert\nabla\mu\Vert^2)\Vert\psi\Vert_{[W^{2,p}(\Omega)]^2}^8+\Vert D\textbf{u}\Vert^2\Vert\psi\Vert_{[W^{2,p}(\Omega)]^2}^4+(1+\Vert\nabla\mu\Vert^2)\Vert D\textbf{u}\Vert^\frac{8}{1-2\varepsilon}\right)
\nonumber\\
&+\frac{1}{8}\Vert\nabla\partial_t\phi\Vert^2 \leq C\left(1+\Vert\nabla\textbf{u}\Vert^2+\Vert\phi\Vert_{H^2(\Omega)}^{\frac{8}{1-6\varepsilon}}\right)\Vert \partial_t\textbf{u}\Vert^2_{\textbf{V}_\sigma^{-1/2-\varepsilon}}\nonumber\\&+C\left((1+\Vert\nabla\mu\Vert^2)\Vert\psi\Vert_{[W^{2,p}(\Omega)]^2}^8+\Vert D\textbf{u}\Vert^2\Vert\psi\Vert_{[W^{2,p}(\Omega)]^2}^4+(1+\Vert\nabla\mu\Vert^2)\Vert D\textbf{u}\Vert^{12}\right)\nonumber\\
&+\frac{1}{8}\Vert\nabla\partial_t\phi\Vert^2.
\label{estu}
\end{align}
\textbf{Fifth estimate.}
Using Theorem \ref{stokes} once more, with $p=s$ and $r=\infty$ for any $p\in\left(2,\frac{4}{1+2\varepsilon}\right]$, we obtain
\begin{align*}
\Vert\textbf{u}\Vert_{[W^{2,p}(\Omega)]^2}
&\leq \Vert\partial_t\textbf{u}\Vert_{[L^p(\Omega)]^2}+\Vert(\textbf{v}\cdot \nabla)\textbf{u}\Vert_{[L^p(\Omega)]^2}+\Vert(I+\alpha\textbf{A})^{-2}\textbf{P}\text{div}(\nabla\phi\otimes\nabla\phi)\Vert_{[L^p(\Omega)]^2}\\
&+\Vert(I+\alpha\textbf{A})^{-2}\textbf{P}\text{div}(\lambda(\phi)\nabla\psi^T\nabla\psi)\Vert_{[L^p(\Omega)]^2}+\left\Vert\frac{\nu(\phi)}{2k(\phi)}(1-\phi)\textbf{u}\right\Vert_{[L^p(\Omega)]^2}\\
&+\Vert\nabla\varphi\Vert_{[L^\infty(\Omega)]^2}\Vert D\textbf{u}\Vert_{[L^p(\Omega)]^4}.
\end{align*}
On the other hand, by \eqref{lp}, \eqref{energy3} and Sobolev embedding $W^{2,p}(\Omega)\hookrightarrow W^{1,\infty}(\Omega)$, we have
$$
\Vert(I+\alpha\textbf{A})^{-2}\textbf{P}\text{div}(\nabla\phi\otimes\nabla\phi)\Vert_{[L^p(\Omega)]^2}\leq C\Vert\nabla\phi\Vert_{[L^\infty(\Omega)]^2}\Vert D^2\phi\Vert_{[L^p(\Omega)]^4}\leq C\Vert\phi\Vert_{W^{2,p}(\Omega)}^2.
$$
Thus, by the embedding $H^2(\Omega)\hookrightarrow {W^{1,p}(\Omega)}$ and by \eqref{lp}, we infer
\begin{align*}
&\Vert(\textbf{v}\cdot \nabla)\textbf{u}\Vert_{[L^p(\Omega)]^2}\leq C\Vert\textbf{u}\Vert_{L^{2p}(\Omega)}\Vert\nabla\textbf{u}\Vert_{[L^{2p}]^2}\leq C\Vert D\textbf{u}\Vert\Vert D\textbf{u}\Vert^{1/p}\Vert\textbf{u}\Vert_{[H^2(\Omega)]^2}^{1-1/p}\\&\leq C\Vert D\textbf{u}\Vert^{1+1/p}\Vert\textbf{u}\Vert_{[H^2(\Omega)]^2}^{1-1/p}\leq \delta\Vert\textbf{u}\Vert_{[H^2(\Omega)]^2}^{2-2/p}+C(\delta)\Vert D\textbf{u}\Vert^{2+2/p}
,
\end{align*}
for some $\delta>0$.
Moreover, recalling \eqref{lambdaest}, \eqref{lp}, \eqref{energy3}, \eqref{psi3} and Sobolev embeddings, we get
\begin{align*}
&\Vert(I+\alpha\textbf{A})^{-2}\textbf{P}\text{div}(\lambda_\alpha(\phi)\nabla\psi^T\nabla\psi)\Vert_{[L^p(\Omega)]^2}\\&\leq C\Vert\nabla\phi\Vert_{[L^p(\Omega)]^2}\Vert\nabla\psi\Vert_{[L^{\infty}(\Omega)]^4}^2+C\Vert\nabla\psi\Vert_{[L^\infty(\Omega)]^4}\Vert\psi\Vert_{[W^{2,p}(\Omega)]^2}\\&\leq C(1+\Vert\nabla\phi\Vert_{[L^p(\Omega)]^2})\Vert\psi\Vert_{[W^{2,p}(\Omega)]^2}^2.
\end{align*}
Also, by Sobolev embeddings, observe that
$$
\left\Vert\frac{\nu(\phi)}{2k(\phi)}(1-\phi)\textbf{u}\right\Vert_{[L^p(\Omega)]^2}\leq C\Vert\textbf{u}\Vert_{[L^p(\Omega)]^2}\leq C\Vert D\textbf{u}\Vert.
$$
Therefore, by \eqref{lpphi2} and the embedding $W^{2,p}(\Omega)\hookrightarrow W^{1,\infty}(\Omega)$, we deduce
\begin{align*}
&\Vert\nabla\varphi\Vert_{[L^\infty(\Omega)]^2}\Vert D\textbf{u}\Vert_{[L^p(\Omega)]^4}\leq C\Vert\varphi\Vert_{W^{2,p}(\Omega)}\Vert D\textbf{u}\Vert^{2/p}\Vert\textbf{u}\Vert_{[H^2(\Omega)]^2}^{1-2/p}\\&\leq C\Vert\phi\Vert_{W^{2,p}(\Omega)}\Vert D\textbf{u}\Vert^{2/p}\Vert\textbf{u}\Vert_{ [H^2(\Omega)]^2}^{1-2/p}\leq \delta\Vert\textbf{u}\Vert_{[H^2(\Omega)]^{2}}^{2-4/p}+C(\delta)\Vert\phi\Vert_{W^{2,p}(\Omega)}^{2}\Vert D\textbf{u}\Vert^{{4}/{p}}.
\end{align*}
Adding now the results together and exploiting \eqref{interp3}, we find
\begin{align}
\nonumber&\Vert\textbf{u}\Vert_{[W^{2,p}(\Omega)]^2}\\
&\nonumber\leq \delta\left(\Vert\textbf{u}\Vert_{[H^2(\Omega)]^2}^{2-2/p}+\Vert\textbf{u}\Vert_{[H^2(\Omega)]^2}^{2-4/p}\right)+C(\delta)\left( \Vert\partial_t\textbf{u}\Vert_{\textbf{V}_\sigma^{1/2-\varepsilon}}+ \Vert\phi\Vert_{W^{2,p}(\Omega)}^2+\Vert D\textbf{u}\Vert^{2+2/p}\right.\\&\left.+(1+\Vert\nabla\phi\Vert_{L^p(\Omega)})\Vert\psi\Vert_{[W^{2,p}(\Omega)]^2}^2+\Vert D\textbf{u}\Vert+\Vert\phi\Vert_{W^{2,p}(\Omega)}^{2}\Vert D\textbf{u}\Vert^{{4}/{p}}\right).
\label{W2p}
\end{align}
Combining \eqref{W2p} with \eqref{psiw2} and using the embeddings $W^{2,p}(\Omega)\hookrightarrow W^{1,\infty}(\Omega)$, $W^{2,p}(\Omega)\hookrightarrow H^2(\Omega)\hookrightarrow W^{1,p}(\Omega)$, recalling the range of $p$,
we obtain
\begin{align}
&\nonumber\frac{1}{2}\frac{d}{dt}\Vert \psi\Vert_{[W^{2,p}(\Omega)]^2}^2 \\&\leq\nonumber C\left(\delta\left(\Vert\textbf{u}\Vert_{[H^2(\Omega)]^2}^{2-2/p}+\Vert\textbf{u}\Vert_{[H^2(\Omega)]^2}^{2-4/p}\right)+C(\delta)\left( \Vert\partial_t\textbf{u}\Vert_{\textbf{V}_\sigma^{1/2-\varepsilon}}+ \Vert\phi\Vert_{W^{2,p}(\Omega)}^2+\Vert D\textbf{u}\Vert^{2+2/p}\right.\right.\\&\left.\left.+(1+\Vert\nabla\phi\Vert_{L^p(\Omega)})\Vert\psi\Vert_{[W^{2,p}(\Omega)]^2}^2+\Vert D\textbf{u}\Vert+\Vert\phi\Vert_{W^{2,p}(\Omega)}^{2}\Vert D\textbf{u}\Vert^{{4}/{p}}\right)\right)\Vert\psi\Vert_{[W^{2,p}(\Omega)]^2}^2\nonumber\\&\nonumber\leq
C\left(1+\Vert\phi\Vert_{W^{2,p}(\Omega)}^4+\left(\delta^{\frac{p}{p-1}}+\delta^{\frac{p}{p-2}}\right)\Vert\textbf{u}\Vert_{H^2(\Omega)}^{2}+\Vert D\textbf{u}\Vert^{4+4/p}+\Vert D\textbf{u}\Vert^2+\Vert D\textbf{u}\Vert^{8/p}\right.\\&\left.+\Vert\phi\Vert_{W^{2,p}(\Omega)}^2+\Vert\psi\Vert_{[W^{2,p}(\Omega)]^2}^{8}+\Vert\psi\Vert_{[W^{2,p}(\Omega)]^2}^4+\Vert\psi\Vert_{[W^{2,p}(\Omega)]^2}^{2p}+\Vert\psi\Vert_{[W^{2,p}(\Omega)]^2}^{p}\right)\nonumber+\delta\Vert\partial_t\textbf{u}\Vert_{\textbf{V}_\sigma^{1/2-\varepsilon}}^2\nonumber\\&\leq C\left(1+\Vert\phi\Vert_{W^{2,p}(\Omega)}^4+\left(\delta^{\frac{p}{p-1}}+\delta^{\frac{p}{p-2}}\right)\Vert\textbf{u}\Vert_{H^2(\Omega)}^{2}+\Vert D\textbf{u}\Vert^{4+4/p}+\Vert\psi\Vert_{[W^{2,p}(\Omega)]^2}^{8}\right)+\delta\Vert\partial_t\textbf{u}\Vert_{\textbf{V}_\sigma^{1/2-\varepsilon}}^2.
\label{phibis}
\end{align}
\textbf{Final estimate.} We can now conclude by collecting  \eqref{estvelox}, \eqref{estmu}, \eqref{estu}, \eqref{phibis} and adding them together. Choosing $\delta$ sufficiently small, we find
\begin{align*}
&\frac{d}{dt}\mathcal{Y}+\frac{1}{4}\Vert\nabla\partial_t\phi\Vert^2+\frac{1}{4}\Vert \partial_t\textbf{u}\Vert^2+ \frac{1}{16\hat{C}}\Vert\textbf{u}\Vert^2_{[H^2(\Omega)]^2}+\Vert \partial_t\textbf{u}\Vert_{\textbf{V}^{1/2-\varepsilon}_\sigma}^2\left(\frac{\nu_1+\nu_2}{8}-\frac{C_3\vert \nu_1-\nu_2\vert }{4}\right)\\&\leq C\left(1+\mathcal{Y}\right)^{\frac{9}{1-6\varepsilon}},
\end{align*}
where (see \eqref{vcontrol})
$$
\mathcal{Y}:= C+\frac{1}{2}\Vert D\textbf{u}\Vert^2+\frac{1}{2}\Vert\psi\Vert^2_{[W^{2,p}(\Omega)]^2}+\frac{1}{2}\Vert\nabla\mu\Vert^2+(\textbf{v}\cdot \nabla\phi,\mu)+\frac{1}{2}\Vert \partial_t\textbf{u}\Vert^2_{\textbf{V}_\sigma^{-1/2-\varepsilon}}\geq 0.
$$
Indeed, we have exploited \eqref{w2p2} and the right-hand side term with the highest degree is $\Vert\phi\Vert_{H^2(\Omega)}^{\frac{8}{1-6\varepsilon}}\Vert \partial_t\textbf{u}\Vert^2_{\textbf{V}_\sigma^{-1/2-\varepsilon}}$, giving rise to a term like $\Vert\psi\Vert_{[W^{2,p}(\Omega)]^2}^{\frac{16}{1-6\varepsilon}}\Vert \partial_t\textbf{u}\Vert^2_{\textbf{V}_\sigma^{-1/2-\varepsilon}}$, which is controlled, for example, by $C\left(1+\mathcal{Y}\right)^{\frac{9}{1-6\varepsilon}}$.
Notice that the left-hand side is positive as long as we suppose that
\begin{align}
\frac{\vert \nu_1-\nu_2\vert }{\nu_1+\nu_2}<\frac{1}{2C_3},
\end{align}
which is exactly the smallness condition \eqref{small}.
We now need to find an estimate for $\mathcal{Y}(0)$. First, by \eqref{regularities}, \eqref{muAC} and \eqref{ureg}, we can evaluate $\mathcal{Y}$ at $t=0$, obtaining, recalling that $\textbf{v}(0)=(I+\alpha\textbf{A})^{-2}\textbf{u}_{0,\alpha}$ and by \eqref{lp},
\begin{align}
\nonumber\mathcal{Y}(0)&\leq \tilde{C}+\frac{1}{2}\Vert \textbf{u}_{0,\alpha}\Vert^2_{\textbf{V}_\sigma}+\frac{1}{2}\Vert\psi_{0,\alpha}\Vert^2_{[W^{2,p}(\Omega)]^2}+\frac{1}{2}\Vert\mu_0\Vert^2_V\\&+\Vert \textbf{u}_{0,\alpha}\Vert_{[L^4(\Omega)]^2}\Vert \nabla\phi_0\Vert_{[L^4(\Omega)]^2} \Vert \mu_0\Vert+\frac{1}{2}\Vert \partial_t\textbf{u}(0)\Vert^2_{\textbf{V}_\sigma^{-1/2-\varepsilon}}.
\label{y0}
\end{align}
Due to assumptions \ref{Hbegin}-\ref{Hend}, and being $\Vert\textbf{u}_{0,\alpha}\Vert_{\textbf{V}_\sigma}\leq 1+\Vert\textbf{u}_{0}\Vert_{\textbf{V}_\sigma}$ and $\Vert\psi_{0,\alpha}\Vert_{[W^{2,p}(\Omega)]^2}\leq 1+\Vert\psi_{0}\Vert_{[W^{2,p}(\Omega)]^2}$, independently of $\alpha$, we are only left to find a uniform bound for $\Vert \partial_t\textbf{u}(0)\Vert^2_{\textbf{V}_\sigma^{-1/2-\varepsilon}}$.
On account of \eqref{ureg}, we can multiply \eqref{problem}$_1$ by $\textbf{A}^{-1/2-\varepsilon}\partial_t\textbf{u}(0)$, integrate over $\Omega$ and then evaluate it at $t=0$. This gives (recall that $\varphi_0=(I+\alpha A_1 )^{-1}\phi_0$)
\begin{align}
&\nonumber\Vert\partial_t\textbf{u}(0)\Vert_{\textbf{V}_\sigma^{-1/2-\varepsilon}}^2+((\textbf{v}(0)\cdot\nabla)\textbf{u}_{0,\alpha},\textbf{A}^{-1/2-\varepsilon}\partial_t\textbf{u}(0))-(\text{div}(\nu(\varphi_0)D\textbf{u}_{0,\alpha}),\textbf{A}^{-1/2-\varepsilon}\partial_t\textbf{u}(0))\\&\nonumber+((I+\alpha \textbf{A})^{-2}\text{div}(\lambda_\alpha(\phi_0)\nabla\psi_{0,\alpha}^T\nabla\psi),\textbf{A}^{-1/2-\varepsilon}\partial_t\textbf{u}(0) )\\
&\nonumber +((I+\alpha \textbf{A})^{-2}\text{div}(\nabla\phi_0\otimes\nabla\phi_0),\textbf{A}^{-1/2-\varepsilon}\partial_t\textbf{u}(0))\\
&+(\nu(\phi_0)\frac{(1-\phi_0)\textbf{u}_{0,\alpha}}{2k(\phi_0)},\textbf{A}^{-1/2-\varepsilon}\partial_t\textbf{u}(0))=0.
\label{zero}
\end{align}
We now apply a similar argument as in the fourth estimate. In particular we write
\begin{align*}
&-(\text{div}(\nu(\varphi_0)D\textbf{u}_{0,\alpha}),\textbf{A}^{-1/2-\varepsilon}\partial_t\textbf{u}(0))\\&=-(\nu^\prime(\varphi_0)D\textbf{u}_{0,\alpha}\nabla\varphi_0,\textbf{A}^{-1/2-\varepsilon}\partial_t\textbf{u}(0)\nabla\varphi_0)-{\frac{1}{2}(\Delta\textbf{u}_{0,\alpha},\nu(\varphi_0)\textbf{A}^{-1/2-\varepsilon}\partial_t\textbf{u}(0))}\\&=-(\nu^\prime(\varphi_0)D\textbf{u}_{0,\alpha}\nabla\varphi_0,\textbf{A}^{-1/2-\varepsilon}\partial_t\textbf{u}(0)\nabla\varphi_0)\\&+{\frac{1}{2}(\textbf{A}\textbf{u}_{0,\alpha},\nu(\varphi_0)\textbf{A}^{-1/2-\varepsilon}\partial_t\textbf{u}(0))}+\frac{1}{2}(\nu^\prime(\varphi_0)\nabla\varphi_0\tilde{\pi},\textbf{A}^{-1/2-\varepsilon}\partial_t\textbf{u}(0))
,
\end{align*}
where $\tilde{\pi}$ satisfies $-\Delta\textbf{u}_{0,\alpha}+\nabla\tilde{\pi}= \textbf{A}\textbf{u}_{0,\alpha}$.
By Lemma \ref{press},  we have, for $\widetilde{\delta}\in(0,1/4]$ (see \ref{Hbegin}),
\begin{align}
\Vert \tilde{\pi}\Vert\leq C\left(\Vert\textbf{u}_{0,\alpha}\Vert_{\textbf{V}_\sigma}+\Vert\textbf{u}_{0,\alpha}\Vert_{\textbf{V}_\sigma^{3/2+2\widetilde{\delta}}}\right).
\label{p2}
\end{align}
Moreover, again by \cite[Lemma 4.11]{Maxwell}), recalling that $\nu(\cdot)$ is linear on $[-1,1]$ and using \eqref{lpphi2}, we obtain
\begin{align}
&\nonumber\Vert \nu(\varphi_0)\textbf{A}^{-1/2-\varepsilon}\partial_t\textbf{u}(0)\Vert_{[H^{1/2-2\widetilde{\delta}}(\Omega)]^2}\\&\nonumber\leq {C}\left(\Vert \nu(\varphi_0)\Vert_{L^\infty(\Omega)}\Vert\textbf{A}^{-1/2-\varepsilon}\partial_t\textbf{u}(0)\Vert_{[H^{1/2-2\widetilde{\delta}}(\Omega)]^2}+\Vert \nu(\varphi_0)\Vert_{H^2(\Omega)}\Vert\textbf{A}^{-1/2-\varepsilon}\partial_t\textbf{u}(0)\Vert\right)\\&\leq {C}\left(C\Vert\textbf{A}^{-1/2-\varepsilon}\partial_t\textbf{u}(0)\Vert_{[H^{1/2-2\widetilde{\delta}}(\Omega)]^2}+C(1+\Vert \phi_0\Vert_{H^2(\Omega)})\Vert\textbf{A}^{-1/2-\varepsilon}\partial_t\textbf{u}(0)\Vert\right)\nonumber\\&\leq {C}\left(\Vert\textbf{A}^{-1/4-\varepsilon-\widetilde{\delta}}\partial_t\textbf{u}(0)\Vert+C(1+\Vert \phi_0\Vert_{H^2(\Omega)})\Vert\textbf{A}^{-1/2-\varepsilon}\partial_t\textbf{u}(0)\Vert\right)\nonumber\\&\leq C(1+\Vert \phi_0\Vert_{H^2(\Omega)})\Vert\partial_t\textbf{u}(0)\Vert_{\textbf{V}_\sigma^{-1/2-\varepsilon}}
.
\label{Maxwell2}
\end{align}
Therefore, by \eqref{lpphi2}, \eqref{interp4bis3} and recalling that $\nu^\prime\in W^{1,\infty}(\R)$, we deduce
\begin{align*}
&\vert (\nu^\prime(\varphi_0)D\textbf{u}_{0,\alpha}\nabla\varphi_0,\textbf{A}^{-1/2-\varepsilon}\partial_t\textbf{u}(0)\nabla\varphi_0)\vert \leq C\Vert D\textbf{u}_{0,\alpha}\Vert\Vert\textbf{A}^{-1/2-\varepsilon}\partial_t\textbf{u}(0)\Vert\Vert\nabla\varphi_0\Vert_{[L^\infty(\Omega)]^2}\\&\leq C\Vert\nabla \textbf{u}_{0,\alpha}\Vert\Vert\phi_0\Vert_{W^{2,p}(\Omega)}\Vert\partial_t\textbf{u}(0)\Vert_{\textbf{V}_\sigma^{-1/2-\varepsilon}},
\end{align*}
and (see \eqref{p2})
\begin{align*}
&\left\vert\frac{1}{2}(\nu^\prime(\varphi_0)\nabla\varphi_0\tilde{\pi},\textbf{A}^{-1/2-\varepsilon}\partial_t\textbf{u}(0))\right\vert\leq C\Vert\textbf{A}^{-1/2-\varepsilon}\partial_t\textbf{u}(0)\Vert\Vert\tilde{\pi}\Vert\Vert\nabla\varphi_0\Vert_{[L^\infty(\Omega)]^2}\\&\leq C\Vert\partial_t\textbf{u}(0)\Vert_{\textbf{V}_\sigma^{-1/2-\varepsilon}}\left(\Vert\textbf{u}_{0,\alpha}\Vert_{\textbf{V}_\sigma}+\Vert\textbf{u}_{0,\alpha}\Vert_{\textbf{V}_\sigma^{3/2+2\widetilde{\delta}}}\right)\Vert\phi_0\Vert_{W^{2,p}(\Omega)}.
\end{align*}
On the other hand, by \eqref{lpphi}, \eqref{lpphi2} and \eqref{Maxwell2}, we have
\begin{align*}
&\left\vert{\frac{1}{2}(\textbf{A}\textbf{u}_{0,\alpha},\nu(\varphi_0)\textbf{A}^{-1/2-\varepsilon}\partial_t\textbf{u}(0))}\right\vert\leq C\Vert\textbf{A}^{1/4-\widetilde{\delta}}\textbf{P}(\nu(\varphi_0)\textbf{A}^{-1/2-\varepsilon}\partial_t\textbf{u}(0))\Vert\Vert\textbf{A}^{3/4+\widetilde{\delta}}\textbf{u}_{0,\alpha}\Vert\\&\leq C\Vert\nu(\varphi_0)\textbf{A}^{-1/2-\varepsilon}\partial_t\textbf{u}(0)\Vert_{H^{1/2-2\widetilde{\delta}}}\Vert\textbf{A}^{3/4+\widetilde{\delta}}\textbf{u}_{0,\alpha}\Vert\\&\leq C(1+\Vert \phi_0\Vert_{H^2(\Omega)})\Vert\partial_t\textbf{u}(0)\Vert_{\textbf{V}_\sigma^{-1/2-\varepsilon}}\Vert\textbf{u}_{0,\alpha}\Vert_{\textbf{V}_\sigma^{3/2+2\widetilde{\delta}}}.
\end{align*}
Observe again that $\widetilde{\delta}>0$ is needed to ensure that $\textbf{P}(\nu(\varphi_0)\textbf{A}^{-1/2-\varepsilon}\partial_t\textbf{u}(0))\in D(\textbf{A}^{1/4-\widetilde{\delta}})$.

The above bounds combined with applying standard inequalities and  \eqref{w2p}, \eqref{lp} allow us to infer from \eqref{zero} that
\begin{align}
&\nonumber\Vert\partial_t\textbf{u}(0)\Vert_{\textbf{V}_\sigma^{-1/2-\varepsilon}}^2\leq
C\Vert\textbf{u}_{0,\alpha}\Vert_{[L^4(\Omega)]^2}^2\Vert\textbf{u}_{0,\alpha}\Vert_{[W^{1,4}(\Omega)]^2}^2+C\Vert\nabla\psi_{0,\alpha}\Vert_{[L^4(\Omega)]^4}^4\Vert\nabla\phi_0\Vert_{[L^\infty(\Omega)]^2}^2\\&+C\Vert\nabla\psi_{0,\alpha}\Vert_{[L^\infty(\Omega)]^2}\Vert\psi_{0,\alpha}\Vert_{[H^2(\Omega)]^2}^2+C\Vert\nabla\phi_{0}\Vert_{[L^\infty(\Omega)]^2}\Vert\phi_{0}\Vert_{[H^2(\Omega)]^2}^2+C\Vert\textbf{u}_{0,\alpha}\Vert^2\nonumber\\&+C\Vert\nabla \textbf{u}_{0,\alpha}\Vert^2\Vert\phi_0\Vert_{W^{2,p}(\Omega)}^2+C\left(\Vert\textbf{u}_{0,\alpha}\Vert_{\textbf{V}_\sigma}^2+\Vert\textbf{u}_{0,\alpha}\Vert_{\textbf{V}_\sigma^{3/2+2\widetilde{\delta}}}^2\right)\Vert\phi_0\Vert_{W^{2,p}(\Omega)}^2\nonumber\\&+C(1+\Vert \phi_0\Vert_{H^2(\Omega)}^2)\Vert\textbf{u}_{0,\alpha}\Vert_{\textbf{V}_\sigma^{3/2+2\widetilde{\delta}}}^2\leq C,
\end{align}
independently of $\alpha$, since $\Vert\textbf{u}_{0,\alpha}\Vert_{\textbf{V}_\sigma^{3/2+2\widetilde{\delta}}}\leq C$, $\Vert\psi_{0,\alpha}\Vert_{[W^{2,p}(\Omega)]^2}\leq C$ and $\phi_0\in H^3(\Omega)$ (see Remark \ref{separated}). Therefore we can conclude by \eqref{y0} that $\mathcal{Y}(0)$ (which a priori depends on $\alpha$) is bounded uniformly with respect to $\alpha$. Thus, by \cite[Lemma II.4.12]{Boyerlibro}, there exists a maximal time $T_{M,\alpha}>0$, decreasing as $\mathcal{Y}(0)$ increases, such that
\begin{equation}
\mathcal{Y}(t)\leq C.
\label{last3}
\end{equation}
for any $t\in[0,T]$, for any $T<T_{M,\alpha}$, and also
\begin{align}
\Vert\partial_t\phi\Vert_{L^2(0,T;V)}+\Vert\partial_t\textbf{u}\Vert_{L^2(0,T;\textbf{V}_\sigma^{1/2-\varepsilon})}+\Vert\textbf{u}\Vert_{L^2(0,T;\textbf{W}_\sigma)}\leq C,
\label{L22}
\end{align}
for some $C$ independent of $\alpha$. Clearly, the constants involved as well as $T_{M,\alpha}$ may also depend on $\varepsilon$. Being $\mathcal{Y}(0)$ uniformly bounded by a constant independent of $\alpha$, we see that actually there exists a minimum (strictly positive) time $\overline{T}_M=\overline{T}_M(\varepsilon)>0$ such that, independently of $\alpha$, all the solutions are defined over $[0,T_M]$, for any $T_M< \overline{T}_M$. Therefore, \eqref{last3} and \eqref{L22} hold for any $T_M<\overline{T}_M$, independently of $\alpha$.
These bounds allow to conclude, by a standard argument, that the maximal time $T_0$ of the approximated solution actually coincides with $T_M$.
Moreover, for any $T_M< \overline{T}_M$, we have obtained the following regularity for the approximating solution
\begin{align}
\begin{cases}
\textbf{u}_{\alpha}\in L^\infty(0,T_M;\textbf{V}_\sigma)\cap L^2(0,T_M;\textbf{W}_\sigma)\cap H^1(0,T_M;\textbf{V}_\sigma^{1/2-\varepsilon})\cap W^{1,\infty}(0,T_M;\textbf{V}_\sigma^{-1/2-\varepsilon}),\\
\phi_{\alpha}\in L^\infty(0,T_M;V)\cap L^2(0,T_M;H^2(\Omega))\cap H^1(0,T_M;V),\\
{\mu}_{\alpha}\in L^\infty(0,T_M;V),\\
{\psi}_{\alpha}\in L^\infty(0,T_M; [W^{2,p}(\Omega)]^2),
\end{cases}
\label{regularities2}
\end{align}
for any $p\in \left[2,\frac{4}{1+2\varepsilon}\right]$.
Moreover, from \eqref{psit3} we have
$$
\psi_{\alpha}\in W^{1,\infty}(0,T_M;[H^1(\Omega)]^2)
$$
and by \eqref{pot} we get
$$
\phi_{\alpha}\in L^\infty(0,T_M; W^{2,p}(\Omega)),\quad F^\prime(\phi_{\alpha})\in L^\infty(0,T_M;L^p(\Omega)),
$$
for any $p\in \left[2,\frac{4}{1+2\varepsilon}\right]$.
\subsection{Letting $\alpha$ go to $0$}
\label{limits}
Using weak and weak$^*$ compactness,  we can extract a subsequence
$$
\{(\textbf{u}_{\alpha_k}, \phi_{\alpha_k}, \psi_{\alpha_k}) \}_{k\in \mathbb{N}}
$$
which suitably converges as $k \to \infty$ to a triplet $(\textbf{u}, \phi, \psi)$ such that
\begin{align}
\begin{cases}
\textbf{u}\in L^\infty(0,T_M;\textbf{V}_\sigma)\cap L^2(0,T_M;\textbf{W}_\sigma)\cap H^1(0,T_M;\textbf{V}_\sigma^{1/2-\varepsilon})\cap W^{1,\infty}(0,T_M;\textbf{V}_\sigma^{-1/2-\varepsilon}),\\
\phi\in L^\infty(0,T_M;W^{2,p}(\Omega))\cap  H^1(0,T_M;V),\\
{\mu}\in L^\infty(0,T_M;V),\\
{\psi}\in L^\infty(0,T_M; [W^{2,p}(\Omega)]^2)\cap W^{1,\infty}(0,T_M;[H^1(\Omega)]^2),
\end{cases}
\label{regularities3}
\end{align}
for any $p\in\left[2,\frac{4}{1+2\varepsilon}\right]$.
In particular,  by the embeddings $\textbf{W}_\sigma\hookrightarrow\hookrightarrow\textbf{V}_\sigma\hookrightarrow\textbf{H}_\sigma$ and $W^{2,q}(\Omega)\hookrightarrow\hookrightarrow W^{1,q}(\Omega)\hookrightarrow V$ (for any $q\in[2,\infty)$) and the Aubin-Lions Lemma, we get the following strong convergences (always up to subsequences)
\begin{align}
\begin{cases}
\textbf{u}_{\alpha_k}\rightarrow\textbf{u}\quad \text{ in }L^2(0,T_M;\textbf{V}_\sigma),\\
\phi_{\alpha_k}\rightarrow \phi \quad \text{ in }C([0,T_M];W^{1,4}(\Omega)),\\
\psi_{\alpha_k}\rightarrow \psi \quad\text{ in }C([0,T_M];W^{1,4}(\Omega)).
\end{cases}
\label{convergences}
\end{align}
Observe now that
$$
\textbf{w}_{\alpha_k}+\alpha\textbf{A}\textbf{w}_{\alpha_k}=\textbf{u}_{\alpha_k},\qquad
\textbf{v}_{\alpha_k}+\alpha\textbf{A}\textbf{v}_{\alpha_k}=\textbf{w}_{\alpha_k}.
$$
From these equalities we immediately deduce (see \eqref{w2p})
$$
\Vert\textbf{v}_{\alpha_k}-\textbf{u}_{\alpha_k}\Vert\leq \Vert\textbf{v}_{\alpha_k}-\textbf{w}_{\alpha_k}\Vert+\Vert\textbf{u}_{\alpha_k}-\textbf{w}_{\alpha_k}\Vert\leq \alpha\left(\Vert\textbf{A}\textbf{v}_{\alpha_k}\Vert+\Vert\textbf{A}\textbf{w}_{\alpha_k}\Vert\right)\leq \alpha C\Vert\textbf{u}_{\alpha_k}\Vert_{\textbf{W}_\sigma}
$$
and, being $\Vert\textbf{u}_{\alpha_k}\Vert_{L^2(0,T_M;\textbf{W}_\sigma)}\leq C$, from \eqref{convergences} we get, as $k\rightarrow\infty$,
$$
\textbf{v}_{\alpha_k}\rightarrow\textbf{u}\quad \text{ in }L^2(0,T_M;\textbf{H}_\sigma).
$$
Now observe that, due to \eqref{w2p}, the same uniform estimates about $\textbf{u}_{\alpha_k}$ also hold for $\textbf{v}_{\alpha_k}$  (see \eqref{regularities2}).
Thus by the previous convergence, we deduce
\begin{align}
\textbf{v}_{\alpha_k}\rightarrow\textbf{u}\quad \text{ in }L^2(0,T_M;\textbf{V}_\sigma).
\label{vconvergence}
\end{align}
By a similar argument, being $\phi_{\alpha_k}=(I+\alpha_k A_1)^{-1}\varphi_{\alpha_k}$, thanks to \eqref{lpphi2} we infer
$$
\Vert\varphi_{\alpha_k}-\phi_{\alpha_k}\Vert_{L^\infty(0,T_M;H)}\leq C\alpha_k\Vert\phi_{\alpha_k}\Vert_{L^\infty(0,T_M;H)}\leq C\alpha_k
$$
and
 \begin{align}
\varphi_{\alpha_k}\rightarrow\phi\quad \text{ in }L^\infty(0,T_M;H),\quad \nu(\varphi_{\alpha_k})\rightarrow\nu(\phi)\quad\text{ in }L^\infty(0,T_M;H).
\label{phiconvergence}
\end{align}
The above convergences are enough to prove that $(\textbf{u},\ \phi,\ \psi)$ is a solution to \eqref{syst} starting from \eqref{problem}.
We only show the convergence of some nonlinear term, which could be somehow non-standard.
In particular, from \eqref{regularities2} and \eqref{convergences} we deduce
\begin{align}
\text{div}(\nabla\phi_{\alpha_k}\otimes\nabla\phi_{\alpha_k})\rightharpoonup\text{div}\left(\nabla\phi\otimes\nabla\phi\right)\quad\text{ in }L^2(0,T;\textbf{H}_\sigma).
\label{div}
\end{align}
Now, by a similar argument as the one to obtain \eqref{vconvergence}, given $\textbf{w}_{\alpha_{k}}=(I+{\alpha_{k}}\textbf{A})^{-1}\textbf{v}$, we have
$$
\Vert\textbf{v}-\textbf{w}_{\alpha_{k}}\Vert\leq {\alpha_{k}}\Vert\nabla\textbf{v}\Vert\Vert\nabla(\textbf{w}_{\alpha_{k}}-\textbf{v})\Vert,
$$
and if $\textbf{v}\in L^2(0,T_M;\textbf{V}_\sigma)$ we immediately obtain by \eqref{lp1} and \eqref{h1} that
$$
\Vert\textbf{v}-\textbf{w}_{\alpha_{k}}\Vert_{L^2(0,T_M;\textbf{H}_\sigma)}^2\leq {\alpha_{k}} C,
$$
entailing $\textbf{w}_{\alpha_{k}}\rightarrow \textbf{v}$ as $k\rightarrow\infty$. Therefore, for $\textbf{w}\in L^2(0,T_M;\textbf{V}_\sigma)$, being $\Vert (I+{\alpha_{k}}\textbf{A})^{-1}\Vert_{\mathcal{L}(\textbf{H}_\sigma)}\leq C$ (see \eqref{lp1}), we get
\begin{align*}
&\Vert(I+{\alpha_{k}}\textbf{A})^{-2}\textbf{w}-\textbf{w}\Vert\\
&\leq  \Vert(I+{\alpha_{k}}\textbf{A})^{-2}\textbf{w}-(I+{\alpha_{k}}\textbf{A})^{-1}\textbf{w}\Vert+\Vert(I+{\alpha_{k}}\textbf{A})^{-1}\textbf{w}-\textbf{w}\Vert\\&\leq (C+1)\Vert(I+{\alpha_{k}}\textbf{A})^{-1}\textbf{w}-\textbf{w}\Vert.
\end{align*}
Hence we infer that
\begin{align}
(I+{\alpha_{k}}\textbf{A})^{-2}\textbf{w}\rightarrow\textbf{w}\quad\text{in }L^2(0,T_M;\textbf{H}_\sigma).
\label{conv}
\end{align}
Let then $\textbf{w}\in L^2(0,T_M;\textbf{V}_\sigma)$. Observe that
\begin{align*}
&\left\vert\int_0^{T_M}\left((I+{\alpha_{k}}\textbf{A})^{-2}\textbf{P}\text{div}(\nabla\phi_{\alpha_k}\otimes\nabla\phi_{\alpha_k}),\textbf{w}\right)ds-\int_0^{T_M}\left(\text{div}(\nabla\phi\otimes\nabla\phi),\textbf{w}\right)ds\right\vert\\&\leq \left\vert\int_0^{T_M}\left(\text{div}(\nabla\phi_{\alpha_k}\otimes\nabla\phi_{\alpha_k})-\text{div}(\nabla\phi\otimes\nabla\phi),\textbf{w}\right)ds\right\vert\\&+\left\vert\int_0^{T_M}\left(\text{div}(\nabla\phi_{\alpha_k}\otimes\nabla\phi_{\alpha_k}),(I+{\alpha_{k}}\textbf{A})^{-2}\textbf{w}-\textbf{w}\right)ds\right\vert,
\end{align*}
and the two terms on the right-hand side go to zero thanks to \eqref{regularities2}, \eqref{div} and \eqref{conv}. Then, by a density argument, we can extend the result for $\textbf{w}\in L^2(0,T_M;\textbf{H}_\sigma)$, since, given $\{\textbf{w}_n\}_n\in L^2(0,T_M;\textbf{V}_\sigma)$ converging to $\textbf{w}$ in  $L^2(0,T_M;\textbf{H}_\sigma)$, we get, as $n\to\infty$, thanks to \eqref{regularities3},
$$
\int_0^{T_M}\left(\text{div}(\nabla\phi\otimes\nabla\phi),\textbf{w}_n\right)ds\rightarrow\int_0^{T_M}\left(\text{div}(\nabla\phi\otimes\nabla\phi),\textbf{w}\right)ds.
$$
A very similar proof holds for the term $(I+{\alpha_{k}} \textbf{A})^{-2}\textbf{P}\text{div}(\lambda_{\alpha_k}(\phi)\nabla\psi^T\nabla\psi)$, exploiting the convergence properties of $\lambda_{\alpha_{k}},\lambda_{\alpha_{k}}^\prime$ uniformly on compact sets of $\R$.
From these results we can pass to the limit as $k\rightarrow\infty$ in \eqref{problem} and the existence part is proven.
Note that, using \cite[Lemma A.6,][]{contigiorgini}, being $\mu\in L^\infty(0,T_M;V)$ and $\psi\in L^\infty(0,T_M;[W^{2,3}(\Omega)]^2)$, we obtain $\Psi^{\prime\prime}(\phi)\in L^\infty(0,T_M;L^q(\Omega))$, for any $q\in [2,\infty)$. Moreover, the property $\vert\phi(x,t)\vert\leq 1$ for almost any $(x,t)\in \Omega\times(0,T_M)$, is obtained in a standard way, being already $\vert\phi_{\alpha_k}(x,t)\vert\leq 1$ for any $k$.

Concerning the strict separation property: fix $q>2$. Being $F'\in C^2((-1,1))$, we can apply the chain rule to obtain
	\begin{equation*}
	\nabla F'(\phi)=F''(\phi)\nabla\phi,
	\end{equation*}
	but then $\Vert \nabla F'(\phi)\Vert_{L^q(\Omega)} \leq \Vert F''(\phi) \Vert_{L^{2q}(\Omega)}\Vert \nabla\phi\Vert_{L^{2q}(\Omega)}\leq C$ for almost any $t\in(0,T_M)$, since it holds $\phi\in L^\infty(0,T_M;W^{1,2q}(\Omega))$ and $\Psi^{\prime\prime}(\phi)\in L^\infty(0,T_M;L^q(\Omega))$. Therefore we get $F'(\phi)\in L^\infty(0,T_M;W^{1,q}(\Omega))$, implying that $F'(\phi)\in L^\infty(\Omega\times(0,T_M))$ thanks to the embedding $W^{1,q}(\Omega)\hookrightarrow L^\infty(\Omega)$. Thus $\phi$ is strictly separated, namely, there exists $\hat{\delta}>0$ such that
	\begin{equation*}
	\max_{t\in[0,T_M]}\Vert \phi\Vert_{C(\overline{\Omega})}\leq 1-\hat{\delta}.
	\end{equation*}	
	Notice that this control has been possible since $\phi \in C([0,T_M]; C(\overline{\Omega}))$ (see Remark \ref{rmm}).

In conclusion, arguing as in the approximated case to get \eqref{muAC}, we obtain $\mu\in L^2(0,T_M;H^3(\Omega))$ by elliptic regularity.
\subsection{Continuous dependence estimate and uniqueness}
\label{uni}
To prove uniqueness, we also suppose that $\lambda\in C^{1,1}(\R)$.
Let us set
$$\textbf{u}=\textbf{u}_1-\textbf{u}_2,\quad p=p_1-p_2, \quad \phi=\phi_1-\phi_2, \quad \psi=\psi_1-\psi_2,$$
and
$$
\mu=-\Delta\phi+\Psi'(\phi_1)-\Psi'(\phi_2)+\frac{\lambda^{\prime}(\phi_1)}{2}\vert\nabla\psi_1\vert^2-\frac{\lambda^{\prime}(\phi_2)}{2}\vert\nabla\psi_2\vert^2.
$$
Then we can write
\begin{align}
&\nonumber \partial_t\textbf{u}+({\textbf{u}}_1\cdot\nabla)\textbf{u}+({\textbf{u}}\cdot\nabla)\textbf{u}_2+\nabla p-\text{div}(\nu(\phi_1) D\textbf{u})-\text{div}((\nu(\phi_1)-\nu(\phi_1)) D\textbf{u}_2)\\&\nonumber+\text{div}(\lambda({\phi}_1)\nabla{\psi}_1^T\nabla{\psi})+\text{div}((\lambda({\phi}_1)-\lambda({\phi}_2))\nabla{\psi}_1^T\nabla{\psi}_2)\\&\nonumber+\text{div}(\lambda({\phi}_2)\nabla{\psi}^T\nabla{\psi}_2)+\text{div}(\nabla{\phi}_1\otimes\nabla{\phi})\\&+\text{div}(\nabla{\phi}\otimes\nabla{\phi}_2)-\nu({\phi }_2)\frac{1-{\phi}_2}{2k({\phi}_1)}\textbf{u}-\nu({\phi}_1)\frac{{\phi}\textbf{u}_1}{2k({\phi }_1)}\nonumber\\&+\left(\nu({\phi }_2)-\nu({\phi }_2)\right)\frac{1-{\phi}_2}{2k({\phi}_1)}\textbf{u}_1-\nu({\phi }_2)(1-{\phi}_2)\frac{k({\phi }_2)-k({\phi }_1)}{2k({\phi}_1)k({\phi}_2)}\textbf{u}_2=\textbf{0},
\label{vel_s}\\
&\partial_t\phi+\textbf{u}_1\cdot\nabla\phi+\textbf{u}\cdot\nabla\phi_2-\Delta\mu=0,
\label{phi_s}\\
&\partial_t{\psi}+{\textbf{u}}_1\cdot \nabla{\psi}+{\textbf{u}}\cdot \nabla{\psi}_2=0,
\label{psi_s}
\end{align}
almost everywhere in $\Omega\times (0,T_M)$.
We know that
\begin{align}
&\nonumber\Vert\textbf{u}_i\Vert_{L^{3/2}(0,T_M;[W^{2,3}(\Omega)]^2)}+\Vert   \textbf{u}_i\Vert   _{L^\infty(0,T_M;\textbf{V}_\sigma)}+\Vert   \textbf{u}_i\Vert   _{L^2(0,T_M;\textbf{W}_\sigma)}
\\&+\Vert   \phi_i\Vert   _{L^\infty(0,T_M;W^{2,3}(\Omega))}+\Vert   \Psi''(\phi_i)\Vert   _{L^\infty(0,T_M;L^4(\Omega))}\nonumber\\&+\Vert   \psi_i\Vert _{L^\infty(0,T_M;[W^{2,3}(\Omega)]^2)}
\leq C,\quad i=1,2,
\label{stime}
\end{align}
for some constant $C>0$ also depending on $T_M$.
Let us multiply \eqref{vel_s} by $\textbf{u}$, \eqref{phi_s} by $\phi$ and integrate over $\Omega$. Recalling \eqref{nuA}, we obtain
\begin{align}
&\nonumber\frac{ d}{dt}\frac{1}{2}\Vert\textbf{u}\Vert^2+\nu_*\Vert D   \textbf{u}\Vert^2\leq-\underbrace{(({\textbf{u}}\cdot\nabla)\textbf{u}_2,\textbf{u})}_{J_1}-\underbrace{((\nu(\phi_1)-\nu(\phi_1)) D\textbf{u}_2,\nabla\textbf{u})}_{J_2}\\&\nonumber+\underbrace{(\lambda({\phi}_1)\nabla{\psi}_1^T\nabla{\psi},\nabla\textbf{u})}_{J_3}+\underbrace{(((\lambda({\phi}_1)-\lambda({\phi}_2))\nabla{\psi}_1^T\nabla{\psi}_2,\nabla\textbf{u})}_{J_4}\\&\nonumber+\underbrace{(\lambda({\phi}_2)\nabla{\psi}^T\nabla{\psi}_2,\nabla\textbf{u})}_{J_5}+\underbrace{(\nabla{\phi}_1\otimes\nabla{\phi},\nabla\textbf{u})}_{J_6}+\underbrace{(\nabla{\phi}\otimes\nabla{\phi}_2,\nabla\textbf{u})}_{J_7}\\&\nonumber+\underbrace{\left(\nu({\phi }_2)\frac{1-{\phi}_2}{2k({\phi}_1)}\textbf{u},\textbf{u}\right)}_{J_8}+\underbrace{\left(\nu({\phi}_1)\frac{\phi\textbf{u}_1}{2k({\phi }_1)},\textbf{u}\right)}_{J_9}+\underbrace{\left(\left(\nu({\phi }_2)-\nu({\phi }_2)\right)\frac{1-{\phi}_2}{2k({\phi}_1)}\textbf{u}_1,\textbf{u}\right)}_{J_{10}}\nonumber\\&+\underbrace{\left(\nu({\phi }_2)(1-{\phi}_2)\frac{k({\phi }_2)-k({\phi }_1)}{2k({\phi}_1)k({\phi}_2)}\textbf{u}_2,\textbf{u}\right)}_{J_{11}},
\label{vel_s2}\\
&\nonumber\frac{d}{dt}\frac{1}{2}\Vert\phi\Vert^2+\Vert\Delta\phi\Vert^2=-\underbrace{(\textbf{u}\cdot\nabla\phi_2,\phi)}_{J_{12}}+\underbrace{(\Psi^{\prime}(\phi_1)-\Psi^\prime(\phi_2),\Delta\phi)}_{J_{13}}\\&+\underbrace{\left(\lambda^\prime({\phi }_1)-\lambda^\prime({\phi }_2),\frac{1}{2}\vert\nabla{\psi}_1\vert^2\Delta{\phi }\right)}_{J_{14}}+\underbrace{\left(\frac{\lambda^\prime({\phi }_2)}{2}\left(\vert\nabla{\psi}_1\vert^2-\vert\nabla{\psi}_2\vert^2\right),\Delta{\phi}\right)}_{J_{15}}.
\label{phi_s2}
\end{align}
Moreover, we first multiply \eqref{psi_s} by $\psi$ and integrate over $\Omega$ to get
$$
\frac{d}{dt}\frac{1}{2}\Vert\psi\Vert^2+(\textbf{u}\cdot \nabla\psi_2,\psi)=0.
$$
Then we take the gradient of \eqref{psi_s}, multiply it by $\nabla\psi$ and then integrate it over $\Omega$. This gives
$$
\frac{d}{dt}\frac{1}{2}\Vert\nabla\psi\Vert^2+(\nabla\textbf{u}_1\nabla\psi, \nabla\psi)+(\nabla\textbf{u}\nabla\psi_2, \nabla\psi)+(\nabla^2\psi_2\textbf{u}, \nabla\psi).
$$
Adding the two identities above together, we find
\begin{align}
\frac{d}{dt}\frac{1}{2}\Vert\psi\Vert^2_{V}=-\underbrace{(\textbf{u}\cdot \nabla\psi_2,\psi)}_{J_{16}}-\underbrace{(\nabla\textbf{u}_1\nabla\psi, \nabla\psi)}_{J_{17}}-\underbrace{(\nabla\textbf{u}\nabla\psi_2, \nabla\psi)}_{J_{18}}-\underbrace{(\nabla^2\psi_2\textbf{u}, \nabla\psi)}_{J_{19}}.
\end{align}
Therefore, we have
\begin{align}
\frac{1}{2}\frac{d}{dt}\left(\Vert\textbf{u}\Vert^2+\Vert\phi\Vert^2+\Vert\psi\Vert^2_{V}\right)+\nu_*\Vert D   \textbf{u}\Vert^2+\Vert\Delta\phi\Vert^2\leq \sum_{k=1}^{19}J_k.
\label{sumup}
\end{align}
We are thus left to estimate the terms $J_k$. By Gagliardo-Nirenberg and Young's inequalities we get, by \eqref{stime},
$$
\vert J_1\vert \leq \Vert \textbf{u}\Vert_{[L^4(\Omega)]^2}^2\Vert\nabla\textbf{u}_2\Vert\leq C\Vert\textbf{u}\Vert\Vert D\textbf{u}\Vert\leq C\Vert\textbf{u}\Vert^2+\frac{\nu_*}{18}\Vert D\textbf{u}\Vert^2.
$$
Then, thanks to \eqref{V}, we obtain (see again \eqref{stime})
\begin{align*}
\vert J_2\vert &\leq C\Vert\phi\Vert_{L^4(\Omega)}\Vert D\textbf{u}_2\Vert_{[L^4(\Omega)]^4}\Vert D\textbf{u}\Vert\leq C\Vert\phi\Vert_{V}\Vert \textbf{u}_2\Vert_{[H^2(\Omega)]^2}^{1/2}\Vert D\textbf{u}\Vert\\&\leq \frac{\nu_*}{18}\Vert D\textbf{u}\Vert^2+\frac{1}{12}\Vert \Delta\phi\Vert^2+C(1+\Vert\textbf{u}_2\Vert^2_{[H^2(\Omega)]^2})\Vert\phi\Vert^2.
\end{align*}
By similar estimates, recalling \eqref{lambda}, \eqref{stime} and the embedding $W^{2,3}(\Omega)\hookrightarrow W^{1,\infty}(\Omega)$, we deduce
\begin{align*}
\vert J_3\vert \leq C\Vert\nabla\psi_1\Vert_{[L^\infty(\Omega)]^4}\Vert\nabla\psi\Vert\Vert D\textbf{u}\Vert\leq C\Vert\nabla\psi\Vert^2+\frac{\nu_*}{18}\Vert D\textbf{u}\Vert^2,
\end{align*}
and, analogously,
\begin{align*}
\vert J_5\vert\leq C\Vert\nabla\psi\Vert^2+\frac{\nu_*}{18}\Vert D\textbf{u}\Vert^2.
\end{align*}
Then, being $\lambda$ globally Lipschitz, using the embedding $W^{2,3}(\Omega)\hookrightarrow W^{1,\infty}(\Omega)$ and \eqref{stime} we have
$$
\vert J_4\vert\leq C\Vert\phi\Vert\Vert\nabla\psi_1\Vert_{[L^\infty(\Omega)]^4}\Vert\nabla\psi_2\Vert_{[L^\infty(\Omega)]^4}\Vert D\textbf{u}\Vert\\
\leq C\Vert\phi\Vert^2+\frac{\nu_*}{18}\Vert D\textbf{u}\Vert^2
$$
By the same embedding, \eqref{V}, \eqref{stime} and Young's inequality, we infer
\begin{align*}
\vert J_6\vert&\leq C\Vert\nabla\phi_2\Vert_{[L^\infty(\Omega)]^2}\Vert\nabla\phi\Vert\Vert D \textbf{u}\Vert \leq C\Vert\nabla\phi\Vert^2+\frac{\nu_*}{18}\Vert D\textbf{u}\Vert^2\\
&\leq C\Vert\phi\Vert^2+\frac{1}{12}\Vert\Delta\phi\Vert^2+\frac{\nu_*}{18}\Vert D\textbf{u}\Vert^2,
\end{align*}
and
$$
\vert J_7\vert\leq C\Vert\phi\Vert^2+\frac{1}{12}\Vert\Delta\phi\Vert^2+\frac{\nu_*}{18}\Vert D\textbf{u}\Vert^2.
$$
Concerning $J_8$, on account of $\Vert\phi_1\Vert_{L^\infty(\Omega\times(0,T))}\leq 1$, we get
$$
\vert J_6\vert\leq C\Vert\textbf{u}\Vert^2,
$$
whereas for $J_9,J_{10},J_{11}$, by H\"older's inequality, the Sobolev embedding $\textbf{V}_\sigma\hookrightarrow [L^4(\Omega)]^2$ and \eqref{stime}, being $\nu,k\in W^{1,\infty}(\R)$,
we infer that
\begin{align*}
\vert J_9\vert+\vert J_{10}\vert+\vert J_{11}\vert&\leq C\Vert\phi\Vert(\Vert\textbf{u}_1\Vert_{[L^4(\Omega)]^2}+\Vert\textbf{u}_2\Vert_{[L^4(\Omega)]^2})\Vert\textbf{u}\Vert_{[L^4(\Omega)]^2}\\&\leq C(\Vert D\textbf{u}_1\Vert^2+\Vert D\textbf{u}_2\Vert^2)\Vert\phi\Vert^2+\frac{\nu_*}{18}\Vert D\textbf{u}\Vert^2\leq C\Vert\phi\Vert^2+\frac{\nu_*}{18}\Vert D\textbf{u}\Vert^2.
\end{align*}
Using again the embedding $W^{2,3}(\Omega)\hookrightarrow W^{1,\infty}(\Omega)$ and \eqref{stime}, we find
$$
\vert J_{12}\vert\leq \Vert\nabla\phi_2\Vert_{[L^\infty(\Omega)]^2}\Vert\textbf{u}\Vert\Vert\phi\Vert\leq C\Vert\textbf{u}\Vert^2+C\Vert\phi\Vert^2
$$
Observe that, since $\Psi^{\prime\prime}$ is convex, \eqref{stime} and standard embeddings allow us to deduce (see also \eqref{V})
\begin{align*}
J_{13}&= (\Psi'(\phi_1)-\Psi'(\phi_2),\Delta\phi)=\left(\phi\int_0^1\left\{s\Psi''(\phi_1)+(1-s)\Psi''(\phi_2)\right\}ds,\Delta\phi\right)\\&\leq (\Vert   \Psi''(\phi_1)\Vert   _{L^4(\Omega)}+\Vert   \Psi''(\phi_2)\Vert   _{L^4(\Omega)})\Vert   \phi\Vert   _{L^4(\Omega)}\Vert   \Delta\phi\Vert   \leq  C\Vert   \phi\Vert   _V\Vert   \Delta\phi\Vert\\&\leq C\Vert\phi\Vert^{1/2}\Vert\Delta\phi\Vert^{3/2}+C\Vert\phi\Vert\Vert\Delta\phi\Vert\leq C\Vert\phi\Vert^2+\frac{1}{12}\Vert\Delta\phi\Vert^2.
\end{align*}
Consider now $J_{14}$. On account of $\lambda^\prime$ globally Lipschitz, owing to $W^{2,3}(\Omega)\hookrightarrow W^{1,\infty}(\Omega)$ and \eqref{stime}, we have
$$
\vert J_{14}\vert\leq C\Vert\phi\Vert\Vert\nabla\psi_1\Vert^2_{[L^\infty(\Omega)]^2}\Vert\Delta\phi\Vert\leq C\Vert\phi\Vert^2+\frac{1}{12}\Vert\Delta\phi\Vert^2.
$$
Then, by the embedding $W^{2,3}(\Omega)\hookrightarrow W^{1,\infty}(\Omega)$ and again \eqref{stime}, and being $\lambda^\prime$ bounded (see \eqref{lambda}),
\begin{align*}
\vert J_{15}\vert\leq C\left(\Vert\nabla\psi_1\Vert_{[L^\infty(\Omega)]^4}+\Vert\nabla\psi_2\Vert_{[L^\infty(\Omega)]^4}\right)\Vert\nabla\psi\Vert\Vert\Delta\phi\Vert\leq C\Vert\psi\Vert_V^2+\frac{1}{12}\Vert\Delta\phi\Vert^2.
\end{align*}
Arguing similarly, we get
$$
\vert J_{16}\vert\leq C\Vert\nabla\psi_2\Vert_{[L^\infty(\Omega)]^4}\Vert\textbf{u}\Vert\Vert\psi\Vert\leq C\Vert\textbf{u}\Vert^2+C\Vert\psi\Vert_V^2,
$$
whereas, by H\"older's inequality, we obtain
$$
\vert J_{17}\vert\leq \Vert\nabla\textbf{u}_1\Vert_{[L^\infty(\Omega)]^4}\Vert\nabla\psi\Vert^2\leq C\Vert \nabla\textbf{u}_1\Vert_{[L^\infty(\Omega)]^4}\Vert\psi\Vert^2_V.
$$
On the other hand, by \eqref{stime} and $W^{2,3}(\Omega)\hookrightarrow W^{1,\infty}(\Omega)$, we deduce
$$
\vert J_{18}\vert\leq C\Vert\nabla\psi_2\Vert\Vert D\textbf{u}\Vert\Vert\psi\Vert_V\leq C\Vert\psi\Vert_V^2+\frac{\nu_*}{18}\Vert D\textbf{u}\Vert^2.
$$
Finally, using \eqref{stime} once more and H\"older and Gagliardo-Nirenberg's inequalities, we find
\begin{align*}
\vert J_{19}\vert&\leq C\Vert\psi_2\Vert_{[W^{2,3}(\Omega)]^2}\Vert \textbf{u}\Vert_{[L^{6}(\Omega)]^2} \Vert\nabla\psi\Vert\leq C\Vert\textbf{u}\Vert^{1/3}\Vert D\textbf{u}\Vert^{2/3}\Vert\psi\Vert_V\\&\leq C\Vert\psi\Vert^2_V+\frac{\nu_*}{18}\Vert D\textbf{u}\Vert^{2}+C\Vert\textbf{u}\Vert^2.
\end{align*}
Collecting all the above estimates, we deduce the inequality
\begin{align*}
\frac{d}{dt}\mathcal{H}+\frac{\nu_*} {2}\Vert D   \textbf{u}\Vert^2+\frac{1}{2}\Vert\Delta\phi\Vert^2\leq C\left(1+\Vert \nabla\textbf{u}_1\Vert_{[L^\infty(\Omega)]^4}+\Vert\textbf{u}_2\Vert^2_{[H^2(\Omega)]^2}\right)\mathcal{H},
\end{align*}
for almost any $t\in (0,T_M)$, where
$$
\mathcal{H}(t):=\Vert\textbf{u}(t)\Vert^2+\Vert\phi(t)\Vert^2+\Vert\psi(t)\Vert^2_{V}.
$$
Thus, thanks to \eqref{stime} and the embedding $[W^{2,3}(\Omega)]^2\hookrightarrow [W^{1,\infty}(\Omega)]^2$,
Gronwall's Lemma entails a continuous dependence estimate with respect to the initial data and, in particular, the
uniqueness of the local (in time) solution.

\section{Proofs of the Theorems in Section \ref{three}}
\label{proofs}
We give here the proofs of the auxiliary results stated in Section \ref{three}.
\subsection{Proof of Lemma \ref{transport}}
\label{proof1}
The existence of a map $X(s,t,x)$ with the desired properties and regularity in $s$ and $x$ is a well known result. We need to establish the regularity of $X(0,t,x)$ with respect to $t$.
Observe that, using the Einstein notation, for $i,j,k=1,2$,
$$
\partial_t X_k(0,t,x)=\textbf{v}(x,t)+\int_{t}^0\partial_l\textbf{v}_k(X(\tau,t,x),\tau)\partial_t X_l(\tau,t,x)d\tau,
$$
$$
\partial_j X_k(0,t,x)=-\delta_{jk}+\int_{t}^0\partial_l\textbf{v}_k(X(\tau,t,x),\tau)\partial_j X_l(\tau,t,x)d\tau,
$$

\begin{align*}
\partial_t\partial_{j} X_k(0,t,x)&=\int_{t}^0\partial_n\partial_l\textbf{v}_k(X(\tau,t,x),\tau)\partial_j X_l(\tau,t,x)\partial_t X_n(\tau,t,x)d\tau\\&+\int_{t}^0\partial_l\textbf{v}_k(X(\tau,t,x),\tau)\partial_t\partial_j X_l(\tau,t,x)d\tau,
\end{align*}
\begin{align*}
\partial_i\partial_{j} X_k(0,t,x)&=\int_{t}^0\partial_n\partial_l\textbf{v}_k(X(\tau,t,x),\tau)\partial_j X_l(\tau,t,x)\partial_i X_n(\tau,t,x)d\tau\\&+\int_{t}^0\partial_l\textbf{v}_k(X(\tau,t,x),\tau)\partial_i\partial_j X_l(\tau,t,x)d\tau.
\end{align*}
Considering, e.g., the first equality, we observe that if we look for $z_k$, $k=1,2$, such that
$$
\partial_t z_k(0,t,x)=\textbf{v}(x,t)+\int_{t}^0\partial_l\textbf{v}_k(X(\tau,t,x),\tau)z_l(\tau,t,x)d\tau,
$$
following \cite[Lemma 1.2,][]{Lady} we deduce that it exists and it is unique in $C([0,T];C(\overline{\Omega}))$ and also $z_k=\partial_t X_k$ for $k=1,2$. In a similar way, note that we only need $\textbf{v}	\in C([0,T];C^2(\overline{\Omega}))$, we infer that $\partial_j X_k$, $\partial_t\partial_{j} X_k$ and $\partial_i\partial_{j} X_k$ are also continuous in $[0,T]\times 	\overline{\Omega}$ (i.e, with respect to $t$ and $x$), meaning that, for $k=1,2$,
$$
X_k(0,\cdot,\cdot)\in C([0,T];C^2(\overline{\Omega}))\cap C^1([0,T];C^1(\overline{\Omega})).
$$
Thus, we get \eqref{psireg}. If we now argue as in the proof of \cite[Lemma 1.1,][]{Lady}, we obtain \eqref{estW1p} and \eqref{estW2p}.
If, in addition, we have $\textbf{v}	\in C([0,T];C^3(\overline{\Omega}))$, then, for $i,j,k=1,2$,
\begin{align*}
\partial_t\partial_i\partial_{j} X_k(0,t,x)&=\int_{t}^0\partial_m\partial_n\partial_l\textbf{v}_k(X(\tau,t,x),\tau)\partial_j X_l(\tau,t,x)\partial_i X_n(\tau,t,x)\partial_t X_m(\tau,t,x)d\tau\\&+\int_{t}^0\partial_n\partial_l\textbf{v}_k(X(\tau,t,x),\tau)\partial_t\partial_j X_l(\tau,t,x)\partial_i X_n(\tau,t,x)d\tau\\&+\int_{t}^0\partial_n\partial_l\textbf{v}_k(X(\tau,t,x),\tau)\partial_j X_l(\tau,t,x)\partial_t\partial_i X_n(\tau,t,x)d\tau
\\&+\int_{t}^0\partial_m\partial_l\textbf{v}_k(X(\tau,t,x),\tau)\partial_i\partial_j X_l(\tau,t,x)\partial_t X_m(\tau,t,x)d\tau\\&+\int_{t}^0\partial_l\textbf{v}_k(X(\tau,t,x),\tau)\partial_t\partial_i\partial_j X_l(\tau,t,x)d\tau,
\end{align*}
from which, arguing again as in the proof of \cite[Lemma 1.2,][]{Lady}, we deduce that, for $k=1,2$,
$$
X_k(0,\cdot,\cdot)\in C^1([0,T];C^2(\overline{\Omega})),
$$
and this implies \eqref{psireg2} since $\psi_0$ is smooth enough.
In conclusion, in force of \eqref{psireg2}, we can rigorously repeat verbatim the estimate \cite[(3.17),][]{GGW} and obtain \eqref{psiw2p}, for any $p\in(2,\infty)$. Indeed, we have at least $\Vert\psi\Vert_{[W^{2,p}(\Omega)]^2}\in AC([0,T])$, which is enough to get that estimate.
\subsection{Proof of Theorem \ref{ns}}
\label{proof2}
The uniqueness of the solutions is immediate, being the problem linear (it is enough to exploit $\textbf{u}_1-\textbf{u}_2$ as test function in the equation solved by the difference of two solutions $\textbf{u}_1$ and $\textbf{u}_2$), whereas the proof of existence is a standard application of a Galerkin scheme. We consider the family $\{\textbf{w}_j\}_{j\geqslant1}$ of the eigenfunctions of the Stokes operator $\textbf{A}$ as a Galerkin basis in $\textbf{V}_\sigma$ and $\textbf{H}_\sigma$. Then we define the $n$-dimensional subspace
$$
\mathbb{W}_n:=Span(\textbf{w}_1,\ldots,\textbf{w}_n)
$$
and the related orthogonal projector on $\mathbb{W}_n$ in $\textbf{H}_\sigma$, that is, $P_n := P_{\mathbb{W}_n}$. We then look for a function of the form
\begin{align*}
\textbf{u}_n(t)=\sum_{i=1}^{n}{\gamma}_i(t)\textbf{w}_i\in\mathbb{W}_n,
\end{align*}
where $\gamma_i$ are real valued functions, satisfying
\begin{align}
\label{u}
\begin{cases}
(\partial_t\textbf{u}_n,\textbf{w})+((\textbf{v}\cdot \nabla)\textbf{u}_n,\textbf{w})+(\nu(\varphi) D\textbf{u}_n,\nabla\textbf{w})+(g\textbf{u}_n,\textbf{w})=((I+\alpha^2\textbf{A})^{-2}\textbf{h},\textbf{w}),\; \forall \textbf{w}\in\mathbb{W}_n\\	
\textbf{u}_n(0)=P_n(\textbf{u}_0)
\end{cases}
\end{align}
for every $t\in(0,T)$.   From now on, the generic constant $C>0$ could depend on $\alpha$, but neither on $n$ nor on $t$.
Being the problem linear, we can globally solve the Cauchy problem for the system in the unknowns $\gamma_i$ and find a unique maximal solution  $\gamma^{(n)}\in C^1([0, T ), \R^n)$, for any $T>0$. We can now derive some uniform estimates. In particular, considering $\textbf{w}=\textbf{u}_n$ and applying Young's inequality together with \eqref{nuA} and \eqref{h1v}, we get
$$
\frac{1}{2}\frac{d}{dt}\Vert\textbf{u}_n\Vert^2+{\nu_*}\Vert D\textbf{u}_n\Vert^2+(g,\vert\textbf{u}_n\vert^2)\leq C\Vert\textbf{h}\Vert^2,
$$
from which we deduce, being $g\geq0$,
\begin{equation}
\Vert\textbf{u}_n\Vert_{L^\infty(0,T;\textbf{H}_\sigma)}+\Vert\textbf{u}_n\Vert_{L^2(0,T;\textbf{V}_\sigma)}\leq C.
\label{firstestimates}
\end{equation}
If we now consider $\textbf{w}=\partial_t\textbf{u}_n\in \mathbb{W}_n$ in \eqref{u}, then we obtain
\begin{align*}
&\Vert\partial_t\textbf{u}_n\Vert^2-(\nu^\prime(\varphi)D\textbf{u}_n\nabla\varphi,\partial_t\textbf{u}_n)-\frac{1}{2}(\nu(\varphi)\Delta\textbf{u}_n,\partial_t\textbf{u}_n)+((\textbf{v}\cdot \nabla)\textbf{u}_n,\partial_t\textbf{u}_n)+(g\textbf{u}_n,\partial_t\textbf{u}_n )\\&=((I+\alpha^2\textbf{A})^{-2}\textbf{h}, \partial_t\textbf{u}_n).
\end{align*}
Now, being $\nu\in W^{1,\infty}(\R)$ and recalling the assumptions on $\varphi$, we have
\begin{align*}	
(\nu^\prime(\varphi)D\textbf{u}_n\nabla\varphi,\partial_t\textbf{u}_n)&\leq C\Vert D\textbf{u}_n\Vert_{[L^4(\Omega)]^4}\Vert\nabla\varphi\Vert_{[L^4(\Omega)]^2}\Vert\partial_t\textbf{u}_n\Vert\\&\leq C\Vert\textbf{A}\textbf{u}_n\Vert^{1/2}\Vert D\textbf{u}_n\Vert^{1/2}\Vert\nabla\varphi\Vert^{1/2}\Vert\varphi\Vert^{1/2}_{H^2(\Omega)}\Vert\partial_t\textbf{u}_n\Vert\\
&\leq C\Vert\textbf{A}\textbf{u}_n\Vert^{1/2}\Vert D\textbf{u}_n\Vert^{1/2}\Vert\varphi\Vert^{1/2}_{H^2(\Omega)}\Vert\partial_t\textbf{u}_n\Vert\\&\leq \frac {1}{10}\Vert\partial_t\textbf{u}_n\Vert^2
+ C\Vert \textbf{A}\textbf{u}_n\Vert^2+C\Vert\varphi\Vert^{2}_{H^2(\Omega)}\Vert D\textbf{u}_n\Vert^{2}.
\end{align*}
On the other hand, recalling \eqref{nuA}, we get
\begin{align*}
\frac{1}{2}(\nu(\varphi)\Delta\textbf{u}_n,\partial_t\textbf{u}_n)\leq C\Vert\textbf{A}\textbf{u}_n\Vert\Vert\partial_t\textbf{u}_n\Vert\leq
C\Vert\textbf{Au}_n\Vert^2+\frac{1}{10}\Vert\partial_t\textbf{u}_n\Vert^2,
\end{align*}
Therefore, we find
$$
((\textbf{v}\cdot \nabla)\textbf{u}_n,\partial_t\textbf{u}_n)\leq C\Vert\textbf{v}\Vert_{[L^\infty(\Omega)]^2}\Vert D\textbf{u}_n\Vert\Vert\partial_t\textbf{u}_n\Vert\leq \frac{1}{10}\Vert\partial_t\textbf{u}_n\Vert^2+C\Vert\textbf{v}\Vert_{[L^\infty(\Omega)]^2}^2\Vert D\textbf{u}_n\Vert^2.
$$
Moreover, observe that (see \eqref{firstestimates})
$$
(g\textbf{u}_n,\partial_t\textbf{u}_n )\leq \Vert g\Vert_{L^\infty(\Omega)}\Vert \textbf{u}_n\Vert\Vert\partial_t\textbf{u}_n\Vert\leq \frac{1}{10}\Vert\partial_t\textbf{u}_n\Vert^2+C
$$
and, by Young's inequality and \eqref{dual}, that
$$
((I+\alpha^2\textbf{A})^{-2}\textbf{h}, \partial_t\textbf{u}_n)\leq \frac{1}{10}\Vert\partial_t\textbf{u}_n\Vert^2+C_\alpha\Vert\textbf{h}\Vert_{\textbf{V}_\sigma^\prime}^2
$$
Summing up, we have, independently of $n$,
\begin{align}
\nonumber
&\frac{1}{2}\Vert\partial_t\textbf{u}_n\Vert^2\\
&\leq C\Vert \textbf{A}\textbf{u}_n\Vert^2+C\Vert\varphi\Vert^{2}_{H^2(\Omega)}\Vert D\textbf{u}_n\Vert^{2}+ C\Vert\textbf{v}\Vert_{[L^\infty(\Omega)]^2}^2\Vert D\textbf{u}_n\Vert^2+C_\alpha\Vert\textbf{h}\Vert_{\textbf{V}_\sigma^\prime}^2+C.
\label{dtn}
\end{align}
We now choose $\textbf{w}=\textbf{A}\textbf{u}_n\in \mathbb{W}_n$ in \eqref{u}. This gives, being $(I+\alpha^2\textbf{A})^{-2}$ self-adjoint,
\begin{align}
&\nonumber(\partial_t\textbf{u}_n,\textbf{Au}_n)-(\text{div}(\nu(\varphi)D\textbf{u}_n),\textbf{A}\textbf{u}_n)+((\textbf{v}\cdot \nabla)\textbf{u}_n,\textbf{Au}_n)\\&+(g\textbf{u}_n,\textbf{Au}_n)=(\textbf{h},(I+\alpha^2\textbf{A})^{-2}\textbf{Au}_n).
\label{Au}
\end{align}
Recalling that (see \cite[Appendix B]{Giorginitemam}), there exists $p_n\in L^2(0,T;V)$ such that
$$
-\Delta\textbf{u}_n+\nabla p_n=\textbf{A}\textbf{u}_n,
$$
almost everywhere in $\Omega\times(0,T)$, with
\begin{align}
\Vert p_n\Vert\leq C\Vert\nabla\textbf{u}_n\Vert^{1/2}\Vert\textbf{A}\textbf{u}_n\Vert^{1/2},\qquad \Vert p_n\Vert_V\leq C\Vert\textbf{A}\textbf{u}_n\Vert.
\label{pn}
\end{align}
Hence, we get
\begin{align*}
&-(\text{div}(\nu(\varphi)D\textbf{u}_n),\textbf{A}\textbf{u}_n)=-(\nu^\prime(\varphi)D\textbf{u}_n\nabla\varphi,\textbf{A}\textbf{u}_n)-\frac{1}{2}(\nu(\varphi)\Delta\textbf{u}_n,\textbf{Au}_n)\\&=-(\nu^\prime(\varphi)D\textbf{u}_n\nabla\varphi,\textbf{A}\textbf{u}_n)+\frac{1}{2}(\nu(\varphi)\textbf{Au}_n,\textbf{Au}_n)-\frac{1}{2}(\nu(\varphi)\nabla p_n,\textbf{Au}_n)\\&\geq-(\nu^\prime(\varphi)D\textbf{u}_n\nabla\varphi,\textbf{A}\textbf{u}_n)+\frac{\nu_*}{2}\Vert\textbf{Au}_n\Vert^2+\frac{1}{2}(\nu^\prime(\varphi)\nabla\varphi p_n,\textbf{Au}_n).
\end{align*}
Exploiting the regularity of $\varphi$, being $\nu\in W^{1,\infty}(\R)$, we obtain
\begin{align*}
&(\nu^\prime(\varphi)D\textbf{u}_n\nabla\varphi,\textbf{A}\textbf{u}_n)\leq  C\Vert D\textbf{u}_n\Vert_{[L^4(\Omega)]^4}\Vert\nabla\varphi\Vert_{[L^4(\Omega)]^2}\Vert\textbf{A}\textbf{u}_n\Vert\\&\leq C\Vert\textbf{A}\textbf{u}_n\Vert^{1/2}\Vert D\textbf{u}_n\Vert^{1/2}\Vert\nabla\varphi\Vert^{1/2}\Vert\varphi\Vert^{1/2}_{H^2(\Omega)}\Vert\textbf{A}\textbf{u}_n\Vert\leq C\Vert\textbf{A}\textbf{u}_n\Vert^{3/2}\Vert D\textbf{u}_n\Vert^{1/2}\Vert\varphi\Vert^{1/2}_{H^2(\Omega)}\\&\leq \frac{\nu_*}{16}\Vert \textbf{A}\textbf{u}_n\Vert^2+C\Vert\varphi\Vert^{2}_{H^2(\Omega)}\Vert D\textbf{u}_n\Vert^{2},		
\end{align*}
and, by \eqref{pn},
\begin{align*}
\frac{1}{2}(\nu^\prime(\varphi)\nabla\varphi p_n,\textbf{Au}_n)&\leq C\Vert\nabla\varphi\Vert_{[L^4(\Omega)]^2}\Vert p_n\Vert_{L^4(\Omega)}\Vert\textbf{Au}_n\Vert\\&\leq C\Vert\varphi\Vert_{H^2(\Omega)}^{1/2}\Vert p_n\Vert^{1/2}\Vert p_n\Vert^{1/2}_V\Vert\textbf{Au}_n\Vert\leq C\Vert\varphi\Vert_{H^2(\Omega)}^{1/2}\Vert D\textbf{u}_n\Vert^{1/4}\Vert\textbf{A}\textbf{u}\Vert^{3/4}\\&\leq \frac{\nu_*}{16}\Vert\textbf{A}\textbf{u}\Vert^{2}+C\Vert\varphi\Vert_{H^2(\Omega)}^4\Vert D\textbf{u}_n\Vert^2.
\end{align*}
Therefore, we infer that
$$
((\textbf{v}\cdot \nabla)\textbf{u}_n,\textbf{Au}_n)\leq C\Vert\textbf{v}\Vert_{[L^\infty(\Omega)]^2}\Vert D\textbf{u}_n\Vert\Vert\textbf{Au}_n\Vert\leq \frac{\nu_*}{16}\Vert\textbf{Au}_n\Vert^2+C\Vert\textbf{v}\Vert_{[L^\infty(\Omega)]^2}^2\Vert D\textbf{u}_n\Vert^2
$$
and (see \eqref{calpha1})
\begin{align*}
&(\textbf{h},(I+\alpha^2\textbf{A})^{-2}\textbf{Au}_n)\leq \Vert\textbf{h}\Vert_{\textbf{V}_\sigma^\prime}\Vert(I+\alpha^2\textbf{A})^{-2}\textbf{Au}_n\Vert_{\textbf{V}_\sigma}\\&\leq C_\alpha\Vert\textbf{h}\Vert_{\textbf{V}_\sigma^\prime}\Vert\textbf{Au}_n\Vert\leq \frac{\nu_*}{16}\Vert\textbf{Au}_n\Vert^2+C\Vert\textbf{h}\Vert_{\textbf{V}_\sigma^\prime}^2.
\end{align*}
Thus, we deduce from \eqref{Au} that
\begin{align*}
\frac{1}{2}\frac{d}{dt}\Vert\nabla\textbf{u}\Vert^2+\frac{\nu_*}{4}\Vert\textbf{Au}_n\Vert^2\leq C\left(1+\Vert\varphi\Vert_{H^2(\Omega)}^4+\Vert\textbf{v}\Vert_{[L^\infty(\Omega)]^2}^2\right)\Vert\nabla\textbf{u}_n\Vert^2+C\Vert\textbf{h}\Vert_{\textbf{V}_\sigma^\prime}^2.
\end{align*}
Adding this estimate and \eqref{dtn} multiplied by $\omega=\frac{\nu_*}{16C}$ together, we obtain
\begin{align*}
&\frac{1}{2}\frac{d}{dt}\Vert\nabla\textbf{u}\Vert^2+\frac{\omega}{2}\Vert\partial_t\textbf{u}_n\Vert^2+\frac{\nu_*}{8}\Vert\textbf{Au}_n\Vert^2\\&\leq C\left(1+\Vert\varphi\Vert_{H^2(\Omega)}^4+\Vert\textbf{v}\Vert_{[L^\infty(\Omega)]^2}^2\right)\Vert\nabla\textbf{u}_n\Vert^2+C\left(1+\Vert\textbf{h}\Vert_{\textbf{V}_\sigma^\prime}^2\right),
\end{align*}
and Gronwall's Lemma yields
\begin{equation}
\Vert\partial_t\textbf{u}_n\Vert_{L^2(0,T;\textbf{H}_\sigma)}+\Vert \textbf{u}_n\Vert_{L^\infty(0,T;\textbf{V}_\sigma)}+\Vert\textbf{A}\textbf{u}_n\Vert_{L^2(0,T;[L^2(\Omega)]^2)}\leq C_\alpha(T)
\label{est2}
\end{equation}
so that (see \eqref{est2})
$$
\Vert\textbf{u}_n\Vert_{L^2(0,T;\textbf{W}_\sigma)}\leq C_\alpha(T).
$$
This concludes the proof, by passing to the limit as $n$ goes to $\infty$ along a suitable subsequence by means of standard compactness arguments.
\subsection{Proof of Theorem \ref{ch}}
\label{proof3}
Let us prove uniqueness first. We consider two weak solutions to \eqref{CH}, $\phi_i$, $i=1,2$, corresponding to initial data $\phi_{0,i}$ such that $\overline{\phi}_{0,1}=\overline{\phi}_{0,2}$. We set $\phi=\phi_1-\phi_2$, $\mu=-\Delta \phi+\Psi^\prime(\phi_1)-\Psi^\prime(\phi_2)+\left(\lambda_\alpha^\prime(\phi_1)-\lambda_\alpha^\prime(\phi_2)\right)f$. Then we have
\begin{align}
&(\partial_t \phi,\varphi)+ (\textbf{v}\cdot \nabla\phi, \varphi)+(\nabla\mu,\nabla\varphi)=0,	\quad\forall\,\varphi \in V.
\label{cahn2}
\end{align}
Therefore, testing against $\varphi=A^{-1}\phi$ (note that $\phi\in V_0$), we get
\begin{align*}
\nonumber&\frac{1}{2}\frac{d}{dt}\Vert\phi\Vert_\sharp^2+
(\textbf{v}\cdot \nabla\phi, A^{-1}\phi)+(\mu,\phi)=0,
\end{align*}
that is,
\begin{equation}
\frac{1}{2}\frac{d}{dt}\Vert\phi\Vert_\sharp^2+\Vert\nabla\phi\Vert^2+
(\textbf{v}\cdot \nabla\phi, A^{-1}\phi)+(F^\prime(\phi_1)-F^\prime(\phi_2),\phi)+\left(\lambda_\alpha^\prime(\phi_1)-\lambda_\alpha^\prime(\phi_2),f\phi\right)=\alpha_0\Vert\phi\Vert^2.
\label{uniqCH}
\end{equation}
By strict convexity of $F$, we immediately deduce
\begin{equation*}
(F^\prime(\phi_1)-F^\prime(\phi_2),\phi)\geq\overline{\alpha}\Vert\phi\Vert^2\geq 0,
\end{equation*}
and, being $\lambda_\alpha^\prime$ globally Lipschitz,
\begin{align*}
&\left(\lambda_\alpha^\prime(\phi_1)-\lambda_\alpha^\prime(\phi_2),f\phi\right)\leq C(f,\phi^2)\leq C\Vert f \Vert_{L^\infty(\Omega)} \Vert\phi\Vert^2\\&\leq C\Vert f \Vert_{L^\infty(\Omega)}\Vert\phi\Vert_{\sharp}\Vert\nabla\phi\Vert\leq C\Vert f\Vert^2_{L^\infty(\Omega)}\Vert\phi\Vert_\sharp^2+\frac{1}{4}\Vert\nabla\phi\Vert^2.
\end{align*}
Therefore, by Young's and Poincar\'{e}'s inequalities, we get
\begin{align*}
&\vert(\textbf{v}\cdot \nabla\phi, A^{-1}\phi)\vert\leq  \Vert\textbf{v}\Vert_{[L^\infty(\Omega)]^2}\Vert\phi\Vert\Vert\phi\Vert_{\sharp}\\&\leq C\Vert\textbf{v}\Vert_{[L^\infty(\Omega)]^2}\Vert\nabla\phi\Vert\Vert\phi\Vert_{\sharp}\leq \frac{1}{4}\Vert\nabla\phi\Vert^2+C\Vert\textbf{v}\Vert_{[L^\infty(\Omega)]^2}^2\Vert\phi\Vert_{\sharp}^2.
\end{align*}
Then, collecting the above estimates, we deduce from \eqref{uniqCH} that
$$
\frac{1}{2}\frac{d}{dt}\Vert\phi\Vert_\sharp^2+\frac{1}{2}\Vert\nabla\phi\Vert^2\leq C\left(\Vert\textbf{v}\Vert_{[L^\infty(\Omega)]^2}^2+\Vert f\Vert^2_{L^\infty(\Omega)}\right)\Vert\phi\Vert^2_\sharp,
$$
and uniqueness follows via Gronwall's Lemma.

Let us now prove the existence of a solution. For any $\xi>0$ we introduce $\Psi_\xi(s)=F_\xi(s)-\frac{\alpha_0}{2}s^2$, with $F_\xi$ defined as in \cite[(4.3)]{Giorginitemam} and we recall that $F_\xi$ enjoys the following properties: there exist $0<\overline{\xi}<\gamma\leq1$, $\hat{C}>0$ and $L>0$ such that
\begin{enumerate}
	\item $\Psi_\xi(s)\geq\frac{\alpha_0}{2}s^2-\hat{C}\geq -\hat{C},\; F_\xi(s)\geq 0,\quad  \forall\,\xi\in(0,\overline{\xi}],\;\forall s\in\R$;
	
	\label{Psi1}
	\item  $\Psi_\xi\in\mathcal{C}^{2}(\R), \quad \forall \,\xi \in (0,\overline{\xi}]$;
	\label{apprentr}
	\item $F^{(j)}_\xi(s)=F^{(j)}(s), \quad \forall\,s\in [-1+\xi, 1-\xi]$, $j=0,1,2$;
	\label{2}
	\item as $\xi\rightarrow0$, $F_\xi(s)\rightarrow F(s)$ for all $s\in\R$, $ \vert F'_\xi(s) \vert \rightarrow \vert F'(s) \vert $ for $s\in(-1,1)$ and $F_\xi'$ converges uniformly to $F'$ on any compact subset of $(-1,1)$; \item $ \vert F'_\xi(s) \vert \rightarrow+\infty$ for every $ \vert s \vert \geq1$;
\item for $\xi\in (0,\overline{\xi}]$, $F_\xi(s)\leq F(s)$, for every $s\in[-1,1]$ and $\vert F'_\xi(s) \vert \leq  \vert F'(s) \vert$, for every $s\in (-1,1)$;
	\label{4}
	\item $F_\xi''(s) \geq \overline{\alpha},\quad \forall\, s\in\R,\;\forall\xi\in(0,\overline{\xi}]$;
	\item $\Psi_\xi^{\prime\prime}(z)\leq L, \quad\forall\, z\in \R$.
\end{enumerate}
Following \cite{Giorginitemam}, we introduce the globally Lipschitz function $h_r:\R\rightarrow\R,\ r\in\N$ such that
\begin{equation*}
g_r(z)=\begin{cases}
-r,\quad z<-r,\\z, \quad z\in[-r,r],\\r,\quad z>r.
\end{cases}
\end{equation*}
Then we define $\tilde\mu_{0,r}=g_r\circ\tilde{\mu}_0$, where $\tilde{\mu}_0=-\alpha \Delta\phi_0+F'(\phi_0)=\mu_0+\alpha_0\phi_0$.
Since $\tilde{\mu}_0\in V$, we have $\tilde{\mu}_{0,r}\in  V$, for any $r>0$, and $\nabla\tilde{\mu}_{0,r}=\nabla\tilde{\mu}_0\cdot \chi_{{\color{black}[-r,r]}}(\tilde{\mu}_0)$. This in turn gives
\begin{equation}
\Vert   \tilde{\mu}_{0,r}\Vert   _V\leq\Vert   \tilde{\mu}_0\Vert   _V.
\label{pb}
\end{equation}
For a given $r\in\N$, we consider the Neumann boundary value problem
\begin{equation}
\begin{cases}
-\alpha\Delta\phi_{0,r}+F'(\phi_{0,r})=\tilde{\mu}_{0,r},\quad \text{in }\Omega, \\
\partial_{\textbf{n}}\phi_{0,r}=0, \quad  \text{on }\partial\Omega.
\end{cases}
\label{problem1}
\end{equation}
Recalling \cite[Lemma A.1]{Giorginitemam}, we know that there exists a unique (strong) solution to (\ref{problem1}) such that $\phi_{0,r}\in H^2(\Omega)$, $F'(\phi_{0,r})\in H$. In addition, by \cite[Thm.A.2,][]{Giorginitemam} and (\ref{pb}) we get
\begin{equation}
\Vert   \phi_{0,r}\Vert   _{ H^2(\Omega)}\leq C(1+\Vert   \tilde{\mu}_0\Vert   ).
\label{p1}
\end{equation}
Also, observe that $\tilde{\mu}_{0,r}\rightarrow\tilde{\mu}_0$ in $ H$ entails $\phi_{0,r}\rightarrow\phi_0$ in $ V$ as $r \to \infty$ (see \cite[Lemma A.1]{Giorginitemam}).
As a consequence, there exists an $\tilde{m}\in(0,1)$, independent of $r$, and $\tilde{k}$ sufficiently large such that
\begin{equation}
\Vert   \phi_{0,r}\Vert_V\leq 1+\Vert   \phi_0\Vert_V,\quad \vert \overline{\phi}_{0,r} \vert \leq \tilde{m}<1,\quad \forall r>\tilde{k}.
\label{phik}
\end{equation}
On the other hand, from \cite[Thm.A.2,][]{Giorginitemam}, with $f=\tilde{\mu}_{0,r}$, we obtain
\begin{equation*}
\Vert   F'(\phi_{0,r})\Vert   _{L^\infty(\Omega)}\leq \Vert   \tilde{\mu}_{0,r}\Vert   _{L^\infty(\Omega)}\leq {\color{black}r}.
\end{equation*}
In conclusion, being $\phi_{0,r}\in C(\overline{\Omega})$, we can say that there exists $\delta=\delta(r)>0$ such that
\begin{equation}
\Vert   \phi_{0,r}\Vert   _{L^\infty(\Omega)}\leq 1-\delta.
\label{deltas}
\end{equation}
Observe now that $F''$ is continuous on $(-1,1)$, thus it is bounded on compact sets, so that
\begin{equation*}
\nabla F'(\phi_{0,r})=F''(\phi_{0,r})\nabla \phi_{0,r}\in H.
\end{equation*}
Then, being $F'(\phi_{0,r})\in H$, we deduce that $F'(\phi_{0,r})\in V$.
Thus $\phi_{0,r}\in H^3(\Omega)$.
We also need to note that, for any $\xi\in (0,\xi^*)$, where $\xi^*=min\left\{\frac{1}{2}\delta(r),\overline{\xi}\right\}$, since $F(z)=F_\xi(z)$
for all $z\in[-1+\xi,1-\xi]$ (see \cite{frigeri}), we infer from (\ref{deltas}) that $-\alpha\Delta\phi_{0,r}+F'_\xi(\phi_{0,r})=\tilde{\mu}_{0,r}$, which entails
\begin{equation*}
\Vert   -\alpha\Delta\phi_{0,r}+F_\xi'(\phi_{0,r})\Vert_V\leq\Vert   \tilde{\mu}_0\Vert_V.
\end{equation*}

We now introduce a Galerkin scheme for an approximated problem with $\Psi_\xi$ in place of $\Psi$.
Consider the family $\{\psi_j\}_{j\geqslant1}$ of the  eigenfunctions of the Laplace operator with homogeneous
Neumann boundary conditions as a basis in $V$. Then we define the $n$-dimensional subspaces
$$
\mathbb{Z}_n:=Span(\psi_1,\ldots,\psi_n),
$$
where $\psi_1\equiv1/\sqrt{ \vert \Omega \vert }$ and the related orthogonal projector on these subspaces in $H$, that is, $\tilde{P}_n := P_{\mathbb{Z}_n}$.
We then look for two functions of the form
\begin{align*}
\phi^n_{r,\xi}(t)=\sum_{i=1}^{n}\beta_i(t)\psi_i\in\mathbb{Z}_n\qquad \mu^n_{r,\xi}(t)=\sum_{i=1}^{n}\delta_i(t)\psi_i\in\mathbb{Z}_n,
\end{align*}
where $\beta_i,\delta_i$ are some real valued functions, which solve
the following ODE Cauchy problem
\begin{align}
&\label{eqch}(\partial_t\phi^n_{r,\xi},v)+(\nabla\mu^n_{r,\xi},\nabla v)+(\textbf{v}\cdot\nabla\phi^n_{r,\xi},v)=0\ \ \ \ \forall v\in\mathbb{Z}_n,\\
& \phi^n_{r,\xi}(0)=\tilde{P}_n(\phi_{0,r}),\\
\label{mu}
&\mu^n_{r,\xi}
=-\Delta\phi^n_{r,\xi}+\tilde{P}_n(\Psi_\xi'(\phi^n_{r,\xi})+\lambda^\prime_\alpha(\phi^n_{r,\xi})f),
\end{align}
for every $t\in(0,T)$.

Note that the chosen basis for $ V$ is still a complete family in
$$
\left\{u\in H^3(\Omega): \partial_{\textbf{n}}u=0, \quad \text{on }\partial\Omega\right\}.
$$
Thus we have
\begin{equation*}
\phi^n_{r,\xi}(0)\rightarrow\phi_{0,r}\ \ \ \text{ in } H^3(\Omega).
\end{equation*}
In turn, by the embedding $H^3(\Omega)\hookrightarrow L^\infty(\Omega)$, we get
\begin{equation*}
\phi^n_{r,\xi}(0)\rightarrow\phi_{0,r}\ \ \ \text{ in } L^\infty(\Omega).
\end{equation*}
Hence there exists $\overline{n}=\overline{n}(r)$ such that

\begin{equation}
\Vert   \phi^n_{r,\xi}(0)\Vert   _{\infty}\leq\frac{1}{2}\delta(r)+\Vert   \phi_{0,r}\Vert   _\infty\leq1-\frac{1}{2}\delta(r), \quad\forall\, n>\overline{n}.
\label{eqq}
\end{equation}
{\color{black}For any $r>r_0$ ($r_0$ independent of $n$ and $\xi$) we fix $\xi\in(0,\xi^*(r))$ and $n>\overline{n}(r)$.}
From now on, the constant $C>0$ could depend on $\alpha$, but not on $\xi,r,n,t$.
Recalling that $\Psi_\xi'$ and $\lambda^\prime_\alpha$
are locally Lipschitz, we can solve the Cauchy problem for the system in the unknowns $ \beta_i$ and find a unique maximal solution $\beta^{(n)}\in C^1([0, t_n ), \R^n)$. Then, from \eqref{mu}
we deduce $\delta^{(n)}\in C^1([0, t_n ), \R^n)$.

We can now obtain some uniform (with respect to $n$) bounds in order to guarantee that $t_n=+\infty$. Let us take $v=\mu^n_{r,\xi}\in \mathbb{Z}_n$. This gives
\begin{align}
\frac{d}{dt} \left(\frac{1}{2}\Vert\nabla\phi^n_{r,\xi}\Vert^2+\Psi_\xi(\phi^n_{r,\xi})\right)+(\textbf{v}\cdot\nabla\phi^n_{r,\xi},\mu^n_{r,\xi})
+\Vert\nabla\mu^n_{r,\xi}\Vert^2+(\lambda^\prime_\alpha(\phi^n_{r,\xi})f,\partial_t\phi^n_{r,\xi})=0.
\label{gron2}
\end{align}
Moreover, by comparison, we have that
\begin{align}
\Vert\partial_t\phi^n_{r,\xi}\Vert_{V^\prime}^2\leq C\left(\Vert\nabla\mu^n_{r,\xi}\Vert^2+\Vert\nabla\phi^n_{r,\xi}\Vert^2\Vert\textbf{v}\Vert_{[L^4(\Omega)]^2}^2\right)
\label{Vprime}
\end{align}
Hence, for $\omega=	\frac{1}{4C}$, we can write
\begin{align}
\nonumber&	\frac{d}{dt} \left(\frac{1}{2}\Vert\nabla\phi^n_{r,\xi}\Vert^2+\Psi_\xi(\phi^n_{r,\xi})\right)+(\textbf{v}\cdot\nabla\phi^n_{r,\xi},\mu^n_{r,\xi})+\frac{3}{4}\Vert\nabla\mu^n_{r,\xi}\Vert^2\\&+(\lambda^\prime_\alpha(\phi^n_{r,\xi})f,\partial_t\phi^n_{r,\xi})+\omega	\Vert\partial_t\phi^n_{r,\xi}\Vert_{V^\prime}^2\leq C\Vert\nabla\phi^n_{r,\xi}\Vert^2\Vert\textbf{v}\Vert_{[L^4(\Omega)]^2}^2.
\label{gron}
\end{align}
Note that, $v\equiv 1\in \mathbb{Z}_n$ implies $\overline{\phi^n_{r,\xi}}\equiv\overline{\phi}_0$. Thus, by Young's and Poincar\'{e}'s inequalities, we get
$$
\vert(\textbf{v}\cdot\nabla\phi^n_{r,\xi},\mu^n_{r,\xi})\vert\leq \Vert\textbf{v}\Vert_{[L^\infty(\Omega)]^2}\Vert\nabla\mu^n_{r,\xi}\Vert\Vert\phi^n_{r,\xi}\Vert\leq C\Vert\textbf{v}\Vert_{[L^\infty(\Omega)]^2}^2(\Vert\nabla\phi^n_{r,\xi}\Vert^2+1)+\frac{1}{4}\Vert\nabla\mu^n_{r,\xi}\Vert^2,
$$
and, using Young's inequality, being $\lambda^\prime_\alpha,\lambda^{\prime\prime}_\alpha$ bounded, we also have
\begin{align*}
&(\lambda^\prime_\alpha(\phi^n_{r,\xi})f,\partial_t\phi^n_{r,\xi})\leq C\Vert \lambda^\prime_\alpha(\phi^n_{r,\xi})f\Vert^2_{V}+\frac{\omega}{2}\Vert\partial_t\phi^n_{r,\xi}\Vert^2_{V^\prime}\\&\leq C\left(\Vert\lambda^\prime_\alpha(\phi^n_{r,\xi})f\Vert^2+\Vert\lambda^{\prime\prime}_\alpha(\phi^n_{r,\xi})\nabla\phi^n_{r,\xi}f\Vert^2+\Vert\lambda^\prime_\alpha(\phi^n_{r,\xi})\nabla f\Vert^2\right)+\frac{\omega}{2}\Vert\partial_t\phi^n_{r,\xi}\Vert^2_{V^\prime}\\&\leq C\left(\Vert f\Vert^2+\Vert f\Vert_{L^\infty(\Omega)}^2\Vert\nabla\phi^n_{r,\xi}\Vert^2+\Vert\nabla f\Vert^2\right)+\frac{\omega}{2}\Vert\partial_t\phi^n_{r,\xi}\Vert^2_{V^\prime}
\end{align*}

Therefore, we can apply Gronwall's Lemma and the conservation of mass to obtain, independently of $n,r,\xi$,
\begin{align}
\Vert\phi^n_{r,\xi}\Vert_{L^\infty(0,T;V)}+\Vert\partial_t\phi^n_{r,\xi}\Vert_{L^2(0,T;V^\prime)}+\Vert\nabla\mu^n_{r,\xi}\Vert_{L^2(0,T;H)}\leq C,
\label{es1}
\end{align}
so that the approximating solution is defined in $[0,+\infty)$.

Let us now find a bound for $\overline{\mu}^n_{r,\xi}$. Multiplying \eqref{mu} by $\phi^n_{r,\xi}-\overline{\phi}^n_{r,\xi}$ and integrate over $\Omega$, we find
\begin{align*}
(\mu^n_{r,\xi},\phi^n_{r,\xi}-\overline{\phi}^n_{r,\xi})&=\alpha\Vert   \nabla\phi^n_{r,\xi}\Vert   ^2+(F_\xi'(\phi^n_{r,\xi}),\phi^n_{r,\xi}-\overline{\phi}^n_{r,\xi})-\alpha_0(\phi^n_{r,\xi},\phi^n_{r,\xi}-\overline{\phi}^n_{r,\xi})\\&+(\lambda_\alpha^\prime(\phi^n_{r,\xi})f,\phi^n_{r,\xi}-\overline{\phi}^n_{r,\xi}).
\end{align*}
Applying standard inequalities, we deduce (see \eqref{es1})
\begin{align*}
(F_\xi'(\phi^n_{r,\xi}),\phi^n_{r,\xi}-\overline{\phi}^n_{r,\xi})&\leq C_0^2\left(\Vert   \nabla\mu^n_{r,\xi}\Vert   \ \Vert   \nabla\phi^n_{r,\xi}\Vert  \right )-\alpha\Vert   \nabla\phi^n_{r,\xi}\Vert   ^2+2\alpha_0\left(\Vert   \phi^n_{r,\xi}\Vert   ^2+\Vert   \phi^n_{r,\xi}-\overline{\phi}^n_{r,\xi}\Vert   ^2\right)\\&+C\Vert f\Vert\leq C (1+\Vert   \nabla\mu^n_{r,\xi}\Vert\Vert   \phi^n_{r,\xi}-\overline{\phi}^n_{r,\xi}\Vert+\Vert f\Vert).
\end{align*}
It is well known (see, e.g., \cite{frigeri}) that there exists $C_0=C(\overline{\phi}_0)>0$ such that, for $0<\xi\leq\overline{\xi}$,
\begin{equation}
\int_\Omega \vert F_\xi'(\phi^n_{r,\xi}) \vert dx\leq C_0 \left( \left\vert \int_\Omega F_\xi'(\phi^n_{r,\xi})(\phi^n_{r,\xi}-\overline{\phi}_0)dx \right\vert +1\right),
\label{fprimo}
\end{equation}
provided that $\overline{\phi}_0\in(-1,1)$. Therefore, by \eqref{es1} and the previous estimate, we obtain
\begin{align*}
\vert \overline{\mu}^n_{r,\xi} \vert &\leq \frac{1}{ \vert \Omega \vert }\int_\Omega \vert \Psi_\xi'(\phi^n_{r,\xi}) \vert +C\Vert f\Vert_{L^1(\Omega)}\leq\frac{1}{ \vert \Omega \vert }\left(\int_\Omega \vert F_\xi'(\phi^n_{r,\xi}) \vert +\alpha_0\int_\Omega \vert \phi^n_{r,\xi} \vert \right)+C\Vert f\Vert_{L^1(\Omega)}\\&\leq\frac{C}{ \vert \Omega \vert }\left( \left\vert \int_\Omega F_\xi'(\phi^n_{r,\xi})(\phi^n_{r,\xi}-\overline{\phi}^n_{r,\xi}) \right\vert + 1 \right)+C\Vert f\Vert_{L^1(\Omega)}\\&\leq C\left(1+\Vert   \nabla\mu^n_{r,\xi}\Vert   \right)+C(\Vert f\Vert_{L^1(\Omega)}+\Vert f\Vert),
\end{align*}
which implies
\begin{align}
\Vert\mu^n_{r,\xi}\Vert_{L^2(0,T;V)}\leq C.
\label{mu2}
\end{align}
We can now multiply \eqref{mu} by $-\Delta\phi^n_{r,\xi}$ and integrate over $\Omega$, to get
\begin{align}
\nonumber
&\Vert\Delta\phi^n_{r,\xi}\Vert^2+(\Psi_\xi^{\prime\prime}(\phi^n_{r,\xi})\nabla\phi^n_{r,\xi},\nabla\phi^n_{r,\xi})\\
&=(\nabla\mu^n_{r,\xi},\nabla\phi^n_{r,\xi})+(\lambda^{\prime\prime}_\alpha(\phi^n_{r,\xi})\nabla\phi^n_{r,\xi}f+\lambda_\alpha^\prime(\phi^n_{r,\xi})\nabla f,\nabla\phi^n_{r,\xi}),
\label{h4}
\end{align}
which implies, thanks to Young's inequality, \eqref{es1}, the boundedness of $\lambda_\alpha^\prime,\ \lambda_\alpha^{\prime\prime}$ and the boundedness from below of $\Psi_\xi^{\prime\prime}$,
\begin{align}
\Vert\Delta\phi^n_{r,\xi}\Vert^2&\leq C\Vert\nabla\phi^n_{r,\xi}\Vert^2+\Vert\nabla\mu^n_{r,\xi}\Vert\Vert\nabla\phi^n_{r,\xi}\Vert+C\Vert\nabla\phi^n_{r,\xi}\Vert^2\Vert f\Vert_{L^\infty(\Omega)}+C\Vert\nabla f\Vert\Vert\nabla\phi^n_{r,\xi}\Vert\nonumber\\&\leq C(1+\Vert\nabla\mu^n_{r,\xi}\Vert+\Vert f\Vert_{L^\infty(\Omega)}+\Vert\nabla f\Vert),
\label{h5}
\end{align}
which implies (see \eqref{es1})
\begin{equation}
\Vert\phi^n_{r,\xi}\Vert_{L^4(0,T;H^2(\Omega))}\leq C,
\label{h2}
\end{equation}
for some $C>0$ that is independent of $n,r,\xi$.
Estimates \eqref{es1}, \eqref{mu2} and \eqref{h2} are enough to pass to the limit first as $n\rightarrow\infty$, then as $\xi\to0$ and eventually as $r\rightarrow\infty$,
by well known compactness arguments, obtaining thus a weak solution to \eqref{CH} with the desired regularity.
Notice that $\Vert\phi\Vert_{L^\infty(\Omega\times(0,T))}\leq 1$ can be proven by a standard method (see, e.g., \cite{Giorginitemam}).

We now discuss the stronger set of assumptions on the initial data and on $f,\textbf{v}$. Let us consider again the Galerkin setting described above.

We test \eqref{eqch} with $v=\partial_t\mu^n_{r,\xi}\in \mathbb{Z}_n$, This yields
\begin{align*}
\frac{1}{2}\frac{d}{dt}\Vert\nabla\mu^n_{r,\xi}\Vert^2+(\partial_t\mu^n_{r,\xi},\partial_t\phi^n_{r,\xi})+(\textbf{v}\cdot\nabla\phi^n_{r,\xi},\partial_t\mu^n_{r,\xi})=0.
\end{align*}
Observe that $(\textbf{v}\cdot \nabla\phi^n_{r,\xi},\mu^n_{r,\xi})\in AC([0,T]).$ Thus we can write
$$
(\textbf{v}\cdot\nabla\phi^n_{r,\xi},\partial_t\mu^n_{r,\xi})=\frac{d}{dt}(\textbf{v}\cdot \nabla\phi^n_{r,\xi},\mu^n_{r,\xi})-(\partial_t\textbf{v}\cdot\nabla\phi^n_{r,\xi},\mu^n_{r,\xi})-(\textbf{v}\cdot\nabla\partial_t\phi^n_{r,\xi},\mu^n_{r,\xi}),
$$
and
\begin{align*}
(\partial_t\mu^n_{r,\xi},\partial_t\phi^n_{r,\xi})&=\Vert\nabla\partial_t\phi^n_{r,\xi}\Vert^2+(\Psi_\xi^{\prime\prime}(\phi^n_{r,\xi})\partial_t\phi^n_{r,\xi},\partial_t\phi^n_{r,\xi})\\&+(\partial_t f\lambda_\alpha^\prime(\phi^n_{r,\xi}),\partial_t\phi^n_{r,\xi})+(f\lambda_\alpha^{\prime\prime}(\phi^n_{r,\xi})\partial_t\phi^n_{r,\xi},\partial_t\phi^n_{r,\xi}).
\end{align*}
By Sobolev embeddings and Young's inequality, we have
\begin{align*}
\vert(\partial_t\textbf{v}\cdot\nabla\phi^n_{r,\xi},\mu^n_{r,\xi})\vert&\leq \Vert\nabla\mu^n_{r,\xi}\Vert\Vert\partial_t\textbf{v}\Vert_{[L^4(\Omega)]^2}\Vert\nabla\phi^n_{r,\xi}\Vert_{[L^4(\Omega)]^2}\\&\leq C\Vert\nabla\mu^n_{r,\xi}\Vert^2\Vert\partial_t\textbf{v}\Vert_{[L^4(\Omega)]^2}^2+C\Vert\phi^n_{r,\xi}\Vert_{H^2(\Omega)}^2.
\end{align*}
Similarly, by Poincar\'{e}'s inequality, we deduce
$$
\vert(\textbf{v}\cdot\nabla\partial_t\phi^n_{r,\xi},\mu^n_{r,\xi})\vert\leq \Vert\nabla\mu^n_{r,\xi}\Vert\Vert\textbf{v}\Vert_{[L^\infty(\Omega)]^2}\Vert\partial_t\phi^n_{r,\xi}\Vert\leq C\Vert\nabla\mu^n_{r,\xi}\Vert^2\Vert\textbf{v}\Vert_{[L^\infty(\Omega)]^2}^2+\frac{1}{8}\Vert\nabla\partial_t\phi^n_{r,\xi}\Vert^2.
$$
Note that
$$
(\Psi_\xi^{\prime\prime}(\phi^n_{r,\xi})\partial_t\phi^n_{r,\xi},\partial_t\phi^n_{r,\xi})\geq -\alpha\Vert\partial_t\phi^n_{r,\xi}\Vert^2,
$$
but, by interpolation, we get
$$
\alpha\Vert\partial_t\phi^n_{r,\xi}\Vert^2\leq \frac{1}{8}\Vert\partial_t\nabla\phi^n_{r,\xi}\Vert^2+\frac{\alpha^2}{2}\Vert\partial_t\phi^n_{r,\xi}\Vert_{V^\prime}.
$$
Therefore, from \eqref{Vprime} and \eqref{h2}, we infer
\begin{align*}
&(\Psi_\xi^{\prime\prime}(\phi^n_{r,\xi})\partial_t\phi^n_{r,\xi},\partial_t\phi^n_{r,\xi})\geq -\frac{1}{8}\Vert\partial_t\nabla\phi^n_{r,\xi}\Vert^2-\frac{\alpha^2}{2}\Vert\partial_t\phi^n_{r,\xi}\Vert_{V^\prime}\\&\geq -\frac{1}{8}\Vert\partial_t\nabla\phi^n_{r,\xi}\Vert^2-C\left(\Vert\nabla\mu^n_{r,\xi}\Vert^2+\Vert\textbf{v}\Vert_{[L^4(\Omega)]^2}^2\right).
\end{align*}
Poincar\'{e}'s inequality gives
\begin{align*}
&(\partial_t f\lambda_\alpha^\prime(\phi^n_{r,\xi}),\partial_t\phi^n_{r,\xi})\leq C\Vert\partial_t f\Vert\Vert\nabla\partial_t\phi^n_{r,\xi}\Vert\leq \frac{1}{8}\Vert\nabla\partial_t\phi^n_{r,\xi}\Vert^2 +C\Vert\partial_t f\Vert^2,
\end{align*}
whereas, using also interpolation, we find
\begin{align*}
&(f\lambda_\alpha^{\prime\prime}(\phi^n_{r,\xi})\partial_t\phi^n_{r,\xi},\partial_t\phi^n_{r,\xi})\\
&\leq C\Vert\partial_t\phi^n_{r,\xi}\Vert^2\Vert f\Vert_{L^\infty(\Omega)}\leq C\Vert\partial_t\phi^n_{r,\xi}\Vert_{V^\prime}\Vert\nabla\partial_t\phi^n_{r,\xi}\Vert\Vert f\Vert_{L^\infty(\Omega)}\\&\leq \frac{1}{8}\Vert\nabla\partial_t\phi^n_{r,\xi}\Vert^2+C\Vert\partial_t\phi^n_{r,\xi}\Vert_{V^\prime}^2\Vert f\Vert_{L^\infty(\Omega)}^2\\&\leq \frac{1}{8}\Vert\nabla\partial_t\phi^n_{r,\xi}\Vert^2+C\left(\Vert\nabla\mu^n_{r,\xi}\Vert^2+\Vert\textbf{v}\Vert_{[L^4(\Omega)]^2}^2\right)\Vert f\Vert_{L^\infty(\Omega)}^2.
\end{align*}
Collecting the above inequalities, we deduce the following
\begin{align}
\nonumber &\frac{d}{dt}\left(\frac{1}{2}\Vert\nabla\mu^n_{r,\xi}\Vert^2+(\textbf{v}\cdot \nabla\phi^n_{r,\xi},\mu^n_{r,\xi})\right)+\frac{1}{2}\Vert\nabla\partial_t\phi^n_{r,\xi}\Vert^2\\&\leq C\Vert\nabla\mu^n_{r,\xi}\Vert^2\left(1+\Vert\partial_t\textbf{v}\Vert_{[L^4(\Omega)]^2}^2+\Vert\textbf{v}\Vert_{[L^\infty(\Omega)]^2}^2+\Vert f\Vert_{L^\infty(\Omega)}^2\right)\nonumber\\&+C\left(\Vert\phi^n_{r,\xi}\Vert_{H^2(\Omega)}^2+\Vert\textbf{v}\Vert_{[L^4(\Omega)]^2}^2+\Vert\partial_t f\Vert^2+\Vert\textbf{v}\Vert_{[L^4(\Omega)]^2}^2\Vert f\Vert_{L^\infty(\Omega)}^2\right).
\label{gronw}
\end{align}
Let us set
$$
\Lambda^n_{r,\xi}:=\frac{1}{2}\Vert\nabla\mu^n_{r,\xi}\Vert^2+(\textbf{v}\cdot \nabla\phi^n_{r,\xi},\mu^n_{r,\xi})
$$
and observe that, by \eqref{h2} and being $\Vert\textbf{v}\Vert_{L^\infty(0,T;L^4(\Omega)]^2)}\leq {C}$, we have
$$
\vert(\textbf{v}\cdot \nabla\phi^n_{r,\xi},\mu^n_{r,\xi})\vert\leq \Vert\textbf{v}\Vert_{[L^4(\Omega)]^2}\Vert\phi^n_{r,\xi}\Vert_{L^4(\Omega)}\Vert\nabla\mu^n_{r,\xi}\Vert\leq \tilde{C}+\frac{1}{4}\Vert\nabla\mu^n_{r,\xi}\Vert^2,
$$
for some $\tilde{C}>0$ independent of $n$. Hence we infer
\begin{align}
\Lambda^n_{r,\xi}\geq -\tilde{C}+\frac{1}{4}\Vert\nabla\mu^n_{r,\xi}\Vert^2.
\end{align}
Therefore, we can rewrite \eqref{gronw} as follows
\begin{align}
\nonumber &\frac{d}{dt}\left(\Lambda^n_{r,\xi}+\tilde{C}\right)+\frac{1}{2}\Vert\nabla\partial_t\phi^n_{r,\xi}\Vert^2\\&\leq C\left(\Lambda^n_{r,\xi}+\tilde{C}\right)\left(1+\Vert\partial_t\textbf{v}\Vert_{[L^4(\Omega)]^2}^2+\Vert\textbf{v}\Vert_{[L^\infty(\Omega)]^2}^2+\Vert f\Vert_{L^\infty(\Omega)}^2\right)\nonumber\\&+C\left(\Vert\phi^n_{r,\xi}\Vert_{H^2(\Omega)}^2+\Vert\textbf{v}\Vert_{[L^4(\Omega)]^2}^2+\Vert\partial_t f\Vert^2+\Vert\textbf{v}\Vert_{[L^4(\Omega)]^2}^2\Vert f\Vert_{L^\infty(\Omega)}^2\right),
\label{lambdan}
\end{align}
which, taking the assumptions on $\textbf{v}$ and $f$ into account, implies, by Gronwall's Lemma
\begin{align}
0\leq \Lambda^n_{r,\xi}(t)+\tilde{C}\leq  C(T)(\Lambda^n_{r,\xi}(0)+\tilde{C})+C(T), \quad \forall\, t\in (0,T).
\label{mun}
\end{align}
We are left to prove a uniform bound on $\Lambda^n_{r,\xi}(0)$. To this aim, observe that, by Sobolev embeddings, being $\textbf{v}\in C([0,T];\textbf{V}_\sigma)$ and $\textbf{v}(0)=\textbf{v}_0$,
we have
\begin{align*}
\Lambda^n_{r,\xi}(0)&=(\textbf{v}_0\cdot \nabla\phi^n_{r,\xi}(0),\mu^n_{r,\xi}(0))+\frac{1}{4}\Vert\nabla\mu^n_{r,\xi}(0)\Vert^2\\
&\leq \Vert\textbf{v}_0\Vert_{\textbf{V}_\sigma}\Vert\phi_0\Vert_{V}\Vert\nabla\mu^n_{r,\xi}(0)\Vert+\frac{1}{4}\Vert\nabla\mu^n_{r,\xi}(0)\Vert^2.
\end{align*}
On the other hand, notice that, for any ${w}\in \mathbb{Z}_n$,
\begin{align*}
(\nabla\mu^n_{r,\xi}(t),\nabla w)&=(-\nabla\Delta\phi^n_{r,\xi}(t),\nabla w)+(\nabla\tilde{P}_n(\Psi_\xi^\prime(\phi^n_{r,\xi}(t))),\nabla w)\\
&+(\nabla(\lambda_\alpha^\prime(\phi^n_{r,\xi}(t))f(t)),\nabla w),
\end{align*}
for any $t\in[0,T]$. Therefore, setting $t=0$ and $w= \nabla\mu^n_{r,\xi}(0)\in \mathbb{Z}_n$, Young's inequality entails
\begin{align}
\Vert\nabla\mu^n_{r,\xi}(0)\Vert^2\leq C\left(\Vert\nabla\Delta\tilde{P}_n\phi_0+\nabla\Psi^\prime_\alpha(\tilde{P}_n\phi_0)\Vert^2+\Vert\nabla(\lambda_\alpha^\prime(\tilde{P}_n\phi_0)f_0)\Vert^2\right).
\label{mu1}
\end{align}
On account of the boundedness of $\lambda_\alpha^\prime,\lambda_\alpha^{\prime\prime}$ and the properties of $\tilde{P}_n$, we get
\begin{align*}
\Vert\nabla(f_0\lambda_\alpha^\prime(\tilde{P}_n\phi_0))\Vert&\leq\Vert\nabla f_0\lambda_\alpha^\prime(\tilde{P}_n\phi_0)\Vert+\Vert f_0\lambda_\alpha^{\prime\prime}(\tilde{P}_n\phi_0)\nabla\tilde{P}_n\phi_0\Vert\\&\leq C\Vert f_0\Vert_{V}+C\Vert f_0\Vert_{L^4(\Omega)}\Vert\nabla\tilde{P}_n\phi_0\Vert_{[L^4(\Omega)]^2}\\
&\leq  C\Vert f_0\Vert_{V}+C\Vert f_0\Vert_{L^4(\Omega)}(\Vert\tilde{P}_n\phi_0\Vert+\Vert\Delta\tilde{P}_n\phi_0\Vert)\\&\leq C\Vert f_0\Vert_{V}+C\Vert f_0\Vert_{L^4(\Omega)}(\Vert\phi_0\Vert+\Vert\Delta\phi_0\Vert)\leq C,
\end{align*}
since $H^2(\Omega)\hookrightarrow W^{1,4}(\Omega)$. Concerning the first term on the right-hand side of \eqref{mu1}, we can argue exactly as in \cite[(4.38)-(4.39)]{Giorginitemam} to infer that, for any fixed $r>r_0$, $\lambda\in(0,\lambda^*(r_0))$ and $n>\overline{n}(r_0,\lambda^*(r_0))$:
\begin{equation*}
\tilde{\Lambda}^n_{r,\xi}(0)\leq C,
\end{equation*}
independently of $n,r,\xi$.
From \eqref{lambdan} and \eqref{mun} we thus get, independently of $n,r,\xi$,
\begin{align}
\Vert\mu^n_{r,\xi}\Vert_{L^\infty(0,T;V)}+\Vert\partial_t\phi^n_{r,\xi}\Vert_{L^2(0,T;V)}\leq C,
\label{mu4}
\end{align}
and in conclusion we can pass to the limit first in $n\rightarrow\infty$, then in $\xi\to0$ and then in $r\rightarrow\infty$, by compactness arguments, to obtain a solution $(\phi,\mu)$ to \eqref{CH} with the desired regularity. Concerning the additional regularity, by \cite[Thm A.2]{Giorginitemam}, we deduce that, for any $q\in[2,\infty)$,
\begin{align*}
\Vert F^\prime(\phi)\Vert_{L^q(\Omega)}\leq C(\Vert\mu\Vert_{L^q(\Omega)}+\Vert\phi\Vert_{L^q(\Omega)}+\Vert f\Vert_{L^q(\Omega)}),
\end{align*}
but being $\mu\in L^\infty(0,T;V)$, $\phi\in L^4(0,T;H^2(\Omega))$ and $f\in L^\infty(0,T;V)$, we infer that
\begin{align*}
\Vert F^\prime(\phi)\Vert_{L^\infty(0,T;L^q(\Omega))}\leq C,
\end{align*}
from which, by elliptic regularity, we deduce that
\begin{align*}
\Vert \phi\Vert_{L^\infty(0,T;W^{2,q}(\Omega))}\leq C.
\end{align*}
Then, again by \cite[Thm A.2]{Giorginitemam} and by \eqref{mu4}, we infer, for any $q\in[2,\infty)$, that
\begin{align}
\Vert F^{\prime\prime}(\phi) \Vert_{L^\infty(0,T;L^q(\Omega))}\leq C,
\label{Psisec}
\end{align}
and this ends the proof.

\textbf{Acknowledgments.} The authors thank Andrea Giorgini for several discussions and helpful comments. The authors have been partially funded by MIUR-PRIN research grant n. 2020F3NCPX and they are members of Gruppo Nazionale per l'Analisi Ma\-te\-ma\-ti\-ca, la Probabilit\`{a} e le loro Applicazioni (GNAMPA), Istituto Nazionale di Alta Matematica (INdAM).


\begin{thebibliography}{widestlabel}
	
	

        \bibitem{Abels} Abels, H. (2009), \textit{On a diffuse interface model for two-phase flows of viscous, incompressible fluids
			with matched densities}, Arch. Ration. Mech. Anal. \textbf{194}, 463-506.
		
		\bibitem{ADG} Abels, H., Depner, D., Garcke, H. (2013), \textit{Existence of weak solutions for a diffuse interface model for two-phase
			flows of incompressible fluids with different densities}, J. Math. Fluid Mech. \textbf{15}, 453-480.
		
		\bibitem{AF} Abels, H., Feireisl, E. (2008), \textit{On a diffuse interface model for a two-phase flow of compressible viscous fluids}, Indiana Univ. Math. J. \textbf{57}, 659-698.
	
		
		\bibitem{xu} Alber, M., Kim, O.V., Litvinov, R.I., Weisel, J.W., Xu, S., Xu, Z. (2017),\textit{ Model predictions of deformation, embolization
		and permeability of partially obstructive blood clots under variable shear flow},  J. R. Soc. Interface \textbf{14}, pp.13.
		
		\bibitem{tierra} Alber, M.S., Nerenberg, R., Pavissich, J.P., Tierra, G., Xu., Z. (2015), \textit{Multicomponent model of deformation and
		detachment of a biofilm under fluid flow},  J. R. Soc. Interface \textbf{12}, pp.13.
		
		\bibitem{Anand} Anand, M., Rajagopal, K.R. (2004), \textit{A shear-thinning viscoelastic fluid model for describing the flow of blood}, Int. J. Cardiovascular Medicine and Science \textbf{4},  59-68.
		
		
		\bibitem{AMW} Anderson, D.M., McFadden, G.B., Wheeler, A.A. (1998), \textit{Diffuse-interface methods in fluid mechanics},
		Annu. Rev. Fluid Mech. \textbf{30}, 139-165.
	
		
		\bibitem{Boyerlibro} Boyer, F.,  Fabrie, P., \textit{Mathematical Tools for the Study of the Incompressible Navier-Stokes Equations and Related Models}, Springer-Verlag, New York, 2013.
		

		\bibitem{brezisbourg} Bourguignon, J.P., Brezis, H. (1974), \textit{Remarks on the Euler Equation}, J. Funct. Anal. \textbf{15}, 341-363.
		
		
	    \bibitem{BTP} Brangwynne, C.P., Tompa, P., Pappu, R.V. (2015), \textit{Polymer physics of intracellular phase transitions}, Nat. Phys. \textbf{11}, 899-904.
		
		\bibitem{Brezis} Br\'{e}zis, H., \textit{Op\'{e}rateurs Maximaux Monotones et Semi-groupes de Contractions dans les
			Espaces de Hilbert}, North-Holland Mathematics Studies vol. 5, Amsterdam, 1973.

        \bibitem{Cag} Caginalp, G. (1986), \textit{An analysis of a phase field model of a free boundary}, Arch. Rational Mech. Anal. \textbf{92}, 205--245.
		
		\bibitem{CH}Cahn, J.W., Hilliard, J.E., \textit{Free energy of a nonuniform system I. Interfacial free energy} (1958), J. Chem. Phys. \textbf{28}, 258-267.
		
		
		\bibitem{Chapin} Chapin, J.C., Hajjar, K.A. (2015),\textit{ Fibrinolysis and the control of blood coagulation}, Blood Rev. \textbf{29}, 17-24.
		
		\bibitem{Colman}  Colman R.W., Marder V.J., Clowes A.W., George J.N., Goldhaber S.Z., eds.,
        \textit{Hemostasis and thrombosis: basic principles and clinical practice}, Lippincott Williams \& Wilkins, Philadelphia, 2006.
		
	
		\bibitem{contigiorgini} Conti, M., Giorgini, A. (2020), \textit{Well-posedness for the Brinkman-Cahn-Hilliard system with unmatched viscosities}, J. Differential Equations \textbf{268},
        6350-6384.
		
	
		\bibitem{Dong} Dong, S., Shen J. (2012), \textit{A time-stepping scheme involving constant coefficient matrices for phase-field simula-
			tions of two-phase incompressible flows with large density ratios}, J. Comput. Phys. \textbf{231}: 5788-5804.
		
		\bibitem{em}Emelianov, S., Chen, X., O’donnell, M., Knipp, B., Myers, D., Wakefield, T. et al. (2002), \textit{Triplex ultrasound: elasticity
		imaging to age deep venous thrombosis}, Ultrasound Med. Biol. \textbf{28}, 757-767.

         \bibitem{E} Emmerich, H., \textit{The diffuse interface approach in materials science,} Springer, Berlin Heidelberg, 2003.
		
		  \bibitem{Fedosov} Fedosov, D.A., \textit{Hemostasis is a highly multiscale process: Comment on "Modeling thrombosis in silico:
			Frontiers, challenges, unresolved problems and milestones" by A.V. Belyaev et al.}, Phys. Life Rev. \textbf{26-27}: 108-109.
		
		
		\bibitem{Fineschi} Fineschi V, Neri, M., Pomara, C., Riezzo, I., Turillazzi, E. (2009), \textit{Histological age determination of venous thrombosis: a neglected forensic task in fatal pulmonary thrombo-embolism}, Forensic Sci. Int. \textbf{186}, 22-28.
		
		\bibitem{Formaggia} Formaggia, L.,  Quarteroni, A., Veneziani, A., \textit{Cardiovascular Mathematics: Modeling and Simulation of the Circulatory System}, Springer-Verlag, Milano, 2009.
		
		
		\bibitem{frigeri} Frigeri, S., Grasselli, M. (2012), \textit{Nonlocal Cahn-Hilliard-Navier-Stokes systems with singular
			potential}, Dyn. Partial Differ. Equ. \textbf{24}, 827-856.
		
		\bibitem{frigeri2}  Frigeri, S., Grasselli, M. (2012), \textit{Global and trajectories attractors for a nonlocal Cahn-Hilliard-Navier-Stokes system}, J. Dynam. Differential
		Equations \textbf{24}, 827-856.
		
		\bibitem{gal2} Gal, C.G., Grasselli, M. (2010), \textit{Asymptotic behavior of a Cahn-Hilliard-Navier-Stokes in 2D}, Ann.
		Inst. H. Poincar\'e Anal. Non Lin\'eaire \textbf{27}, 401-436.
		
		\bibitem{Giorgini} Giorgini, A. (2021), \textit{Well-posedness of the two-dimensional Abels-Garcke-Grün model for two-phase flows with unmatched densities}, Calc. Var. Partial Differential Equations \text{60}, 100,
        Paper No., 40 pp.
		
		\bibitem{GGW} Giorgini, A., Grasselli, M., Wu, H. (2022), \textit{Diffuse interface models for incompressible binary fluids and the mass-conserving Allen-Cahn approximation},  J. Funct. Anal. 283 (2022), Paper No.109631, 86 pp.
		
		
		\bibitem{Giorginitemam} Giorgini, A., Miranville, A., Temam, R. (2019), \textit{Uniqueness and regularity for the Navier-Stokes-Cahn-Hilliard system}, SIAM J. Math. Anal. \textbf{51}, 2535-2574.
		
		\bibitem{GP} Grasselli, M., Poiatti, A. (2022), \textit{The Cahn-Hilliard-Boussinesq system with singular
		potential}, Commun. Math. Sci. \textbf{20}, 897-946
		
	    \bibitem{salgado}	Guermond, J.-L., Salgado, A. (2011), \textit{A note on the Stokes operator and its powers}, J. Appl. Math. Comput. \textbf{36}, 241-250.
		
		
		\bibitem{GPV} Gurtin, M.E. Polignone D., Vi\~{n}als, J. (1996), \textit{Two-phase binary fluids and immiscible fluids described by an order parameter}, Math. Models Methods Appl. Sci. \textbf{6}, 815-831.
		
		
        \bibitem{HeWu} He, J., Wu, H. (2021), \textit{Global well-posedness of a Navier–Stokes–Cahn–Hilliard system with chemotaxis and singular potential in 2D}, J. Differential Equations \textbf{297}, 47-80.

		
		\bibitem{Saal} Hieber, M., Saal, J., \textit{The Stokes Equation in the $L^p$ -Setting: Well-Posedness and Regularity Properties}, Handbook of mathematical analysis in mechanics of viscous fluids, 117-206, Springer, Cham, 2018.
		
		\bibitem{Lee} Humphrey, J., Lee, A., Lee, Y.U., Rausch, M., \textit{Histological and biomechanical changes in a mouse model of
			venous thrombus remodeling} (2015), Biorheology \textbf{52}: 235-245.
		
		
	    \bibitem{Karniadakis} Karniadakis, G.E., Humphrey, J.D., Li, H., Yazdani, A., Zheng, X. (2020):
	    \textit{A three-dimensional phase-field model for multiscale modeling of thrombus biomechanics in blood vessels}, PLoS Comput. Biol. \textbf{16}, pp.24.

	    \bibitem{Zeng} Karniadakis, G., Zheng, X. (2016), \textit{A phase-field/ALE method for simulating fluid-structure interactions in two-
		phase flow}, Comput. Methods Appl. Mech. Engrg. \textbf{309}, 19-40.
	

        \bibitem{KTT} Kim, W., Tawry, K., Temam, R. (2022), \textit{Local well-posedness of a three-dimensional phase-field model for thrombus and blood flow},	
         Rev. R. Acad. Cienc. Exactas F\'{\i}s. Nat. Ser. A Mat. RACSAM, \textbf{116}, 149, 23 pp.
		
        \bibitem{Lin}Lin, F.H., Liu, C., Zhang, P. (2005), \textit{On hydrodynamics of viscoelastic fluids}, Comm. Pure Appl. Math. \textbf{5}, 1437-1471.

        \bibitem{Lin2} Lin, F.H., Zhang, P. (2008), \textit{On the initial-boundary value problem of the incompressible viscoelastic fluid system}, Comm.
        Pure Appl. Math. \textbf{61}, 539-558.

        \bibitem{Lady} Ladyzhenskaja, O.A., Solonnikov, V.A. (1975), \textit{The unique solvability of an initial-boundary value problem for viscous incompressible inhomogeneous fluids, boundary value problems of mathematical physics, and related questions of the theory of functions}, 8. Zap. Naučn. Sem. Leningrad. Otdel. Mat. Inst. Steklov. (LOMI) \textbf{52}, 52-109.

        \bibitem{LM} Lions, J.-L., Magenes, E.: \textit{Problèmes aux Limites non Homogènes et Applications}, vol. 1. Dunod, Paris, 1968.

        \bibitem{Lions} Lions, J.-L., Magenes, E., \textit{Non-homogeneous boundary value problems and applications}, Springer-Verlag, Berlin, 1972.

        \bibitem{Mack} Mackman, N., Tilley, R.E., Key, N.S. (2007), \textit{Role of the extrinsic pathway of blood coagulation in hemostasis and thrombosis},
        Arterioscler. Thromb. Vasc. Biol. \textbf{27}, 1687-1693.

	    \bibitem{Maxwell} Maxwell, D., \textit{Initial Data for Black Holes and Rough Spacetimes, PhD thesis, University of Washington}, 2004.

        \bibitem{Miranville} Miranville, A., Zelik, S. (2004), \textit{Robust exponential attractors for Cahn-Hilliard type equations with singular potentials}, Math.
        Methods. Appl. Sci. \textbf{2}7 , 545-582.	


        \bibitem{temam} Temam, R., \textit{Navier-Stokes equations}, Theory and Numerical Analysis, North-Holland, Amsterdam, 1984.

        \bibitem{temam2} Temam, R., \textit{Infinite-dimensional dynamical systems in mechanics and physics}, Springer-Verlag,
        New York, 1997.

        \bibitem{Xu} Xu, Z., Chen, N., Kamocka, M.M., Rosen, E.D., Alber M. (2008), \textit{A multiscale model of thrombus development}, J. R. Soc. Interface \textbf{5}: 705-722.


\end{thebibliography}
\end{document}